\documentclass[12pt,leqno]{amsart}
\usepackage{amssymb}
\usepackage{amsmath,amssymb}
\newtheorem{theorem}{Theorem}[section]
\newtheorem{corollary}{Corollary}[section]

\newtheorem{proposition}{Proposition}[section]

\begin{document}
\theoremstyle{plain}
\newtheorem{MainThm}{Theorem}
\newtheorem{thm}{Theorem}[section]
\newtheorem{clry}[thm]{Corollary}
\newtheorem{prop}[thm]{Proposition}
\newtheorem{lem}[thm]{Lemma}
\newtheorem{deft}[thm]{Definition}
\newtheorem{hyp}{Assumption}
\newtheorem*{ThmLeU}{Theorem (J.~Lee, G.~Uhlmann)}

\theoremstyle{definition}
\newtheorem{rem}[thm]{Remark}
\newtheorem*{acknow}{Acknowledgments}
\numberwithin{equation}{section}
\newcommand{\eps}{{\varphi}repsilon}
\renewcommand{\d}{\partial}
\newcommand{\re}{\mathop{\rm Re} }
\newcommand{\im}{\mathop{\rm Im}}
\newcommand{\R}{\mathbf{R}}
\newcommand{\C}{\mathbf{C}}
\newcommand{\N}{\mathbf{N}}
\newcommand{\D}{C^{\infty}_0}
\renewcommand{\O}{\mathcal{O}}
\newcommand{\dbar}{\overline{\d}}
\newcommand{\supp}{\mathop{\rm supp}}
\newcommand{\abs}[1]{\lvert #1 \rvert}
\newcommand{\csubset}{\Subset}
\newcommand{\detg}{\lvert g \rvert}
\newcommand{\ppp}{\partial}
\newcommand{\dd}{\mbox{div}\thinspace}

\title
[Navier-Stokes equations and Lam\'e system]
{Global uniqueness in inverse boundary value problems
for Navier-Stokes equations and Lam\'e system in
two dimensions}

\author{
O.~Yu.~Imanuvilov and \,
M.~Yamamoto }
\thanks{ Department of Mathematics, Colorado State
University, 101 Weber Building, Fort Collins, CO 80523-1874, U.S.A.
E-mail: {\tt oleg@math.colostate.edu}.  Partially supported by NSF grant DMS 1312900}\,

\thanks{ Department of Mathematical Sciences, University
of Tokyo, Komaba, Meguro, Tokyo 153, Japan e-mail:
myama@ms.u-tokyo.ac.jp}

\date{}

\maketitle

\begin{abstract}
We consider inverse boundary value problems for the Navier-Stokes equations
and the isotropic Lam\'e system in two dimensions.
The uniqueness without any smallness assumptions on unknown coefficients,
which is called global uniqueness, was longstanding open problems
for the Navier-Stokes equations and the isotropic Lam\'e system in
two dimensions.  We prove the global uniqueness
for both inverse boundary value problems.
Our methodology are common for both systems and the key is the construction
of complex geometric optics solutions
after decoupling the systems into weakly coupling systems.
\end{abstract}
\section{Introduction and main results}
In this article, we consider inverse boundary value problems for the
two dimensional Navier-Stokes equations and Lam\'e system where we determine
spatially varying viscosity and two Lam\'e coefficients in the
Navier-Stokes equations and the Lam\'e system respectively by the
Dirichlet-to-Neumann maps on the whole boundary.

Throughout this paper, let $\Omega \subset \Bbb R^2$ be a bounded
domain with smooth boundary $\partial\Omega$ and $\nu=(\nu_1,\nu_2)$ be
the outward unit normal vector to $\partial\Omega$.
We set $x = (x_1, x_2) \in \Bbb R^2$ and
$\beta = (\beta_1, \beta_2) \in (\Bbb {N}_+)^2$, and
$\vert \beta\vert = \beta_1 + \beta_2$,
$\partial_x^{\beta} = \frac{\partial^{\beta_1}}{\partial x_1^{\beta_1}}
\frac{\partial^{\beta_2}}{\partial x_2^{\beta_2}}$.
Let $(\cdot,\cdot)$ be the scalar product in $\Bbb R^2$.

We start with the formulation of the inverse boundary value problem for
the Navier-Stokes equations.  In the domain $\Omega$, we consider
the stationary Navier-Stokes equations:
$$
G_\mu(x,D)({\bf u},p) := \Biggl(\sum_{j=1}^2\left(-2\frac{\partial}{\partial
x_j}(\mu(x)\epsilon_{1j}({\bf u}))+u_j\frac{\partial u_1}{\partial
x_j}\right) + \frac{\partial p}{\partial x_1},
$$
$$
\sum_{j=1}^2\left(-2\frac{\partial}{\partial
x_j}(\mu(x)\epsilon_{2j}({\bf u}))+u_j\frac{\partial u_2}{\partial
x_j}\right) + \frac{\partial p}{\partial x_2}\Biggr) =0\quad\mbox{in}
\,\,\Omega,
$$
where ${\bf u}=(u_1,u_2)$ is a velocity field, $p$ is a pressure
and $\epsilon_{ij}({\bf u})
=\frac 12 (\frac{\partial u_i}{\partial
x_j}+\frac{\partial u_j}{\partial x_i}).$ Assume that

\begin{equation}\label{ox}
\mu(x)>0\quad\mbox{on}\,\,\overline{\Omega}, \quad
\mu\in C^m(\overline\Omega).
\end{equation}
Here and henceforth $m \in \Bbb N$ is sufficiently large
(e.g., $m=10$ is sufficient).

We define the Dirichlet-to-Neumann map $\Lambda_\mu$:
\begin{equation}\label{z1}
\Lambda_\mu({\bf f}) = \left(\frac{\partial {\bf u}}{\partial \nu},
p\right)\vert_{\partial\Omega},
\end{equation}
where
$$
G_\mu(x,D)({\bf u},p)=0\quad\mbox{in}\,\,\Omega,
\quad {\bf u} = {\bf f} \quad \mbox{on $\partial\Omega$},
\quad \mbox{div}\,{\bf u}=0, \quad  {\bf u}\in W_2^2(\Omega), \quad
p \in W_2^1(\Omega)
$$
and
$$
D(\Lambda_\mu)\subset \left\{{\bf f}\in W_2^\frac 32(\partial\Omega);
\thinspace \int_{\partial\Omega}(\nu, {\bf f}) d\sigma=0\right\}.
$$
The first subject of this paper is the following inverse boundary
value problem:
\\
{\it Using the  Dirichlet-to-Neumann map $\Lambda_{\mu}$, determine the
coefficient $\mu$.}

Our first main result is the global uniqueness:
\begin{theorem}\label{vokal}
We assume that $\mu_1,\mu_2\in C^{10}(\overline\Omega)$ and
$\partial_x^{\beta}\mu_1 = \partial_x^{\beta}\mu_2$
on $\partial\Omega$ for each multi-index $\beta$ with
$\vert \beta\vert \le 10$.
If $\Lambda_{\mu_1} = \Lambda_{\mu_2}$, then
$\mu_1=\mu_2$ in $\Omega$.
\end{theorem}

The uniqueness result of the above theorem holds without any assumption
on smallness unknown coefficients or closeness of these coefficients to
constants.  We call such uniqueness the global uniqueness.
To the authors' best knowledge, there are no global uniqueness results
for the Navier-Stokes equations in two dimensions.  On the other hand,
in the three dimensional case, the global uniqueness was proved in
Heck, Li and Wang \cite{Heck} for the Stokes equations and
in Li and Wang \cite{Li} for the Navier-Stokes equations.
\\

Next we consider the inverse boundary value problem for the two dimensional
Lam\'e system.  Assume that
$$
\mu(x)>0, \thinspace (\lambda+\mu)(x)>0\quad
\mbox{on}\,\,\overline\Omega.
$$
We set
$$
\mathcal L_{\mu,\lambda}(x,D) {\bf u} =
\left(\sum_{j=1}^2 \partial_{x_j}(\mu(\partial_{x_j}u_1
+ \partial_{x_1} u_j)),
\sum_{j=1}^2 \partial_{x_j}(\mu(\partial_{x_j}u_2
+ \partial_{x_2} u_j)) \right) + \nabla(\lambda\mbox{div}\thinspace{\bf u})
$$
$$
= \mu\Delta {\bf u} + (\mu+\lambda)\nabla\mbox{div}\thinspace{\bf u}
+ (\mbox{div}\thinspace {\bf u})\nabla\lambda
$$
$$
+ ((\nabla\mu, \nabla u_1), (\nabla\mu, \nabla u_2))
+ ((\nabla\mu, \partial_{x_1}{\bf u}), (\nabla\mu, \partial_{x_2}{\bf u})).
$$
Here ${\bf u}(x)=(u_1(x), u_2(x))$
describes displacement, and we call the functions $\lambda$ and
$\mu$ the Lam{\' e} coefficients.

We define the
Dirichlet-to-Neumann map $\Lambda_{\mu,\lambda}$ as follows:
\begin{equation}\label{z1XV}
\Lambda_{\mu,\lambda}f =
\left(\mu\partial_{\nu}u_1 + \mu(\partial_{x_1}{\bf u}, \nu)
+ \lambda(\mbox{div}\thinspace {\bf u})\nu_1, \thinspace
\mu\partial_{\nu}u_2 + \mu(\partial_{x_2}{\bf u}, \nu)
+ \lambda(\mbox{div}\thinspace {\bf u})\nu_2 \right)
\vert_{\partial\Omega},
\end{equation}
where
\begin{equation}\label{lammas}
\mathcal L_{\mu,\lambda}(x,D) {\bf u} = 0 \quad\mbox{in $\Omega$},
\quad {\bf u}\vert_{\partial\Omega} = f.
\end{equation}
The second subject is the uniqueness in determining $\lambda, \mu$ by
$\Lambda_{\mu,\lambda}$.

Then we can prove
\begin{theorem}\label{lame}
Assume that
$\lambda_1,\lambda_2, \mu_1, \mu_2\in C^{10}(\overline\Omega).$   Then $\Lambda_{\mu_1,\lambda_1}
= \Lambda_{\mu_2,\lambda_2}$ implies
that $\lambda_1 = \lambda_2$ and $\mu_1 =\mu_2$ in
$\Omega$.
\end{theorem}

The global uniqueness for the Lam\'e system has been an open
problem in spite of the physical significance, and Theorem 1.2
affirmatively solves for the two dimensional case.
On the other hand, we can refer to non-global uniqueness as follows.

{\bf Two dimensional case.}

In \cite{NU1} Nakamura and Uhlmann  proved that
the Dirichlet-to-Neumann map uniquely determines
the Lam\'e coefficients, assuming that $\lambda, \mu$ are sufficiently close to
positive constants. They assume that $\lambda, \mu \in W^{31}_\infty(\Omega)$,
and our regularity assumption is $\lambda, \mu \in C^{10}(\overline{\Omega})$,
although it may be more relaxed.

Imanuvilov and Yamamoto \cite{IYlame} proved the uniqueness
by the Dirichlet-to-Neumann map limited to an arbitrarily sub-boundary,
provided that $\mu_1, \mu_2$ are some constants.
These results assume some assumptions on smallness or constants of
unknown coefficients and the assumptions are quite restrictive.
To our best knowledge, Theorem 2 is the first result on the global
uniqueness in two dimensions.

{\bf Three dimensional case.}

We refer to
Eskin and Ralston \cite{Es1}, Imanuvilov, Uhlmann and Yamamoto
\cite{IUYlame}.
In \cite{Es1}, the uniqueness for $\lambda$ and $\mu$ is proved
provided that $\nabla \mu$ is small in suitable norm, while
in \cite{IUYlame}, the uniqueness corresponding to the uniqueness
\cite{IYlame} in two dimensions is proved.
In Nakamura and Uhlmann \cite{NU2}, the authors attempted
to prove the global  uniqueness
for determination of Lam\'e coefficients in $C^{\infty}(\overline{\Omega})$.
However the global uniqueness in the three dimensional
case remains a significant unsolved problem.

As for other types of inverse boundary value problems for
Lam\'e system, see Akamatsu, Nakamura and Steinberg \cite{ANS},
Ikehata \cite{IK} and Nakamura and Uhlmann \cite{NU3}.

The inverse boundary value problem was originally considered in
Calder\'on \cite{C} which considered the determination of
$\sigma$ in conductivity equation div $(\sigma(x)\nabla u(x)) = 0$ by the
Dirichlet-to-Neumann map.
In the case of the higher dimensional case, that is,
the spatial dimensions are equal to or greater than
$3$, the uniqueness for conductivity equation and Schr\"odinger
equation is proved in Sylverster and Uhlmann \cite{SU}, and in two dimensions,
see Nachman \cite{N}.
The inverse boundary value problems have attracted a lot of
attention.  Here we do not compose a complete list of references and
we refer to two surveys Uhlmann \cite{U} which is a survey as of 2009 and
Imanuvilov and Yamamoto \cite{IY5}.  The latter mostly presents results on
global uniqueness by the Dirichlet-to-Neumann map restricted on sub-boundary
which were published after 2009.

This paper is composed of seven sections.
In Section 2, we established properties of integral operator in $\Bbb C$
and oscillatory integral operator, and proved a basic Carleman estimate for
a first-order equation.
In Section 3 we consider second-order elliptic systems
whose principal parts are the two dimensional Laplacian and
construct special solutions (Propositions 3.5 and 3.6) which are used for
establishing adequate asymptotic behavior of some bilinear form.
Section 4 is devoted to proof of Proposition 4.6 which
yields integral-differential equalities
from asymptotic behavior of some quadratic forms generated by
a second-order elliptic operator and the special solutions considered in
Proposition 3.5.
In Section 5 we construct complex geometric optics solutions to the
Navier-Stokes equations by reducing the equations to a weakly coupling
elliptic system considered in Section 3 and in Section 6 we complete the proof
of Theorem 1.1.
Section 7 is the proof of Theorem 1.2 for Lam\'e system in a way
similar to one developed in  Sections 5 and 6.

\section{Preliminary results}\label{sec1}

Throughout the paper, we use the following notations.
\\

\noindent {\bf Notations.} $i=\sqrt{-1}$, $x_1, x_2, \xi_1, \xi_2
\in {\Bbb R}^1$, $z=x_1+ix_2$, $\zeta=\xi_1+i\xi_2$, $\overline{z}$
denotes the complex conjugate of $z \in \Bbb C$. We identify $x =
(x_1,x_2) \in {\Bbb R}^2$ with $z = x_1 +ix_2 \in {\Bbb C},$
$\partial_z = \frac 12(\partial_{x_1}-i\partial_{x_2})$,
$\partial_{\overline z}= \frac12(\partial_{x_1}+i\partial_{x_2}),$
$\beta=(\beta_1,\beta_2), \vert \beta\vert=\beta_1+\beta_2.$
$D=(\frac{1}{i}\frac{\partial}{\partial
x_1},\frac{1}{i}\frac{\partial}{\partial x_2})$,
$\partial{x_j} = \frac{\partial}{\partial x_j}$, $j=1,2$.
For $\widehat x\in \Bbb R^2$ and $r>0$, let $B(\widehat x, r)
= \{x \in \Bbb R^2\vert \thinspace \vert x-\widehat x\vert < r\}$
and $S(\widehat x, r)
= \{x \in \Bbb R^2\vert \thinspace \vert x-\widehat x\vert = r\}$.

Let $W^m_p(\Omega)$, $p\ge 1$ denote a usual Sobolev space and
let $\mathaccent'27{W}_p^m(\Omega)$ denote the closure of
$C^{\infty}_0(\Omega)$ in $W_p^m(\Omega)$.
We set $(u,v)_{L^2(\Omega)} =
\int_{\Omega} u v dx$ for functions $u, v$, while by
$(a,b)$ we denote the scalar product in $\Bbb R^2$ if there is no
fear of confusion. For $f:{\Bbb R}^2\rightarrow {\Bbb R}^1$, the
symbol $f''$ denotes the Hessian matrix  with entries
$\frac{\partial^2 f}{\partial x_k\partial x_j},$ $\mathcal L(X,Y)$
denotes the Banach space of all bounded linear operators from a
Banach space $X$ to another Banach space $Y$.
Let $E$ be the
$3\times 3$ unit matrix.
We set $ \Vert u\Vert_{W^{k,\tau}_2(\Omega)}
= (\Vert u\Vert_{W_2^k(\Omega)}^2 + \vert \tau\vert^{2k}\Vert
u\Vert^2_{L^2(\Omega)})^{\frac{1}{2}}$.
Let $C^*$ denote the adjoint operator of an operator $C$ in a
Hilbert space under consideration, and in particular, $C^*$ is the
transpose for a matrix $C$.  We set
$\vec{e}_1 = (1,0,0)$  $\vec{e}_2 = (0,1,0)$ and $\vec{e}_3 = (0,0,1)$.
\\

Let us introduce the operators:
$$
\partial_{\overline z}^{-1}g=-\frac 1\pi\int_\Omega
\frac{g(\xi_1,\xi_2)}{\zeta-z}d\xi_1d\xi_2,\quad
\partial_{ z}^{-1}g=-\frac 1\pi\int_\Omega
\frac{g(\xi_1,\xi_2)}{\overline\zeta-\overline z}d\xi_1d\xi_2.
$$

Then we have (e.g., p.47, 56, 72 in \cite{VE}):
\begin{proposition}\label{Proposition 3.0}
{\bf A)} Let $m\ge 0$ be an integer number and $\alpha\in (0,1).$
Then $\partial_{\overline z}^{-1},\partial_{ z}^{-1}\in \mathcal
L(C^{m+\alpha}(\overline{\Omega}),C^{m+\alpha+1}
(\overline{\Omega})).$
\newline
{\bf B}) Let $1\le p\le 2$ and $ 1<\gamma<\frac{2p}{2-p}.$ Then
 $\partial_{\overline z}^{-1},\partial_{ z}^{-1}\in
\mathcal L(L^p( \Omega),L^\gamma(\Omega)).$
\newline
{\bf C})Let $1< p<\infty.$ Then  $\partial_{\overline z}^{-1},
\partial_{ z}^{-1}\in \mathcal L(L^p( \Omega),W^1_p(\Omega)).$
\end{proposition}

Consider the operator
$$
Tf=\int_{\partial\Omega}\frac{f}{z-\xi}d\sigma.
$$

\begin{proposition}\label{drak} The operator $T$ is continuous operator
from $L^p(\partial\Omega)$ to $L^2(\Omega)$ for any $p>2.$
\end{proposition}

{\bf Proof.}  Using a partition of unity if necessary and changing the
coordinates, without loss of the generality we can assume that $\text{supp}\,
f \subset \Gamma$ where $\Gamma\subset \{ x_2=0\}.$  Let $\frac 1q=1-\frac 1p.$
$$
\Vert Tf\Vert_{L^2(\Omega)}
= \left\Vert \int_{\partial\Omega}\frac{f}{z-\xi}d\sigma\right\Vert
_{L^2(\Omega)}
\le \left\Vert \int_{\partial\Omega}\frac{\vert f\vert}{\vert z-\xi\vert}
d\sigma\right\Vert_{L^2(\Omega)}
$$
Applying the H\"older inequality we have
\begin{eqnarray}
 \left\Vert \int_{\partial\Omega}\frac{\vert f\vert}{\vert z-\xi\vert}
d\sigma\right\Vert_{L^2(\Omega)}\\
 \le \left\Vert \left(\int_{-\infty}^{+\infty}\frac{1}{\vert z-\xi\vert^q}
d\xi_1\right)^\frac 1q
\right\Vert_{L^2(\Omega)}
\Vert f\Vert_{L^p(\partial\Omega)}
\le C\left(\int_{-K}^{K} \frac{1}{\vert x_2\vert^{2(q-1)/q}}dx_2\right)
^\frac 12 \Vert f\Vert_{L^p(\partial\Omega)}.\nonumber
\end{eqnarray}
Since $q\in (1,2)$, the first integral in the above product is convergent.
The proof of the Proposition \ref{drak} is complete.
$\blacksquare $

Let us introduce the functional $\frak F_{\widetilde x,\tau}:
C^4(\overline\Omega)\rightarrow \Bbb R^1$:
$$
\frak F_{\widetilde x,\tau}u=\frac{\pi}{2}\left(\frac{u(\widetilde x)}{\tau}+\frac{-\partial_{zz}^2u(\widetilde x)+\partial^2_{\overline z\overline z}u(\widetilde x)}{4\tau^2}+\frac{\partial^4_{zzzz}u(\widetilde x)-2\partial^4_{zz\overline z\overline z}u(\widetilde x)+\partial^4_{\overline z\overline z\overline z\overline z}u(\widetilde x)}{32\tau^3}\right).
$$

We set
$$
\Phi(z) = (z - \widetilde{z})^2 - \kappa,
$$
where $\widetilde{z} = \widetilde{x}_1 + i \widetilde{x}_2$ and
$\kappa \in \Bbb R^1$ is an arbitrary constant.

We have
\begin{proposition}\label{osel}
Let $\Phi(z)=(z-\widetilde z)^2$ and $u\in C^{10}_0(\Omega)$ be some function.
Then the following asymptotic formula is true:
\begin{equation}\label{murzik9}
\int_{\Omega}ue^{\tau(\Phi-\overline \Phi)}dx=\frak F_{\widetilde x,\tau}u
+ o\left(\frac{1}{\tau^3}\right)\quad \mbox{as}\,\,\tau\rightarrow +\infty.
\end{equation}
\end{proposition}

{\bf Proof.} Without loss of generality, we can assume that $\widetilde z=0$.
By Theorem 7.7.5 of \cite{Her}, we have
\begin{equation}\label{murzik}
\int_{\Omega}we^{\tau(\Phi-\overline \Phi)}dx
= o\left(\frac{1}{\tau^3}\right)\quad\mbox{as}\,\,\tau\rightarrow +\infty
\end{equation}
for any function $w\in C^{10}_0(\Omega)$ such that $\partial^\beta_x w(0)=0$
for all $\vert \beta\vert\le 4$.
Let $e_1\in C^{\infty}_0(\Omega)$ be a function such that $e_1$ is identically
equal to one in some neighborhood of the origin.

Let
\begin{eqnarray*}
&& T_u(z,\overline z)=u(0)+z\partial_zu(0)+\overline z\partial_{\overline z}u(0)
+\frac{1}{2}( \partial^2_{zz}u(0)z^2+\partial^2_{\overline z\overline z}u(0)\overline z^2
+\partial^2_{\overline z z}u(0)\overline zz)\\
&& +\frac 12\partial^3_{z\overline z\overline z }u(0)z\overline z^2
+ \frac 16 \partial^3_{zzz}u(0)z^3+\frac 12 \partial^3_{\overline zzz}u(0)\overline zz^2
+\frac 16\partial^3_{\overline z\overline z\overline z}u(0)\overline z^3\\
&& + \frac {1}{24}\partial^4_{zzzz}u(0)z^4
+ \frac 14\partial^4_{zz\overline z\overline z} u(0)z^2\overline z^2
+ \frac {1}{24}\partial^4_{\overline z\overline z\overline z\overline z} u(0)\overline z^4
+ \frac 16\partial^4_{zzz\overline z}u(0)z^3\overline z
+ \frac 16\partial^4_{z\overline z\overline z\overline z}u(0)z\overline z^3).
\end{eqnarray*}
Then $T_u(z,\overline{z})$ is the Taylor polynomial of the function $u.$
By (\ref{murzik}), we have
 \begin{equation}\label{zina}
\int_{\Omega}ue^{\tau(\Phi-\overline \Phi)}dx=\int_{\Omega}
e_1T_ue^{\tau(\Phi-\overline \Phi)}dx+\int_{\Omega}(u-T_ue_1)
e^{\tau(\Phi-\overline \Phi)}dx
\end{equation}
$$
=\int_{\Omega}T_ue_1e^{\tau(\Phi-\overline \Phi)}dx
+ o\left(\frac{1}{\tau^3}\right)\quad\mbox{as}\,\,\tau\rightarrow +\infty.
$$

First we show that some terms from the first integral on the right-hand side
of (\ref{zina}) is estimated by $o\left(\frac{1}{\tau^3}\right)$.
Integrating by parts  and using the stationary phase argument we obtain
\begin{eqnarray}
\int_{\Omega}(z\partial_zu(0)+\overline z\partial_{\overline z}u(0))e_1
e^{\tau(\Phi-\overline \Phi)}dx \nonumber\\
= -\frac {1}{2\tau}\int_{\Omega}
(\partial_zu(0)\partial_z e_1-\partial_{\overline z}u(0)\partial
_{\overline z}e_1)e^{\tau(\Phi-\overline \Phi)}dx
= o\left(\frac{1}{\tau^3}\right)\quad\mbox{as}\,\,\tau\rightarrow +\infty,
\end{eqnarray}
\begin{eqnarray}
\int_{\Omega}z\overline z\partial^2_{z\overline z} u(0)e_1
e^{\tau(\Phi-\overline \Phi)}dx
=-\frac {1}{4\tau^2}\int_\Omega \partial^2_{z\overline z}
u(0)e_1\partial^2_{z\overline z}
e^{\tau(\Phi-\overline \Phi)}dx         \nonumber\\
= -\frac{ 1}{4\tau^2}\int_\Omega \partial^2_{z\overline z}u(0)
\partial^2_{z\overline z}e_1e^{\tau(\Phi-\overline \Phi)}dx
= o\left(\frac{1}{\tau^3}\right)
\quad\mbox{as}\,\,\tau\rightarrow +\infty,
\end{eqnarray}
\begin{eqnarray}
\int_{\Omega}
(\frac 12\partial^3_{z\overline z\overline z }u(0)z\overline z^2
+ \frac 16\partial^3_{zzz}u(0)z^3
+ \frac 12\partial^3_{\overline zzz}u(0)\overline zz^2
+ \frac 16\partial^3_{\overline z\overline z\overline z}u(0)
\overline z^3)e_1e^{\tau(\Phi-\overline \Phi)}dx\nonumber\\
= \int_{\Omega} (\partial^3_{z\overline z\overline z }u(0)\overline z\frac{1}
{8\tau^2}
\partial^2_{z\overline z}+\partial^3_{zzz}u(0)z^2\frac{1}{12\tau}\partial_{z}
+\partial^3_{\overline zzz}u(0) z\frac{1}{8\tau^2}\partial^2_{z\overline z}
+\partial^3_{\overline z\overline z\overline z}u(0)
\overline z^2\frac{1}{12\tau}\partial_{\overline z})
e_1e^{\tau(\Phi-\overline \Phi)}dx\nonumber\\
= o\left(\frac{1}{\tau^3}\right) \quad\mbox{as}\,\,\tau\rightarrow +\infty,
\end{eqnarray}
and
\begin{eqnarray}
\int_{\Omega}(\partial^4_{zzz\overline z}u(0)z^3\overline z
+ \partial^4_{z\overline z\overline z\overline z}u(0)z\overline z^3)
e_1e^{\tau(\Phi-\overline \Phi)}dx
                                            \nonumber\\
= \frac{-1}{4\tau^2}\int_{\Omega}(\partial^4_{zzz\overline z}u(0)z^2
\partial_{z\overline z}+\partial^4_{z\overline z\overline z\overline z}u(0)
\overline z^2\partial_{z\overline z})e_1e^{\tau(\Phi-\overline \Phi)}dx
= o\left(\frac{1}{\tau^3}\right)\quad\mbox{as}\,\,\tau\rightarrow +\infty.
\end{eqnarray}

The next two terms contribute to the asymptotic formula (\ref{murzik9}).
Integrating by parts we obtain
\begin{eqnarray}\label{zima}
\frac 12\int_{\Omega}(z^2\partial^2_{zz}u(0)
+\overline z^2\partial^2_{\overline z\overline z}u(0))
e_1e^{\tau(\Phi-\overline \Phi)}dx
=-\frac {1}{4\tau}\int_{\Omega}(\partial^2_{zz}u(0)\partial_z (ze_1)
- \partial^2_{\overline z\overline z}u(0)\partial_{\overline z}
(\overline z e_1))e^{\tau(\Phi-\overline \Phi)}dx                 \nonumber\\
= -\frac {1}{4\tau}\int_{\Omega}(\partial^2_{zz}u(0)e_1
-\partial^2_{\overline z\overline z}u(0) e_1)e^{\tau(\Phi-\overline \Phi)}dx
+ o\left(\frac{1}{\tau^3}\right)
\quad\mbox{as}\,\,\tau\rightarrow +\infty
\end{eqnarray}
and
\begin{eqnarray}\label{zima1}
\int_{\Omega}(\frac{1}{24}\partial^4_{zzzz}u(0)z^4+\frac{1}{24}
\partial^4_{\overline z\overline z\overline z\overline z}u(0)\overline z^4
+\frac 14\partial^4
_{ z z\overline z\overline z}
u(0)z^2\overline z^2)e_1e^{\tau(\Phi-\overline \Phi)}dx\\
= \int_{\Omega}(\frac{1}{48\tau}\partial^4_{zzzz}u(0)z^3\partial_{z}
-\frac{1}{48\tau}\partial^4_{\overline z\overline z\overline z\overline z}
u(0)\overline z^3
\partial_{\overline z}
-\frac{1}{16\tau^2}\partial^4_{ z z\overline z\overline z}u(0)z
\overline z\partial^2_{z\overline z})e_1e^{\tau(\Phi-\overline \Phi)}dx
                                                  \nonumber\\
= \int_{\Omega}(\frac{-3}{48\tau}\partial^4_{zzzz}u(0)z^2
+ \frac{3}{48\tau}\partial^4_{\overline z\overline z\overline z\overline z}
u(0)\overline z^2
-\frac{1}{16\tau^2}\partial^4_{ z z\overline z\overline z}u(0))e_1
e^{\tau(\Phi-\overline \Phi)}dx
+ o\left(\frac{1}{\tau^3}\right)                      \nonumber\\
= \int_{\Omega}(\frac{1}{32\tau^2}\partial^4_{zzzz}u(0)
+\frac{1}{32\tau^2}\partial^4_{\overline z\overline z\overline z\overline z}
u(0)
-\frac{1}{16\tau^2}\partial^4_{ z z\overline z\overline z}u(0))e_1
e^{\tau(\Phi-\overline \Phi)}dx+o(\frac{1}{\tau^3})
\,\,\mbox{as}\,\,\tau\rightarrow +\infty.\nonumber
\end{eqnarray}

Applying the stationary phase argument (see e.g  \cite{BH}) to (\ref{zima}),
(\ref{zima1}) and the integral $\int_\Omega u(0)e_1
e^{\tau(\Phi-\overline \Phi)}dx$,
we obtain (\ref{murzik9}).
$\blacksquare$




\bigskip

We recall that
$$
\Phi(z)=(z-\widetilde z)^2-\kappa,\quad
\widetilde z=\widetilde x_1+i\widetilde x_2. $$
Assume in addition
\begin{equation}\label{nol}
\kappa=1+sup_{z,\widetilde z\in \Omega}\mbox{Re}\, \Phi(z)
\end{equation}
and we define two operators:
\begin{equation}\label{anna}
\mathcal R_{\tau}g = \frac 12e^{\tau(\Phi-\overline {\Phi})}
\partial_{\overline z}^{-1}(g
e^{\tau(\overline{\Phi}- {\Phi})}),\,\, \widetilde {\mathcal
R}_{\tau}g = \frac 12e^{\tau(\overline {\Phi}-{\Phi})}
\partial_{ z}^{-1}(g e^{\tau( {\Phi}
-\overline {\Phi})}).
\end{equation} For any $g\in L^2(\Omega)$ we have
\begin{equation}
2\partial_{\overline z}\mathcal R_{\tau}g+2\tau\partial_{\overline z}
\overline{\Phi}\mathcal R_{\tau}g
= g\quad\mbox{and}\quad 2\partial_{ z}\widetilde{\mathcal R}_{\tau}g
+2\tau\partial_{ z}{\Phi}\widetilde{\mathcal R}_{\tau}g=g\quad\mbox{in}
\,\, \Omega.
\end{equation}

We have

\begin{proposition}\label{elkazelenaja}  Let $u\in W^1_p(\Omega)$
for any $p>1.$
Then for any $\epsilon\in (0,1)$ there exists a constant $C(\epsilon)$
such that
\begin{equation}\label{KJ}
\Vert \widetilde{\mathcal R}_{\tau} u\Vert_{L^2(\Omega)}
+\Vert {\mathcal R}_{\tau} u\Vert_{L^2(\Omega)}\le
C(\epsilon)/\tau^{1-\epsilon}.
\end{equation}
\end{proposition}
{\bf Proof.} We prove the estimate (\ref{KJ})  for the operator
$\widetilde {\mathcal R}_\tau.$ The proof of this estimate for the operator
${\mathcal R}_\tau$ is the same. Let $e_1, e_2 \in C^4(\overline\Omega)$
satisfy
$$
e_1(x)+e_2(x)=1\quad \mbox{in}\quad \Omega,  \quad e_1\in C_0^4(\overline
\Omega)
$$
and let $e_2$ have a support concentrated in a small neighborhood of
$\partial\Omega.$
For the function $e_1u$, the estimate
\begin{equation}\label{KJ1}
\Vert \widetilde{\mathcal R}_{\tau} (e_1u)\Vert_{L^2(\Omega)}\le
C(\epsilon)/\tau^{1-\epsilon}.
\end{equation}
was proved in \cite{IY4}.
Integrating by parts we obtain
\begin{eqnarray}\label{globus}
-\frac 1\pi\int_\Omega
\frac{(e_2 u)(\xi_1,\xi_2)e^{\tau( {\Phi}
-\overline {\Phi})}}{\overline\zeta-\overline z}d\xi_1d\xi_2
=\frac{e_2u e^{\tau(\Phi-\overline\Phi)}}{\tau\partial_z\Phi}
+ \frac 1\pi\int_\Omega
\frac{\partial_\zeta(e_2 u)(\xi_1,\xi_2)e^{\tau( {\Phi}
-\overline {\Phi})}}{\tau(\overline\zeta-\overline z)\partial_\zeta\Phi}
d\xi_1d\xi_2\nonumber\\
-\frac{1}{2\tau\pi}\int_{\partial\Omega}\frac{(\nu_1-i\nu_2)
(e_2 u)(\xi_1,\xi_2)e^{\tau( {\Phi}
-\overline {\Phi})}}{(\overline\zeta-\overline z)\partial_\zeta\Phi}
d\sigma.
\end{eqnarray}

Since the function $e_2$ is identically equal to zero in a neighborhood
of the critical points of the function $\Phi$, the function
$\frac{e_2u }{\partial_z\Phi} $ is continuous on $\overline \Omega$ and
the estimate
\begin{equation}\label{globus1}
\left\Vert \frac{e_2u e^{\tau(\Phi-\overline\Phi)}}
{\tau\partial_z\Phi}\right\Vert_{C^0(\overline \Omega)}\le \frac {C}{\tau}
\end{equation}
is seen. Applying Proposition \ref{Proposition 3.0} we obtain
\begin{equation}\label{globus2}
\left\Vert \int_\Omega
\frac{\partial_\zeta(e_2 u)(\xi_1,\xi_2)e^{\tau( {\Phi}
-\overline {\Phi})}}{\tau(\overline\zeta-\overline z)\partial_\zeta\Phi}
d\xi_1d\xi_2\right\Vert_{L^2(\Omega)}\le C/\tau.
\end{equation}
Since the function $\Phi$ does not have any critical points on
$\partial\Omega$, the function $\frac{(\nu_1-i\nu_2)(e_2 u)(\xi_1,\xi_2)
e^{\tau( {\Phi}-\overline {\Phi})}}{\partial_\zeta\Phi}$
is uniformly bounded in $\tau$
in $L^\infty(\partial\Omega).$ Then by Proposition \ref{drak}
\begin{equation}\label{globus3}
\left\Vert\frac{1}{2\tau\pi}\int_{\partial\Omega}
\frac{(\nu_1-i\nu_2)(e_2 u)(\xi_1,\xi_2)e^{\tau( {\Phi}
-\overline {\Phi})}}{(\overline\zeta-\overline z)\partial_\zeta\Phi}d\sigma
\right\Vert_{L^2(\Omega)}\le \frac{C}{\tau}.
\end{equation}

From (\ref{globus}) using (\ref{globus1})-(\ref{globus3}), we have
\begin{equation}\label{KJ11}
\Vert \widetilde{\mathcal R}_{\tau} (e_2u)\Vert_{L^2(\Omega)}\le
C(\epsilon)/\tau.
\end{equation}
Then (\ref{KJ}) and (\ref{KJ11}) imply (\ref{KJ}).
$\blacksquare$

We conclude this section with Carleman estimates for a partial differential
equation, which are fundamental in the arguments after
the next sections.

Consider the boundary value problem
\begin{equation}\label{spc}
\partial_{z}W+AW=f\quad\mbox{in}\,\,\Omega, \quad W\vert_{\partial\Omega}=0.
\end{equation}

Let $\widetilde \beta$ be a smooth function such that
\begin{equation}\label{ex}
\nabla \widetilde \beta (x)\ne 0\quad \mbox{on}\,\,\overline\Omega,
\quad \widetilde \beta(x)\ge 0\,\,\quad\mbox{on}\,\,\Omega,
\quad \mbox{min}_{x\in \overline\Omega}
\widetilde\beta(x)>\frac 34 \mbox{max}_{x\in \overline\Omega}
\widetilde \beta(x).
\end{equation}
We set
$$
\phi_s=e^{s\widetilde \beta}.
$$
\begin{proposition}\label{inga}
Let $A\in C^0(\overline\Omega)$ be a $3\times 3$
matrix.  Then there exist constants $s_0$ and $C$ independent of
$s$ such that for all $s\ge s_0$
\begin{equation}\label{vika1}
\int_\Omega \phi_s s^2  \vert W\vert^2e^{2s\phi_s}dx
\le C\Vert e^{s\phi_s} f\Vert_{L^2(\Omega)}^2, \quad W\in
\mathaccent'27{W}_2^1(\Omega).
\end{equation}
\end{proposition}

{\bf Proof.} Obviously it suffices to consider the case when $A\equiv 0.$
Denote $\widetilde v=We^{s\phi_s}, \widetilde f=fe^{s\phi_s},$
$L_-(x,D,s){\widetilde v}
= \frac 12\partial_{x_1}\widetilde{v} + \frac{1}{2}is
(\partial_{x_2}\phi_s)\widetilde v$
and $L_+(x,D,s){\widetilde v}=\frac {1}{2i} \partial_{x_2}
\widetilde{v} - \frac{1}{2}s(\partial_{x_1}\phi_s)\widetilde v.$
In the new notations we rewrite
equation (\ref{spc}) as

\begin{equation}\label{zanoza}
L_-(x,D,s){\widetilde v}+L_+(x,D,s){\widetilde v}
=\widetilde f
\quad\mbox{in}\quad\Omega.
\end{equation}
Taking the $L^2$- norm of the left- and right-hand sides of
(\ref{spc}), we obtain
\begin{eqnarray}\label{zanoza1!}
\Vert L_+ (x,D,s)\widetilde
 v\Vert^2_{L^2(\Omega)}+2\mbox{Re}\int_\Omega L_+(x,D,s)\widetilde v
\overline{L_-(x,D,s) \widetilde  v}dx \nonumber\\
+ \Vert L_-(x,D,s) \widetilde  v\Vert^2_{L^2(\Omega)}=\Vert \widetilde
f\Vert^2_{L^2(\Omega)}.
\end{eqnarray}
Integrating by parts the second term of (\ref{zanoza1!}), we obtain
\begin{equation}\label{zanoza2!}
2\mbox{Re}\int_\Omega L_+(x,D,s)\widetilde  v\overline{L_-(x,D,s)
\widetilde  v}dx =
\mbox{Re}\int_\Omega([L_+,L_-]\widetilde  v)\overline{ \widetilde  v}dx.
\end{equation}

Short computations yield
\begin{equation}
[L_+,L_-]=\frac{s}{4}\Delta\phi_s
= \left(\frac{s^2}{4}\vert \nabla\widetilde\beta\vert^2
+\frac{s}{4}\Delta\widetilde\beta\right)\phi_s.
\end{equation}
By (\ref{ex}) there exists a positive constant $C$ independent of $s$
such that
\begin{equation}\label{zanoza3!}
[L_+,L_-]\ge Cs^2\phi_s\quad \mbox{in}\,\,\Omega.
\end{equation}

From (\ref{zanoza1!})-(\ref{zanoza3!}) we obtain (\ref{vika1}).
$\blacksquare$

The following was proved in \cite{FI}.
\begin{proposition}
There exist constants $s_0$ and $C$ independent of $s$ such that for all
$s\ge s_0$
\begin{equation}
\int_\Omega \phi^3_s s^4  \vert W\vert^2e^{2s\phi_s}dx
\le C\Vert e^{s\phi_s} \Delta W\Vert_{L^2(\Omega)}^2\quad
\forall\, W\in  \mathaccent'27{W}_2^2(\Omega).
\end{equation}
\end{proposition}

\begin{corollary}\label{vis}
There exist constants $s_0$ and $C$ independent of $s$ such that
for all $s\ge s_0$
\begin{equation}\label{lazy1}
\int_\Omega \sum_{\vert\beta\vert\le 3} \phi_s^{2(3-\vert\beta\vert)}
s^{8-2\vert\beta\vert}\vert \partial_x^\beta W\vert^2e^{2s\phi_s}dx
\le C\Vert e^{s\phi_s} \Delta^2 W\Vert_{L^2(\Omega)}^2,
\quad\forall W\in \mathaccent'27{W}_2^4(\Omega).
\end{equation}
\end{corollary}
\section{Construction of special solutions to weakly coupling second-order
elliptic systems}\label{sec2}

Let $A_j(x), B_j(x), C_j(x)$, $j=1,2$, be smooth $3\times 3$ matrix
functions.
We set
$$
Q_1(j)=-2\partial_z A_j-B_jA_j+C_j,\quad Q_2(j)=-2\partial_{\overline z}
B_j-A_jB_j+C_j, \quad j=1,2.
$$

Let  $U_0, V_0\in C^{6+\alpha}(\overline \Omega)$
be nontrivial solutions to the differential equations:
\begin{equation}\label{-5}
2\partial_{\overline z}U_{0} +A_1
U_{0}=0\quad\mbox{in}\,\,\Omega
\end{equation}
and
\begin{equation}\label{-55}
2\partial_{ z}V_{0} +B_2
V_{0}=0\quad\mbox{in}\,\,\Omega .
\end{equation}

We have
\begin{proposition}(\cite {IY6}) \label{nikita}
Let
$\vec r_{0,k},\dots ,\vec r_{2,k}\in \Bbb C^3$ be arbitrary vectors
and $x_1,\dots, x_k$ be some points from the domain $\Omega.$
There exists a solution $U_{0}\in C^{6+\alpha}(\overline \Omega)$ to
problem (\ref{-5}) and a solution $V_{0}\in C^{6+\alpha}(\overline \Omega)$
to problem (\ref{-55}) such that
\begin{equation}\label{xoxo1}
\partial_z^j U_{0}(x_k)=\vec r_{j,k}\quad\mbox{and}\quad
\partial^j_{\overline z} V_{0}(x_k)=\vec r_{j,k}, \quad j=0,1,2,
\quad k\in \{1,\dots, 5\}.
\end{equation}
\end{proposition}

Let $B$ be a $3\times 3$ matrix with elements from $C^{5+\alpha}
(\overline\Omega)$ with $\alpha\in (0,1)$ and $\widehat x$ be some fixed
point from $\Omega.$
By Proposition \ref{nikita} for the equation
\begin{equation}\label{nina}
(2\partial_{\overline z}+B)u=0\quad
\mbox{in }\,\,\Omega,
\end{equation}
we  can construct solutions $U_{0,k}$ such that
$$
U_{0,k}(\hat x)=\vec e_k,\quad \forall k\in\{1,2,3\}.
$$
Consider the matrix
$$\Pi (x)=(U_{0,1}(x),U_{0,2}(x),U_{0,3}(x)).
$$
Then
$$
\left(\partial_{\overline z}+\frac 12\mbox{tr} B\right) \mbox{det} \Pi=0
\quad \mbox{in}\,\,\Omega.
$$
Hence there exists a holomorphic function $q(z)$ such that
$\mbox{det}\,\Pi=q(z)e^{- \frac 12 \partial^{-1}_z(\mbox{tr}\, B)}$
(see \cite{VE}).
By $\mathcal Q$ we denote the set of zeros of the function $q$ on $\overline
\Omega$ :  $\mathcal Q=\{z\in \overline\Omega\vert\, q(z)=0\}.$
Obviously $card\, \mathcal Q<\infty.$ By $\kappa$ we denote the highest order
of zeros of the function $q$ on $\overline \Omega.$

Using  Proposition \ref{nikita}
we construct solutions $\widetilde U_{0,k}$ to problem (\ref{nina}) such that
$$
\widetilde U_{0,k}( x)=\vec e_k\quad k\in\{1,2,3\}\quad
\forall x\in \mathcal Q.
$$
Set $\widetilde \Pi (x)=(\widetilde U_{0,1},\widetilde U_{0,2},
\widetilde U_{0,3}).$ Then there exists a holomorphic function $\widetilde q$
such that $\mbox{det}\,\widetilde \Pi=\widetilde q(z)e^{-\frac 12
\partial^{-1}_z(\mbox{tr}\, B)}.$
Let $\widetilde{\mathcal Q}=\{z\in \overline\Omega\vert \widetilde q(z)=0\}$
and $\widetilde \kappa$ the highest order of zeros of the function
$\widetilde q.$

By $\widetilde U_{0,k}(x) = \vec{e_k}$ for $x \in \mathcal{Q}$, we see that
$$
\widetilde{\mathcal Q}\cap\mathcal Q=\emptyset.
$$
Therefore there exists a holomorphic function $r(z)$ such that

$$r\vert_{\mathcal Q}=0\quad\mbox{and}\quad (1-r(z))
\vert_{\widetilde{\mathcal Q}}=0
$$
and the orders of zeros of the function $r$ on $\mathcal Q$ and the function
$1-r(z)$ on $\widetilde {\mathcal Q}$ are greater than or equal to the
$\max\{\kappa,\widetilde \kappa\}.$

We set
\begin{equation}\label{victory}
P_{B}f=\frac 12 \Pi\partial^{-1}_{\overline z} (\Pi^{-1}rf)+\frac 12
\widetilde \Pi\partial^{-1}_{\overline z} (\widetilde \Pi^{-1}(1-r)f).
\end{equation}
Then
$$
P_{B}^*f=-\frac 12 r (\Pi^{-1})^*\partial^{-1}_{\overline z} (\Pi^*f)
-\frac 12(1-r) (\widetilde \Pi^{-1})^*\partial^{-1}_{\overline z}
(\widetilde \Pi^*f).
$$

We have
\begin{proposition}
The linear operators $P_{B}, P_B^*\in {\mathcal L} (L^2(\Omega),
W^1_2(\Omega))$ solve the differential equations
\begin{equation}\label{oblom11}
(-2\partial_{\overline z}+B^*)P^*_{B}g=g, \quad
(2\partial_{\overline z}+B)P_{B}g=g\quad \mbox{in}\,\,\Omega.
\end{equation}
\end{proposition}
{\bf Proof.}
Since
$\partial_{\overline z}\Pi=-\frac 12 B\Pi$
and $\partial_{\overline z}\widetilde\Pi=-\frac 12 B\widetilde\Pi$,
short computations imply
$$
\partial_{\overline z}P_{B}f=\partial_{\overline z}\{\frac 12
\Pi\partial^{-1}_{\overline z} (\Pi^{-1}rf)+\frac 12
\widetilde \Pi\partial^{-1}_{\overline z} (\widetilde \Pi^{-1}(1-r)f)\}
$$
$$
= \frac 12(\partial_{\overline z} \Pi)\partial^{-1}_{\overline z}
(\Pi^{-1}rf)+\frac 12 (\partial_{\overline z}\widetilde \Pi)
\partial^{-1}_{\overline z} (\widetilde \Pi^{-1}(1-r)f)\}
$$
$$
+ \frac 12 \Pi (\Pi^{-1}rf)
+ \frac 12 \widetilde \Pi (\widetilde \Pi^{-1}(1-r)f)
$$
$$
= -\frac 14 B\Pi\partial^{-1}_{\overline z} (\Pi^{-1}rf)
- \frac 14 B\widetilde \Pi\partial^{-1}_{\overline z}
(\widetilde \Pi^{-1}(1-r)f)
$$
$$
+ \frac 12 rf+\frac 12 (1-r)f=-\frac 12 BP_{B}f+\frac 12 f.
$$
Hence the second equality in (\ref{oblom11}) is proved. In order to prove
the first one observe that
 since
$
\Pi \Pi^{-1}=E
$ on $\overline\Omega\setminus \mathcal Q$. The differentiation of this
identity gives
$$
0=\partial_{\overline z} (\Pi\Pi^{-1})=\partial_{\overline z}
\Pi\Pi^{-1}+\Pi\partial_{\overline z}\Pi^{-1}=-\frac 12 B\Pi\Pi^{-1}
+\Pi\partial_{\overline z}\Pi^{-1}.
$$
This equality can be written as
$
\Pi\partial_{\overline z}\Pi^{-1}=\frac 12 B.
$
Multiplying both sides of this equality by $ \Pi^{-1}$ we have
$$
\partial_{\overline z}\Pi^{-1}=\frac 1 2 \Pi^{-1} B\quad \mbox{on}\quad
\overline\Omega\setminus \mathcal Q.
$$
Next we take the adjoint for the left- and the right-hand sides of the above
equality:
$$
 (\partial_{\overline z}\Pi^{-1})^*
= \partial_{\overline z}(\Pi^{-1})^*=(\frac 1 2 \Pi^{-1} B)^*
=\frac 12 B^*(\Pi^{-1})^*\quad \mbox{on}\quad\overline\Omega\setminus
\mathcal Q.
$$
Observing that $(\Pi^{-1})^*=(\Pi^*)^{-1}$, we obtain
\begin{equation}\label{XZ}
 \partial_{\overline z}(\Pi^*)^{-1}=\frac 12 B^*(\Pi^*)^{-1}\quad
\mbox{on}\quad\overline\Omega\setminus \mathcal Q.
\end{equation}
 Similarly we obtain
 \begin{equation}\label{ZX}
\partial_{\overline z}(\widetilde\Pi^*)^{-1}
=\frac 12 B^*(\widetilde \Pi^*)^{-1}\quad \mbox{on}\quad
\overline\Omega\setminus\widetilde{ \mathcal Q}.
\end{equation}
Setting
$$
\Gamma=e^{\frac 12 \partial^{-1}_z(\mbox{tr}\, B)}\left(\begin{matrix}
\Pi_{11} &\Pi_{21} &\Pi_{31}\\
 \Pi_{12} &\Pi_{22} &\Pi_{32}\\
 \Pi_{13} &\Pi_{23} &\Pi_{33}\end{matrix}\right),\quad \widetilde\Gamma
= e^{\frac 12 \partial^{-1}_z(\mbox{tr}\, B)}\left(\begin{matrix}
\widetilde \Pi_{11} &\widetilde \Pi_{21} &\widetilde \Pi_{31}\\
 \widetilde \Pi_{12} &\widetilde \Pi_{22} &\widetilde \Pi_{32}\\
 \widetilde \Pi_{13} &\widetilde \Pi_{23} &\widetilde \Pi_{33}\end{matrix}
\right),
$$
we can write the matrices $(\Pi^*)^{-1},(\widetilde\Pi^*)^{-1}$ as
$$
 (\Pi^*)^{-1}=\frac{1}{q(z)}\Gamma^*,\quad (\widetilde\Pi^*)^{-1}
=\frac{1}{\widetilde q(z)}\widetilde\Gamma^*.
$$
Then (\ref{XZ}) and (\ref{ZX}) imply
$$
\partial_{\overline z}(\Gamma^*)^{-1}-\frac 12 B^*(\Gamma^*)^{-1}=
\partial_{\overline z}(\widetilde\Gamma^*)^{-1}
-\frac 12 B^*(\widetilde \Gamma^*)^{-1}=0\quad \mbox{on}\quad\overline\Omega.
$$
Since $r/q$ and $(1-r)/\widetilde q$ are smooth functions,
the above equalities yield
\begin{equation}\label{ZZZ}
\partial_{\overline z}(r(z)(\Pi^*)^{-1})
= \frac {r(z)}{2} B^*(\Pi^*)^{-1}, \quad
\partial_{\overline z}((1-r(z))(\widetilde\Pi^*)^{-1})
= \frac {1-r(z)}{2} B^*(\widetilde \Pi^*)^{-1}\quad \mbox{in}\quad\Omega.
\end{equation}
Using (\ref{ZZZ}), we compute
$$
\partial_{\overline z}P_{B}^*f=-\partial_{\overline z}\{\frac 12 r(z)
(\Pi^{-1})^*\partial^{-1}_{\overline z} (\Pi^*f)+\frac 12(1-r(z))
(\widetilde \Pi^{-1})^*\partial^{-1}_{\overline z} (\widetilde \Pi^*f)\}
$$
$$
= -\frac 12 \partial_{\overline z}( r(z)(\Pi^{-1})^*)\partial^{-1}
_{\overline z} (\Pi^*f)-\frac 12 \partial_{\overline z}
( (1-r(z))(\widetilde\Pi^{-1})^*)\partial^{-1}_{\overline z}
(\widetilde \Pi^*f)
$$
$$
-\frac 12  r(z)(\Pi^{-1})^* \Pi^*f-\frac 12 (1-r(z))(\widetilde\Pi^{-1})^*
\widetilde\Pi^*f
$$
$$
= -\frac 14  r(z)B^*(\Pi^{-1})^*\partial^{-1}_{\overline z} (\Pi^*f)-\frac 14
(1-r(z))B^*(\widetilde\Pi^{-1})^*\partial^{-1}_{\overline z}
(\widetilde \Pi^*f)
$$
$$
-\frac 12  r(z)f-\frac 12 (1-r(z))f=\frac 12 B^*P_{B}^*f-\frac 12 f.
$$
The proof of the proposition is complete. $\blacksquare$

Short computations imply
\begin{eqnarray}\label{lebanon1}
[\partial_z,P_B]f=\frac 12 \partial_z\Pi\partial_{\overline z}^{-1}
(\Pi^{-1} rf)-\frac{\Pi}{8\pi}\int_{\partial\Omega}
\frac{(\nu_1-i\nu_2)\Pi^{-1} r f}{z-\zeta} d\zeta\nonumber\\
+ \frac 12 \Pi\partial_{\overline z}^{-1}(\partial_z(\Pi^{-1} r(z))f)
                                            \\
+ \frac 12 \partial_z\widetilde\Pi\partial_{\overline z}^{-1}
(\widetilde\Pi^{-1} (1-r)f)-\frac{\widetilde\Pi}{8\pi}\int_{\partial\Omega}
\frac{(\nu_1-i\nu_2)\widetilde\Pi^{-1} (1-r)f }{z-\zeta} d\zeta
+ \frac 12 \widetilde\Pi\partial_{\overline z}^{-1}
(\partial_z(\widetilde\Pi^{-1} (1-r))f).\nonumber
\end{eqnarray}

In a similar way we construct matrices $\Pi_0, \widetilde \Pi_0$,
antiholomorphic function $r_0(\overline z)$ and operators
\begin{equation}\label{giorgi}
T_{B}f=\frac 12 \Pi_0\partial^{-1}_{z} (\Pi_0^{-1} r_0(\overline z)f)
+ \frac 12 \widetilde \Pi_0\partial^{-1}_{ z}
(\widetilde \Pi_0^{-1}(1-r_0(\overline z))f)
\end{equation}
and
\begin{equation}\label{giorgi1}
T_{B}^*f=-\frac 12  r_0(\overline z) (\Pi_0^{-1})^*\partial^{-1}_{ z}
(\Pi_0^*f)-\frac 12(1-r_0(\overline z)) (\widetilde \Pi_0^{-1})^*
\partial^{-1}_{z} (\widetilde \Pi_0^*f).
\end{equation}

For any matrix $B\in C^{1}(\overline \Omega)$,  the linear
operators $T_{B}$ and $T^*_{B}$ solve the differential equation
\begin{equation}\label{oblom1}
(2\partial_z+B)T_{B}g=g\quad \mbox{in}\,\,\Omega; \quad
(-2\partial_{ z}+B^*)T^*_{B}g=g\quad \mbox{in}\,\,\Omega.
\end{equation}

By (\ref{giorgi}) and (\ref{giorgi1}), we have
\begin{eqnarray}\label{lebanon}
[\partial_{\overline z},T_B]f=\frac 12 \partial_{\overline z}
\Pi_0\partial_{ z}^{-1}(\Pi_0^{-1} r_0f)
-\frac{\Pi_0}{8\pi}\int_{\partial\Omega}
\frac{(\nu_1+i\nu_2)\Pi_0^{-1} r_0 f}{\overline z-\overline \zeta}d\zeta
                                       \nonumber\\
+ \frac 12 \Pi_0(\partial_{z}^{-1}(\partial_{\overline z}
(\Pi_0^{-1}) r)f)\\
+ \frac 12 \partial_{\overline z}\widetilde\Pi_0\partial_{ z}^{-1}
(\widetilde\Pi_0^{-1}(1- r_0)f)-\frac{\widetilde\Pi_0}{8\pi}
\int_{\partial\Omega}
\frac{(\nu_1+i\nu_2)\widetilde\Pi_0^{-1}(1- r_0)f }
{\overline z-\overline \zeta} d\zeta+\frac 12 \widetilde\Pi_0
\partial_{ z}^{-1}(\partial_{\overline z}(\widetilde\Pi_0^{-1}(1- r_0))f).
                               \nonumber
\end{eqnarray}

The formulae (\ref{lebanon}), (\ref{lebanon1}),
Propositions \ref{Proposition 3.0} and \ref{drak} imply the
following proposition.
\begin{proposition}\label{pen} Let $B\in C^5(\overline \Omega)$. Then
\label{zanuda}
$$
[\partial_{\overline z},T_B],[\partial_{ z},P_B]\in
\mathcal L(C^0(\overline\Omega),L^2(\Omega)).
$$
\end{proposition}

Next we investigate the properties of the commutators of the operators
$P_B,T_B$ and $\partial_z,\partial_{\overline z}.$
\begin{proposition}\label{zanuda1}
Let $q\in C^1(\overline\Omega)$ and $q(\widetilde x)=0.$ Assume that
restriction of the  function $\psi$ on $\partial\Omega$ has
a finite number of critical points and all these points are
nondegenerate. Then  there exist positive constants $C$ and $\tau_0$
independent of $\tau $ such that
\begin{equation}\label{zinaida}
\Vert -{\partial_z\Phi} [\partial_z,P_{B}] (q e^{\tau(\Phi-\overline\Phi)})
+ [\partial_z,P_{B}] (q{\partial_z\Phi} e^{\tau(\Phi-\overline\Phi)})\Vert
_{L^2(\Omega)}
\le \frac{C}{\tau^{\frac 12}}\quad \forall \tau\ge\tau_0.
\end{equation}
Moreover
\begin{equation}\label{zinaida1}
[\partial_z,[\partial_z,P_B]]\in {\mathcal L}(W^s_2(\Omega),L^2(\Omega))
\quad \forall s >\frac 12.
\end{equation}
\end{proposition}

{\bf Proof.}
From the explicit formula (\ref{lebanon1}) for the operator
$[\partial_z,P_{B}]$, we have
\begin{eqnarray}
 -{\partial_z\Phi} [\partial_z,P_{B}] (q e^{\tau(\Phi-\overline\Phi)})
+ [\partial_z,P_{B}] (q{\partial_z\Phi} e^{\tau(\Phi-\overline\Phi)})\nonumber\\
= \frac 12 \partial_z\Pi\partial_{\overline z}^{-1}
(\Pi^{-1} r (\Phi(\zeta)-\Phi(z))qe^{\tau(\Phi-\overline\Phi)})
-\frac{\Pi}{8\pi}\int_{\partial\Omega}
\frac{(\nu_1-i\nu_2)\Pi^{-1} r (\Phi(\zeta)-\Phi(z))q
e^{\tau(\Phi-\overline\Phi)}}{z-\zeta} d\zeta\nonumber\\+\frac 12
\Pi\partial_{\overline z}^{-1}(\partial_z(\Pi^{-1} r)(\Phi(\zeta)-\Phi(z))
qe^{\tau(\Phi-\overline\Phi)})\nonumber\\
+ \frac 12 \partial_z\widetilde\Pi\partial_{\overline z}^{-1}
(\widetilde\Pi^{-1} (1-r)(\Phi(\zeta)-\Phi(z))qe^{\tau(\Phi-\overline\Phi)})\nonumber\\
-\frac{\widetilde\Pi}{8\pi}\int_{\partial\Omega}
\frac{(\nu_1-i\nu_2)\widetilde\Pi^{-1} (1-r)(\Phi(\zeta)-\Phi(z))q
e^{\tau(\Phi-\overline\Phi)}}{z-\zeta} d\zeta\nonumber\\
+ \frac 12 \widetilde\Pi\partial_{\overline z}^{-1}(\partial_z
(\widetilde\Pi^{-1}(1- r))(\Phi(\zeta)-\Phi(z))qe^{\tau(\Phi-\overline\Phi)})
                 \nonumber\\
= \frac {1}{\pi} \partial_z\Pi\int_\Omega\Pi^{-1} r q
e^{\tau(\Phi-\overline\Phi)}d\zeta+\frac{\Pi}{4\pi}\int_{\partial\Omega}
(\nu_1-i\nu_2)
\Pi^{-1} r qe^{\tau(\Phi-\overline\Phi)} d\zeta\nonumber\\
+ \frac 1\pi \Pi\int_\Omega\partial_z(\Pi^{-1} r)
(qe^{\tau(\Phi-\overline\Phi)})d\zeta\nonumber\\
+ \frac 1\pi \partial_z\widetilde\Pi\int_\Omega\widetilde\Pi^{-1}
(1-r)qe^{\tau(\Phi-\overline\Phi)}d\zeta
+ \frac{\widetilde\Pi}{4\pi}\int_{\partial\Omega}
\widetilde\Pi^{-1} (1-r) q(\nu_1-i\nu_2)e^{\tau(\Phi-\overline\Phi)}
d\zeta\nonumber\\
+ \frac{ 1}{\pi} \widetilde\Pi\int_\Omega\partial_z
(\widetilde\Pi^{-1} (1-r))(qe^{\tau(\Phi-\overline\Phi)})d\zeta.
\end{eqnarray}
This equality implies (\ref{zinaida}) immediately. In order to prove
(\ref{zinaida1}), we observe that the operator $[\partial_z,P_{B}]$ is a
linear combination of the operators ${\bf B}_{ij}f
= \Pi_i\partial^{-1}_{\overline z}(\Pi_j f)$ and ${\bf T}_{ij}(f)
= \Pi_i\int_{\partial\Omega}\frac{\Pi_j f}{z-\zeta}d\sigma$
where $\Pi_k$ are some regular matrices.
Short computations imply
\begin{eqnarray}
[\partial_z,{\bf B}_{ij}]f=\partial_z\Pi_i\partial^{-1}_{\overline z}(\Pi_j f)
+ \Pi_i\partial^{-1}_{\overline z}(\partial_z\Pi_j f)
- \frac {1}{4\pi}\Pi_i\int_{\partial\Omega}\frac{(\nu_1-i\nu_2)\Pi_j f}
{z-\zeta}d\sigma.
\end{eqnarray}

This formula and Propositions \ref{Proposition 3.0} and \ref{drak} imply
$$
[\partial_z,{\bf B}_{ij}]\in\mathcal L (W_2^s(\Omega),L^2(\Omega))
              \quad \forall s>1/2.
$$
On the other hand, there exists  a function $a_{ij}\in C^0(\partial\Omega)$
such that
$$
[\partial_z,{\bf T}_{ij}]=\int_{\partial\Omega}\frac{a_{ij}f}{z-\zeta}d\sigma.
$$
Applying Proposition \ref{drak} we have
$$
[\partial_z,{\bf T}_{ij}]\in\mathcal L (W_2^s(\Omega),L^2(\Omega))
\quad \forall s>1/2.
$$
The proof of the proposition is complete.
$\blacksquare$

Next we introduce two operators
$$
\widetilde {\mathcal R}_{\tau, B}g = e^{\tau(\overline\Phi-\Phi)}T_B
(e^{\tau(\Phi-\overline\Phi)}g), \quad
{\mathcal R}_{\tau, B}g = e^{\tau(\Phi-\overline\Phi)}
P_B(e^{\tau(\overline\Phi-\Phi)}g).
$$

For any $g \in C^\alpha(\overline\Omega)$, the functions
$\mathcal R_{\tau, B}g$ and $\widetilde {\mathcal R}_{\tau, B}g$
solve the equations:
\begin{equation}\label{0011}
(2\partial_{\overline z}+2\tau \partial_{\overline z}\overline \Phi+B) \mathcal
R_{\tau, B} g=g, \quad (2\partial_{ z}+2\tau
\partial_{ z} \Phi+B)\widetilde {\mathcal R}_{\tau, B} g=g\quad
\mbox{in}\,\, \Omega.
\end{equation}

Consider the system of linear equations
\begin{equation}\label{-10}
P_j(x,D)W := \Delta W+2A_j\partial_z W+2B_j\partial_{\overline
z}W+C_jW=0\quad\mbox{in}\,\,\,\Omega, \thinspace j=1,2.
\end{equation}

We have
\begin{proposition}\label{mursilka0}
Let $A_j,B_j\in C^2_0(\Omega)$, $j=1,2$, the functions $U_0, V_0$ be
given by Proposition \ref{nikita} and the functions
$\mathcal A_1,\mathcal D_2$ be some solutions to equations (\ref{-5})
and (\ref{-55}) respectively.  Then the functions $U$ and $V$ defined by
\begin{eqnarray}\label{mursilka}
U=e^{\tau\Phi}( U_0-U_1)+\sum_{j=2}^\infty(-1)^jU_j e^{\tau\overline \Phi},
\quad  U_j={\mathcal
R}_{\tau,A_1}(T_{B_1}(Q_1(1)U_{j-1})) \quad j\ge 3, \nonumber\\
U_1 =\widetilde{\mathcal
R}_{\tau,B_1}(P_{A_1}(Q_1(1)U_{0})-\mathcal A_1),\quad U_2
={\mathcal R}_{\tau, A_1}(T_{B_1} (Q_2(1)e^{\tau(\Phi-\overline\Phi)}U_1));
\end{eqnarray}
\begin{eqnarray}\label{mursilka1}
V=e^{-\tau\overline \Phi}( V_0-V_1)
+ \sum_{j=2}^\infty(-1)^jV_j e^{-\tau \Phi}, \quad  V_j
=\widetilde {\mathcal R}_{-\tau,B_2}(P_{A_2}(Q_2(2)V_{j-1})) \quad j\ge 3,
                              \nonumber\\
V_1={\mathcal
R}_{-\tau,A_2}(T_{B_2}(Q_2(2)V_{0})-\mathcal D_2), \quad
V_2=\widetilde {\mathcal R}_{-\tau, B_1}(P_{A_2} (Q_1(2)
e^{-\tau(\Phi-\overline\Phi)}V_1)
\end{eqnarray}
are solutions to the system (\ref{-10}) with $j=1,2$ respectively.
\end{proposition}

{\bf Proof.}
All the terms of the infinite series (\ref{mursilka}) and (\ref{mursilka1})
are correctly defined. First we show that  these infinite series are
convergent in $L^q(\Omega)$ with some $q>2.$

Let $p>2$ and $\widetilde p\in (2,p)$ be some numbers.
By Proposition \ref{elkazelenaja} there exist positive constants $C$ and
$\delta=\delta(\widetilde p)$ and $\tau_1$ such that
\begin{equation}\label{gastronom}
\Vert {\mathcal R}_{\tau,A_1}u\Vert_{L^\frac{p\widetilde p}
{p-\widetilde p}(\Omega)}
\le C\Vert u\Vert_{W^1_{\widetilde p}(\Omega)}/\tau^\delta
\quad\forall\tau\ge \tau_1,
\end{equation}
where $C$ is independent of $\tau.$
Hence, (\ref{gastronom}) and Proposition \ref{Proposition 3.0} yield
\begin{eqnarray}\label{novel}
\Vert U_j\Vert_{L^\frac{p\widetilde p}{p-\widetilde p}(\Omega)}
\le \frac{C}{\tau^\delta}\Vert T_{B_1}(Q_1(1)U_{j-1})
\Vert_{W^1_{\widetilde p}(\Omega)}\le \frac{C}{\tau^\delta}
\Vert T_{B_1}\Vert_{\mathcal L(L^{\widetilde p}(\Omega),
W^1_{\widetilde p}(\Omega))}
\Vert Q_1(1)\Vert_{L^p(\Omega)}\Vert U_{j-1}\Vert
_{L^\frac{p\widetilde p}{p-\widetilde p}(\Omega)}\nonumber\\
\le (\frac{C}{\tau^\delta}\Vert T_{B_1}\Vert
_{\mathcal L(L^{\widetilde p}(\Omega), W^1_{\widetilde p}(\Omega))}
\Vert Q_1(1)\Vert_{L^p(\Omega)})^{j-1}\Vert  U_1\Vert
_{L^\frac{p\widetilde p}{p-\widetilde p}(\Omega)}.
\end{eqnarray}
Hence, there exists $\tau_0$ such that for all $\tau\ge \tau_0$
\begin{eqnarray}\label{novel1}
\Vert U_j\Vert_{L^\frac{p\widetilde p}{p-\widetilde p}(\Omega)}
\le \frac{1}{2^j}\Vert U_1\Vert_{L^\frac{p\widetilde p}
{p-\widetilde p}(\Omega)}\quad \forall j\ge 2.
\end{eqnarray}
Therefore the infinite series (\ref{mursilka}) is convergent.
Next we show that the function $U$ solves the equation (\ref{-10}).
Denote
$U(N,\tau,x)=(U_0-U_1)e^{\tau\Phi}+\sum_{j=2}^N(-1)^jU_j
e^{\tau\overline \Phi}.$
Thanks to the representation of the operator $P_1(x,D)$ in the form
\begin{eqnarray}\label{ooo} P_1(x,D)=4\partial_z\partial_{\overline z}+2
A_1\partial_z+2B_1\partial_{\overline z} +C_1\\
= (2\partial_z+B_1)(2\partial_{\overline z}+A_1)-2\partial_z A_1-B_1A_1+C_1
=(2\partial_z+B_1)(2\partial_{\overline z}+A_1)+Q_1(1)
\nonumber\\
= (2\partial_{\overline z}+A_1)(2\partial_z+B_1)-2\partial_{\overline z}
B_1-A_1B_1+C_1=(2\partial_{\overline z}+A_1)(2\partial_z+B_1)+Q_2(1),
                                                   \nonumber
\end{eqnarray}
we have  $P_1(x,D)(U_0e^{\tau \Phi})=Q_1(1)U_0e^{\tau \Phi} $ and
$P_1(x,D)(U_1e^{\tau\Phi})=Q_2(1)U_1e^{\tau\Phi}+Q_1(1)U_0e^{\tau\Phi}.$
$P_1(x,D)(U_2e^{\tau\overline\Phi})= Q_1(1)(U_2e^{\tau\overline\Phi})+Q_2(1)
U_1e^{\tau\Phi}.$ Observe that
\begin{equation}\label{J0}
P_1(x,D)(U_0e^{\tau\Phi}-U_1e^{\tau\Phi}+U_2e^{\tau\overline\Phi})
= Q_1(1)(U_2e^{\tau\overline\Phi}).
\end{equation}

For $j\ge 3$, we have
\begin{eqnarray}\label{J1}
P_1(x,D)(e^{\tau\overline \Phi}U_j)=((2\partial_z+B_1)(2\partial_{\overline z}
+A_1)-2\partial_z A_1-B_1A_1+C_1)(e^{\tau\overline \Phi}U_j)\nonumber\\
= ((2\partial_z+B_1)(2\partial_{\overline z}+A_1)+Q_1(1))
(e^{\tau\overline \Phi}U_j)=Q_1(1)U_{j-1}e^{\tau\overline \Phi}
+ Q_1(1)U_je^{\tau\overline\Phi}.
\end{eqnarray}
Then by (\ref{J0}) and (\ref{J1}), we obtain
\begin{equation}\label{JJ}
P_1(x,D)U(N,\tau,\cdot)=Q_1(1)U_N e^{\tau \Phi}\quad \mbox{in}\,\,\Omega.
\end{equation}
Since by (\ref{novel}) we have
$$
\Vert U_j\Vert_{L^2(\Omega)}\rightarrow 0\quad\mbox{as}\quad
j\rightarrow +\infty,
$$
passing to the limit in (\ref{JJ}) as $N\rightarrow +\infty$,
we complete the proof of the fact that
the function $U$ solves the equation (\ref{-10}).

In order to prove a similar statement for the function $V$,
observe that
$P_2(x,D)(V_0e^{-\tau \overline\Phi})=Q_2(2)V_0e^{-\tau \overline\Phi} $
and $P_2(x,D)(V_1e^{-\tau\overline\Phi})=Q_1(2)V_1e^{-\tau\overline\Phi}
+ Q_2(2)V_0e^{-\tau\overline\Phi}$ and
$P_2(x,D)(V_2e^{-\tau\Phi})= Q_1(2)(U_2e^{-\tau\Phi})
+Q_2(2)V_1e^{-\tau\overline\Phi}U_1.$   We further note that
\begin{equation}\label{J0'}
P_1(x,D)(V_0e^{-\tau\overline\Phi}-V_1e^{-\tau\overline\Phi}
+ V_2e^{-\tau\Phi})=Q_2(2)V_2e^{-\tau\Phi}
\end{equation}
and
\begin{equation}\label{J0''}
P_2(x,D)(V_je^{-\tau\Phi})=Q_2(2)(V_{j-1}e^{-\tau\Phi}).
\end{equation}
Therefore, for the function
$V(N,\tau,x)=(V_0-V_1)e^{-\tau\overline\Phi}+\sum_{j=2}^N(-1)^jV_j
e^{-\tau \Phi}$, we have
\begin{equation}\label{JJ}
P_2(x,D)V(N,\tau,\cdot)=Q_2(1)V_N e^{-\tau\overline \Phi}
\quad \mbox{in}\,\,\Omega.
\end{equation}
Using arguments similar to the above, we complete the proof of the proposition.
 $\blacksquare$

Let
$$
\widetilde \Phi(z)=e^{i\theta}z, \quad \mbox{where $\theta\in [0,2\pi)$}.
$$
We have

\begin{proposition}\label{mursilka00}
Let $A_j,B_j\in C^2_0(\Omega)$, $j=1,2$, the functions $U_0, V_0$ be given
by Proposition \ref{nikita} and the functions $\mathcal A_1,\mathcal D_2$ be
some solutions to equations (\ref{-5}) and (\ref{-55}) respectively.
Then the functions $\widetilde{U}$ and $\widetilde{V}$ defined by
\begin{eqnarray}\label{mursilka++}
\widetilde U=e^{\tau \widetilde \Phi}( \widetilde U_0-\widetilde U_1)
+ \sum_{j=2}^\infty(-1)^j\widetilde U_j e^{\tau\overline{\widetilde \Phi}},
\quad  \widetilde U_j={\mathcal
R}_{\tau,A_1}(T_{B_1}(Q_1(1)\widetilde U_{j-1})) \quad j\ge 3 \nonumber\\
\widetilde  U_1=\widetilde{\mathcal
R}_{\tau,B_1}(P_{A_1}(Q_1(1)\widetilde U_{0})-\mathcal A_1),\quad
\widetilde U_2={\mathcal R}_{\tau, A_1}(T_{B_1} (Q_2(1)
e^{\tau(\widetilde\Phi-\overline{\widetilde \Phi})}\widetilde U_1));\nonumber\\
\widetilde U_0=z^kU_0,
\end{eqnarray}
and
\begin{eqnarray}\label{mursilka+++}
\widetilde V=e^{-\tau\overline{\widetilde \Phi}}(\widetilde  V_0
-\widetilde V_1)+\sum_{j=2}^\infty(-1)^j\widetilde V_j
e^{-\tau \widetilde\Phi}, \quad  \widetilde V_j=\widetilde {\mathcal
R}_{-\tau,B_2}(P_{A_2}(Q_2(2)\widetilde V_{j-1})) \quad j\ge 3\nonumber\\
\widetilde V_1={\mathcal
R}_{-\tau,A_2}(T_{B_2}(Q_2(2)\widetilde V_{0})-\mathcal D_2),
\quad \widetilde V_2=\widetilde {\mathcal R}_{-\tau, B_1}(P_{A_2}
(Q_1(2)e^{-\tau(\widetilde\Phi-\overline{\widetilde\Phi})}\widetilde V_1),
                                      \nonumber\\
\widetilde V_0=z^kV_0
\end{eqnarray}
are solutions to the system (\ref{-10}) with $j=1,2$ respectively.
\end{proposition}

Consider the system of linear equations
\begin{equation}\label{-101}
P_1(x,D)Z=\Delta Z+2A_1\partial_z Z+2B_1\partial_{\overline
z}Z+C_1Z=e^{\tau \Phi}f\quad\mbox{in}\,\,\,\Omega.
\end{equation}

The following Proposition establish the solvability of the equation
(\ref{-101}).

\begin{proposition}\label{vanka}
Let $A_1,B_1\in C^2_0(\Omega).$   Then the function
\begin{equation}\label{mursilkaa}
Z=\sum_{j=0}^\infty(-1)^jZ_j e^{\tau \Phi}, \quad
Z_j=\widetilde{\mathcal
R}_{\tau,B_1}(P_{A_1}(Q_2(1)Z_{j-1}))\quad j\ge 2,
 \quad Z_0=\widetilde{\mathcal
R}_{\tau,B_1}(P_{A_1}f)
\end{equation}
is a solution to the system (\ref{-101}). Moreover the following estimate
holds true:
\begin{equation}\label{zopa}
\Vert e^{-\tau\varphi} Z\Vert_{W^2_2(\Omega)}
+ \vert \tau\vert^2 \Vert e^{-\tau\varphi} Z\Vert_{L^2(\Omega)}\le C(1+\tau^2)
\Vert fe^{-\tau\varphi}\Vert_{L^2(\Omega)}.
\end{equation}
\end{proposition}

{\bf Proof.}
By Proposition \ref{Proposition 3.0}, we have $Z_j\in W^2_2(\Omega)$
for each $j\ge 0.$

Short computations imply the formulae :
$$
P_1(x,D)(e^{\tau\Phi}Z_0)=fe^{\tau\Phi}+Q_2(1)Z_0e^{\tau\Phi}
\quad \mbox{and}\quad
P_1(x,D)(e^{\tau\Phi}Z_j)=Q_2(1)Z_je^{\tau\Phi}+Q_2(1)Z_{j-1}e^{\tau\Phi}.
$$
Using these formulae we obtain
\begin{eqnarray}\label{gnomik}
P_1(x,D)(\sum_{j=0}^N (-1)^jZ_je^{\tau\Phi})=e^{\tau\Phi}f
+ Q_2(1)Z_0e^{\tau\Phi}+P_1(x,D)(\sum_{j=1}^N (-1)^jZ_je^{\tau\Phi})\nonumber\\
= e^{\tau\Phi}f+Q_2(1)Z_0e^{\tau\Phi}+P_1(x,D)
(\sum_{j=1}^N (-1)^j(Q_2(1)Z_je^{\tau\Phi}+Q_2(1)Z_{j-1}e^{\tau\Phi}))\\
=e^{\tau\Phi}f+Q_2(1)Z_0e^{\tau\Phi}-Q_2(1)Z_0e^{\tau\Phi}
+(-1)^NQ_2(1)Z_Ne^{\tau\Phi}
 =e^{\tau\Phi}f+(-1)^NQ_2(1)Z_Ne^{\tau\Phi}.\nonumber
\end{eqnarray}
Repeating an argument similar to (\ref{novel}), we obtain
\begin{eqnarray}\label{novell}
\Vert Z_j\Vert_{L^\frac{p\widetilde p}{p-\widetilde p}(\Omega)}
\le \frac{C}{\tau^\delta}\Vert P_{A_1}(Q_2(1)Z_{j-1})\Vert
_{W^1_{\widetilde p}(\Omega)}\le \frac{C}{\tau^\delta}\Vert P_{A_1}
\Vert_{\mathcal L(L^{\widetilde p}(\Omega), W^1_{\widetilde p}(\Omega))}
\Vert Q_2(1)\Vert_{L^p(\Omega)}\Vert Z_{j-1}\Vert_{L^\frac{p\widetilde p}
{p-\widetilde p}(\Omega)}\nonumber\\
\le (\frac{C}{\tau^\delta}\Vert P_{A_1}\Vert
_{\mathcal L(L^{\widetilde p}(\Omega), W^1_{\widetilde p}(\Omega))}
\Vert Q_2(1)\Vert_{L^p(\Omega)})^{j-1}\Vert  Z_1\Vert
_{L^\frac{p\widetilde p}{p-\widetilde p}(\Omega)}
\end{eqnarray}
and the inequality corresponding to (\ref{novel1}):
\begin{eqnarray}\label{novel11}
\Vert Z_j\Vert_{L^\frac{p\widetilde p}{p-\widetilde p}(\Omega)}
\le \frac{1}{2^j}\Vert Z_1\Vert
_{L^\frac{p\widetilde p}{p-\widetilde p}(\Omega)}
\quad \forall j\ge 2\quad\mbox{and} \quad \forall \tau\ge \tau_0.
\end{eqnarray}
Therefore the infinite series (\ref{mursilkaa}) is convergent in
$L^2(\Omega).$   Passing to the limit in (\ref{gnomik}) as
$N\rightarrow +\infty$ we show that the function $Z$ is determined
by the infinite series (\ref{mursilkaa}) and is the solution to
the problem (\ref{-101}). By (\ref{novel}) we have
\begin{equation}
\Vert e^{-\tau\varphi} Z\Vert_{L^2(\Omega)}
\le C\Vert fe^{-\tau\varphi}\Vert_{L^2(\Omega)}.
\end{equation}
Let us show that the infinite series series defined by (\ref{mursilkaa}) is
convergent in the space $W^2_2(\Omega)$.
Applying  Proposition \ref{Proposition 3.0} we have
\begin{equation}\label{norton}
\Vert Z_j\Vert_{W^2_2(\Omega)}
\le C(\tau^2+1)\Vert Z_{j-2}\Vert_{L^2(\Omega)}\le C(\tau^2+1)\Vert Z
_{j-4}\Vert_{L^\frac{p\widetilde p}{p-\widetilde p}(\Omega)}.
\end{equation}
Let $j\ge 5.$  Using (\ref{novell}) to estimate the right-hand side of
(\ref{norton}) we obtain
\begin{eqnarray}
\Vert Z_j\Vert_{W^2_2(\Omega)}\le C(\tau^2+1)
(\frac{C}{\tau^\delta}\Vert P_{A_1}\Vert_{\mathcal L(L^{\widetilde p}(\Omega),
W^1_{\widetilde p}(\Omega))}
\Vert Q_2(1)\Vert_{L^p(\Omega)})^{j-5}\Vert  Z_1\Vert
_{L^\frac{p\widetilde p}{p-\widetilde p}(\Omega)}\nonumber\\
\le C(\tau^2+1)
(\frac{C}{\tau^\delta}\Vert P_{A_1}\Vert_{\mathcal L(L^{\widetilde p}(\Omega),
W^1_{\widetilde p}(\Omega))}
\Vert Q_2(1)\Vert_{L^p(\Omega)})^{j-5}\Vert  e^{-\tau\varphi}f\Vert
_{L^2(\Omega)}.
\end{eqnarray}
Hence for all sufficiently large $\tau$ there exists $\hat j$ such that
$$
\Vert Z_j\Vert_{W^2_2(\Omega)}\le (\tau^2+1)\left(\frac 12\right)^j\Vert
e^{-\tau\varphi}f\Vert_{L^2(\Omega)}.
$$
The proof of the estimate (\ref{zopa}) is complete.
$\blacksquare$

\section{Derivation of differential equations with respect to
unknown coefficients}

Let ${\bf H}(x,\partial_z,\partial_{\overline z})$ be a second-order
differential operator whose coefficients are smooth and have compact supports
in $\Omega.$
We write such an operator in the form
\begin{equation}\label{poker}
{\bf H}(x,\partial_z,\partial_{\overline z})
=\mathcal C_1(x)\partial^2_{zz}+\mathcal C_0(x)\partial^2_{z\overline z}
+\mathcal C_2(x)\partial^2_{\overline z\overline z}+{\bf B}_1(x)\partial_z
+{\bf B}_2(x)\partial_{\overline z} +{\bf B_0}(x),
\end{equation}
where $\mathcal C_j, {\bf B}_j$ are complex $3\times  3$ matrices.

Let the functions $U,V$ be solutions to equation (\ref{-10}) constructed
in Proposition \ref{mursilka0} by formulae (\ref{mursilka}) and
(\ref{mursilka1}) respectively and the functions
$\widetilde U,\widetilde V$ be solutions to equation (\ref{-10})
constructed in Proposition \ref{mursilka00} by formulae (\ref{mursilka++}) and (\ref{mursilka+++}) respectively. Denote
\begin{equation}\label{mishka}
q_2=(T_{B_2}(Q_2(2)V_0)-\mathcal D_2),\quad  q_1
=(P_{A_1}(Q_1(1)U_0)-\mathcal A_1)	
\end{equation}
and
\begin{equation}\label{mishka66}
\widetilde q_2=(T_{B_2}(Q_2(2)\widetilde V_0)-\overline z^k\mathcal D_2),
\quad  \widetilde q_1=(P_{A_1}(Q_1(1)\widetilde U_0)-z^k\mathcal A_1).
\end{equation}
We introduce the functional
$$
\frak I_\tau r=\int_{\partial\Omega} r\frac{(\nu_1-i\nu_2)}
{2\tau \partial_z\Phi}e^{\tau(\Phi-\overline \Phi)}d\sigma
-\int_{\partial\Omega} (\nu_1-i\nu_2)\partial_z
(\frac{r}{2\tau^2 \partial_z\Phi})e^{\tau(\Phi-\overline \Phi)} d\sigma.
$$
We recall that $\Phi(z) = (z-\widetilde{z})^2 - \kappa$ and
$\widetilde{\Phi}(z) = e^{i\theta}z$ with $\theta \in [0,2\pi)$.
We have

\begin{proposition}\label{zanoza1}
Let the functions $U_j,V_j$ be defined in Proposition \ref{mursilka0}
and $q_j$ given by (\ref{mishka}) and $q_1(\widetilde x)=q_2(\widetilde x)=0.$
The following asymptotic formulae hold true:
\begin{equation}
(U_0e^{\tau\Phi},{\bf H}(x,\partial_z,\partial_{\overline z})
(V_0e^{-\tau\overline \Phi}))_{L^2(\Omega)}
=\frak F_{\widetilde x,\tau}(U_0, {\bf H}(x,\partial_z,\partial_{\overline z}
-\tau{\partial_{\overline z}\overline\Phi})V_0)+o\left(\frac 1\tau\right)
\quad\mbox{as}\,\,\tau\rightarrow +\infty,
\end{equation}
\begin{eqnarray}
(U_0e^{\tau\Phi}, {\bf H}(x,\partial_z,\partial_{\overline z})
(V_1e^{-\tau\overline{ \Phi}}))_{L^2(\Omega)}=\frak F_{\widetilde x,\tau}
(P_{A_2}^*{\bf H}(x,\partial_z+\tau{\partial_z\Phi},\partial_{\overline z})^*
U_0, q_2)
                                         \nonumber\\
+\frak I_\tau(P_{A_2}^*{\bf H}(x,\partial_z+\tau{\partial_z\Phi},
\partial_{\overline z})^*U_0, q_2)+o\left(\frac {1}{\tau}\right)
\quad\mbox{as}\,\,\tau\rightarrow +\infty,
\end{eqnarray}

\begin{eqnarray}
(U_1e^{\tau\Phi},{\bf H}(x,\partial_z,\partial_{\overline z})
(V_0e^{-\tau\overline \Phi}))_{L^2(\Omega)}
= \frak F_{\widetilde x,\tau}(q_1, T_{B_1}^*
{\bf H}(x,\partial_z,\partial_{\overline z}-\tau{\partial_{\bar z}\bar \Phi})V_0)
                                          \nonumber \\
+\frak I_\tau (q_1, T_{B_1}^*{\bf H}(x,\partial_z,\partial_{\overline z}
-\tau{\partial_{\bar z}\bar \Phi})V_0)+o\left(\frac {1}{\tau}\right)
                      \quad\mbox{as}\,\,\tau\rightarrow +\infty,
\end{eqnarray}

\begin{eqnarray}\label{slon}
(U_2 e^{\tau\Phi},{\bf  H}(x,\partial_z,\partial_{\overline z})
(V_0 e^{-\tau\overline\Phi}))_{L^2(\Omega)}
= -\frak F_{\widetilde x,\tau}(q_1,T_{B_1}^*Q_2(1)^*
T_{B_1}^*({\mathcal C}_0\partial_z+{\mathcal C}_2(\partial_{\overline z}
-\tau{\partial_{\bar z}\bar \Phi})+b_1(x))V_0)                      \nonumber\\
-\frak I_\tau(q_1,T_{B_1}^*Q_2(1)^* T_{B_1}^*({\mathcal C}_0\partial_z
+{\mathcal C}_2(\partial_{\overline z}-\tau{\partial_{\bar z}\overline \Phi})+b)V_0)
+ o\left(\frac 1\tau\right)\quad\mbox{as}\,\,\tau\rightarrow +\infty,
\end{eqnarray}
where $b={\bf B_2}-\partial_{\overline z}\mathcal C_2+A_1^*\mathcal C_2;$
\begin{eqnarray}
(U_0 e^{\tau\Phi},{\bf  H}(x,\partial_z,\partial_{\overline z})
(V_2 e^{-\tau\overline\Phi}))_{L^2(\Omega)}
= -\frak F_{\widetilde x,\tau}((P_{A_2}^*Q_1(2)^*P_{A_2}^*({\mathcal C}_1^*
(\partial_{z}+\tau{\partial_z\Phi})+{\mathcal C}_0^*\partial_{\overline z}
+ \widetilde b)U_0 ,q_2)\nonumber\\
+ \frak I_\tau ((P_{A_2}^*Q_1(2)^*P_{A_2}^*({\mathcal C}_1^*(\partial_{z}
+ \tau{\partial_z\Phi})+{\mathcal C}_0^*\partial_{\overline z}  + \widetilde b)U_0 ,q_2)
+o\left(\frac 1\tau\right)\quad\mbox{as}\,\,\tau\rightarrow +\infty,
\end{eqnarray}
where $\widetilde b={\bf B}_1+\mathcal C_1B_2-\partial_z{\mathcal C}_1;$
\begin{eqnarray}
(U_1 e^{\tau\Phi}, {\bf H}(x,\partial_z,\partial_{\overline z})
(V_1 e^{-\tau\overline\Phi}))_{L^2(\Omega)}=
\frak F_{\widetilde x, \tau}(P_{A_2}^*({\mathcal C}_1^*
(\partial_z+\tau {\partial_z\Phi}+B_1)-{\mathcal B}_1^*) q_1) , q_2)\\
- \frak F_{\widetilde x, \tau}(q_1 ,{\mathcal C}_0 q_2 )
+ \frak F_{\widetilde x, \tau}(q_1 ,T_{B_1}^*( {\mathcal C}_2
(\partial_{\overline z}-\tau\partial_{\bar z}\overline\Phi+A_2) +{\mathcal B}_2)q_2 ))                                                     \nonumber\\
+ \frak I_\tau(P_{A_2}^*({\mathcal C}_1^*(\partial_z+\tau {\partial_z\Phi}+B_1)
- {\mathcal B}_1^*) q_1) , q_2)\nonumber\\
- \frak I_\tau (q_1 ,{\mathcal C}_0 q_2  )+\frak I_\tau (q_1 ,T_{B_1}^*((
{\mathcal C}_2(\partial_{\overline z}-\tau\partial_{\bar z}\overline\Phi+A_2) +{\mathcal B}_2)
q_2 ))+o\left(\frac 1\tau\right)\quad\mbox{as}\,\,\tau\rightarrow +\infty,\nonumber
\end{eqnarray}
where $\mathcal B_1=-2\partial_{z}\mathcal C_1 +2B_1^*\mathcal C_1
-\mathcal C_0A_2+{\bf B_1}$ and
$\mathcal B_2=B_1^*\mathcal C_0-\partial_z{\mathcal C}_0
-2{\mathcal C}_2A_2+{\bf B}_2.$
\end{proposition}

{\bf Proof.}
Since the coefficients of the operator ${\bf H}$ have compact supports,
by Proposition \ref{osel}, we have
\begin{eqnarray}
(U_0e^{\tau\Phi}, {\bf H}(x,\partial_z,\partial_{\overline z})
(V_0e^{-\tau\overline \Phi}))_{L^2(\Omega)}
= (U_0e^{\tau\Phi-\tau\overline \Phi},
{\bf H}(x,\partial_z,\partial_{\overline z}
-\tau{\partial_{\bar z}\bar \Phi})V_0)_{L^2(\Omega)}             \nonumber\\
= \frak F_{\widetilde x,\tau}(U_0, {\bf H}(x,\partial_z,\partial
_{\overline z}-\tau{\partial_{\bar z}\overline\Phi})V_0)
+o\left(\frac 1\tau\right)
\quad\mbox{as}\,\,\tau\rightarrow +\infty.
\end{eqnarray}

Using the definitions of the functions  $U_1$ and $V_1$  given in formulae
(\ref{mursilka}) and (\ref{mursilka1}) respectively and applying
the stationary phase argument, we have

\begin{eqnarray}
(U_1e^{\tau\Phi}, {\bf H}(x,\partial_z,\partial_{\overline z})
(V_0e^{-\tau\overline \Phi}))_{L^2(\Omega)}=(U_1e^{\tau(\Phi-\overline \Phi)},
{\bf H}(x,\partial_z,\partial_{\overline z}-\tau{\partial_{\bar z}\bar  \Phi})V_0)
_{L^2(\Omega)}                           \nonumber\\
= (e^{\tau(\Phi-\overline \Phi)}q_1, T_{B_1}^*{\bf H}
(x,\partial_z,\partial_{\overline z}-\tau{\partial_{\bar z}\bar  \Phi})V_0)_{L^2(\Omega)}
= \frak F_{\widetilde x,\tau}(q_1, T_{B_1}^*{\bf H}
(x,\partial_z,\partial_{\overline z}-\tau{\partial_{\bar z}\bar  \Phi})V_0) \nonumber\\
+ \frak I_\tau(q_1, T_{B_1}^*{\bf H}(x,\partial_z,\partial_{\overline z}
-\tau{\partial_{\bar z}\bar  \Phi})V_0)+o\left(\frac {1}{\tau}\right)
\quad\mbox{as}\,\,\tau\rightarrow +\infty,
\end{eqnarray}
and
\begin{eqnarray}
(U_0e^{\tau\Phi}, {\bf H}(x,\partial_z,\partial_{\overline z})
(V_1e^{-\tau\overline \Phi}))_{L^2(\Omega)}
= ({\bf H}(x,\partial_z,\partial_{\overline z})^*U_0e^{\tau\Phi},
V_1e^{-\tau\overline \Phi})_{L^2(\Omega)}\\
= ({\bf H}(x,\partial_z
+ \tau{\partial_z\Phi},\partial_{\overline z})^*U_0, V_1e^{\tau(\Phi-\overline \Phi)})
_{L^2(\Omega)}
= (P_{A_2}^*{\bf H}(x,\partial_z+\tau{\partial_z\Phi},\partial_{\overline z})^*U_0, q_2)
_{L^2(\Omega)}                             \nonumber\\
= \frak F_{\widetilde x,\tau} (P_{A_2}^*
{\bf H}(x,\partial_z+\tau{\partial_z\Phi},\partial_{\overline z})^*U_0, q_2)
                                         \nonumber\\
+ \frak I_\tau (P_{A_2}^*{\bf H}(x,\partial_z+\tau{\partial_z\Phi},
\partial_{\overline z})^*U_0, q_2)+o\left(\frac {1}{\tau}\right)
\quad\mbox{as}\,\,\tau\rightarrow +\infty.\nonumber
\end{eqnarray}

In order to prove the asymptotics (\ref{slon}) it is convenient to
represent the operator ${\bf H}$ in the form
$$
{\bf H}(x,\partial_z,\partial_{\overline z})
= {\mathcal C}_1(x)\partial_{zz}^2+(\partial_{\overline z}-A_1^*)
{\mathcal C}_0(x)\partial_z +(\partial_{\overline z}-A_1^*){\mathcal C}_2(x)
\partial_{\overline z} +B(x,D),
$$
where $B(x,D)=(\partial_{\overline z}-A_1^*)b(x)+b_1(x)\partial_z+b_0(x)$
and $b_0, b_1$ are some smooth functions with compact supports.
Using this representation and Proposition \ref{elkazelenaja} we obtain
\begin{eqnarray}
(U_2 e^{\tau\Phi},{\bf  H}(x,\partial_z,\partial_{\overline z})
(V_0 e^{-\tau\overline\Phi}))_{L^2(\Omega)}\\
= (U_2 e^{\tau\Phi},
({\mathcal C}_1\partial_{zz}^2+(\partial_{\overline z}-A_1^*){\mathcal C}_0
\partial_z +(\partial_{\overline z}-A_1^*){\mathcal C}_2\partial_{\overline z}
+ B(x,D)) (V_0 e^{-\tau\overline\Phi}))_{L^2(\Omega)}    \nonumber\\
= (U_2 e^{\tau\Phi},  ((\partial_{\overline z}-A_1^*){\mathcal C}_0\partial_z
+(\partial_{\overline z}-A_1^*){\mathcal C}_2\partial_{\overline z}
+(\partial_{\overline z}-A_1^*)b_1) (V_0 e^{-\tau\overline\Phi}))_{L^2(\Omega)}
+ o\left(\frac 1\tau\right).\nonumber
\end{eqnarray}

Then using  Proposition \ref{osel}, we have
\begin{eqnarray}
(U_2 e^{\tau\Phi},{\bf  H}(x,\partial_z,\partial_{\overline z})
(V_0 e^{-\tau\overline\Phi}))_{L^2(\Omega)}              \nonumber\\
= -((\partial_{\overline z}+A_1)(U_2 e^{\tau\Phi}),{\mathcal C}_0\partial_z
+{\mathcal C}_2\partial_{\overline z}+b(V_0 e^{-\tau\overline\Phi}))_{L^2(\Omega)}+o\left(\frac 1\tau\right)\\
= -(T_{B_1}(Q_2(1) e^{\tau(\Phi-\overline\Phi)}U_1),({\mathcal C}_0
\partial_z+{\mathcal C}_2(\partial_{\overline z}-\tau{\partial_{\bar z}\bar  \Phi})+b)
V_0 )_{L^2(\Omega)}+o\left(\frac 1\tau\right)            \nonumber\\
= -(T_{B_1}(Q_2(1) T_{B_1}(e^{\tau(\Phi-\overline\Phi)}q_1)),
({\mathcal C}_0\partial_z+{\mathcal C}_2(\partial_{\overline z}
- \tau{\partial_{\bar z}\bar  \Phi})+b)V_0 )_{L^2(\Omega)}+o\left(\frac 1\tau\right)
                                                           \nonumber\\
= -(e^{\tau(\Phi-\overline\Phi)}q_1,T_{B_1}^*Q_2(1)^* T_{B_1}^*({\mathcal C}_0
\partial_z+{\mathcal C}_2(\partial_{\overline z}-\tau{\partial_{\bar z}\bar  \Phi})+b)V_0 )
_{L^2(\Omega)}+o\left(\frac 1\tau\right)                  \nonumber\\
= -\frak F_{\widetilde x,\tau}(q_1,T_{B_1}^*Q_2(1)^* T_{B_1}^*(({\mathcal C}_0
\partial_z+{\mathcal C}_2(\partial_{\overline z}-\tau{\partial_{\bar z}\bar  \Phi})+b)V_0))
                                                      \nonumber\\
- \frak I_\tau(q_1,T_{B_1}^*Q_2(1)^* T_{B_1}^*(({\mathcal C}_0
\partial_z+{\mathcal C}_2(\partial_{\overline z}-\tau{\partial_{\bar z}\bar  \Phi})+b)V_0)
+o\left(\frac 1\tau\right)\nonumber
\end{eqnarray}
and
\begin{eqnarray}
(U_0 e^{\tau\Phi},{\bf  H}(x,\partial_z,\partial_{\overline z})
(V_2 e^{-\tau\overline\Phi}))_{L^2(\Omega)}
= ({\bf  H}^*(x,\partial_z,\partial_{\overline z})
(U_0 e^{\tau\Phi}), V_2 e^{-\tau\overline\Phi})_{L^2(\Omega)}    \nonumber\\
= (((\partial_z-B_2^*){\mathcal C}_1^*\partial_{z}+(\partial_{ z}-B_2^*)
{\mathcal C}_0^*\partial_{\overline z} +{\mathcal C}_2^*(x)
\partial_{\overline z}^2 +(\partial_z-B_2^*) \widetilde b)(U_0 e^{\tau\Phi}),
V_2 e^{-\tau\overline\Phi})_{L^2(\Omega)}+o\left(\frac 1\tau\right)
                                                                \nonumber\\
= -(({\mathcal C}_1^*\partial_{z}+{\mathcal C}_0^*\partial_{\overline z}
+ \widetilde b)(U_0 e^{\tau\Phi}),e^{-\tau\Phi}(P_{A_2}(Q_1(2)V_1)
e^{\tau(\Phi-\overline\Phi)}))_{L^2(\Omega)}+o\left(\frac 1\tau\right)
                                                                  \nonumber\\
= -(({\mathcal C}_1^*(\partial_{z}+\tau{\partial_z\Phi})
+ {\mathcal C}_0^*\partial_{\overline z}  + \widetilde b)U_0 ,P_{A_2}
(Q_1(2)V_1e^{\tau(\Phi-\overline\Phi)}))_{L^2(\Omega)}
+o\left(\frac 1\tau\right)                            \nonumber\\
= -(({\mathcal C}_1^*(\partial_{z}+\tau{\partial_z\Phi})+{\mathcal C}_0^*
\partial_{\overline z} + \widetilde b)U_0 ,P_{A_2}(Q_1(2)P_{A_2}
(q_2e^{\tau(\Phi-\overline\Phi)})))_{L^2(\Omega)}+o\left(\frac 1\tau\right)
                                                          \nonumber\\
= -(P_{A_2}^*Q_1(2)^*P_{A_2}^*({\mathcal C}_1^*(\partial_{z}+\tau{\partial_z\Phi})
+{\mathcal C}_0^*\partial_{\overline z}  + \widetilde b)U_0 ,q_2
e^{\tau(\Phi-\overline\Phi)})_{L^2(\Omega)}+o\left(\frac 1\tau\right)
                                         \nonumber\\
= -\frak F_{\widetilde x,\tau}((P_{A_2}^*Q_1(2)^*P_{A_2}^*({\mathcal C}_1^*
(\partial_{z}+\tau{\partial_z\Phi})+{\mathcal C}^*_0\partial_{\overline z}
+ \widetilde b)U_0 ,q_2)                                       \nonumber\\
- \frak I_\tau(P_{A_2}^*Q_1(2)^*P_{A_2}^*({\mathcal C}_1^*(\partial_{z}
+\tau{\partial_z\Phi})+{\mathcal C}^*_0\partial_{\overline z}  + \widetilde b)U_0 ,q_2)
+o\left(\frac 1\tau\right).\nonumber
\end{eqnarray}
Integrating by parts, using Propositions \ref{elkazelenaja} and \ref{osel},
we obtain
\begin{eqnarray}
(U_1 e^{\tau\Phi}, {\bf H}(x,\partial_z,\partial_{\overline z})
(V_1 e^{-\tau\overline\Phi}))_{L^2(\Omega)}=(U_1 e^{\tau\Phi},
((-\partial_z+B_1^*)^2{\mathcal C}_1+(\partial_z-B_1^*){\mathcal C}_0
(\partial_{\overline z}+A_2)                               \nonumber\\
+ {\mathcal C}_2(\partial_{\overline z}+A_2)^2+(\partial_z-B_1^*){\mathcal B}_1
+ {\mathcal B}_2(\partial_{\overline z}+A_2))(V_1 e^{-\tau\overline\Phi}))_{L^2(\Omega)}+o\left(\frac 1\tau\right)                           \nonumber\\
= ((\partial_z+B_1)^2U_1 e^{\tau\Phi}, {\mathcal C}_1V_1
e^{-\tau\overline\Phi})_{L^2(\Omega)}-((\partial_z+B_1)U_1 e^{\tau\Phi},
{\mathcal C}_0(\partial_{\overline z}+A_2)(V_1 e^{-\tau\overline\Phi}))
_{L^2(\Omega)}                                              \nonumber\\
+ (U_1 e^{\tau\Phi}, ({\mathcal C}_2(\partial_{\overline z}+A_2)^2
+{\mathcal B}_2(\partial_{\overline z}+A_2))(V_1 e^{-\tau\overline\Phi}))
_{L^2(\Omega)}                                             \nonumber\\
-((\partial_z+B_1)(U_1 e^{\tau\Phi}),{\mathcal B}_1V_1
e^{-\tau\overline\Phi})_{L^2(\Omega)}+o\left(\frac 1\tau\right)   \nonumber\\
= ((\partial_z+B_1) q_1 e^{\tau\Phi}, {\mathcal C}_1V_1
e^{-\tau\overline\Phi})_{L^2(\Omega)}-(q_1 e^{\tau\Phi},{\mathcal C}_0
q_2 e^{-\tau\overline\Phi})_{L^2(\Omega)}                     \nonumber\\
+ (U_1 e^{\tau\Phi}, ({\mathcal C}_2(\partial_{\overline z}+A_2)
+ {\mathcal B}_2)(q_2 e^{-\tau\overline\Phi}))_{L^2(\Omega)} \nonumber\\
-( q_1 e^{\tau\Phi},{\mathcal B}_1V_1 e^{-\tau\overline\Phi})_{L^2(\Omega)}
+ o\left(\frac 1\tau\right)                \nonumber\\
= ((\partial_z+\tau {\partial_z\Phi}+B_1) q_1 , {\mathcal C}_1V_1
e^{\tau(\Phi-\overline\Phi)})_{L^2(\Omega)}-(q_1 e^{\tau\Phi},{\mathcal C}_0
q_2 e^{-\tau\overline\Phi})_{L^2(\Omega)}                      \nonumber\\
+ (U_1 e^{\tau\Phi-\tau\overline\Phi}, ({\mathcal C}_2(\partial_{\overline z}
-\tau\partial_{\bar z}{\overline \Phi}+A_2) +{\mathcal B}_2) q_2)_{L^2(\Omega)}       \nonumber\\
-( q_1 ,{\mathcal B}_1V_1 e^{\tau(\Phi-\overline\Phi)})_{L^2(\Omega)}
                                                                 \nonumber\\
= ((\partial_z+\tau {\partial_z\Phi}+B_1) q_1 , {\mathcal C}_1P_{A_2}
(q_2 e^{\tau(\Phi-\overline\Phi)}))_{L^2(\Omega)}               \nonumber\\
-(q_1 e^{\tau(\Phi-\overline\Phi)},{\mathcal C}_0 q_2 )_{L^2(\Omega)}
                                                                  \nonumber\\
+(T_{B_1}(q_1 e^{\tau(\Phi-\overline\Phi)}), ({\mathcal C}_2
(\partial_{\overline z}-\tau\partial_{\bar z}{\overline \Phi}+A_2) +{\mathcal B}_2) q_2)
_{L^2(\Omega)}                                  \nonumber\\
-( q_1 ,{\mathcal B}_1V_1 e^{\tau(\Phi-\overline\Phi)})_{L^2(\Omega)}+o\left(\frac 1\tau\right)                                        \nonumber\\
= (P_{A_2}^*{\mathcal C}_1^*(\partial_z+\tau {\partial_z\Phi}+B_1) q_1 , q_2
e^{\tau(\Phi-\overline\Phi )})_{L^2(\Omega)}                 \nonumber\\
-(q_1 e^{\tau(\Phi-\overline\Phi)},{\mathcal C}_0q_2)_{L^2(\Omega)}\nonumber\\
+(q_1 e^{\tau(\Phi-\overline\Phi)},T_{B_1}^*( {\mathcal C}_2(\partial
_{\overline z}-\tau\partial_{\bar z}\overline \Phi +A_2) +{\mathcal B}_2) q_2 ))_{L^2(\Omega)}
                                                               \nonumber\\
-( q_1,{\mathcal B}_1P_{A_2}(q_2e^{\tau(\Phi-\overline\Phi)}))_{L^2(\Omega)}
                                                                 \nonumber\\
= \frak F_{\widetilde x, \tau}(P_{A_2}^*({\mathcal C}_1^*
(\partial_z+\tau {\partial_z\Phi}+B_1) q_1) , q_2)                  \nonumber\\
- \frak F_{\widetilde x, \tau}(q_1 ,{\mathcal C}_0 q_2 )\nonumber\\
+ \frak F_{\widetilde x, \tau}(q_1 ,T_{B_1}^*(( {\mathcal C}_2
(\partial_{\overline z}-\tau\partial_{\bar z}{\overline \Phi}+A_2) +{\mathcal B}_2)q_2))
                                                            \nonumber\\
-\frak F_{\widetilde x, \tau}(P_{A_2}^*{\mathcal B}_1^* q_1 ,q_2 )\nonumber\\
+\frak I_\tau(P_{A_2}^*{\mathcal C}_1^*(\partial_z+\tau {\partial_z\Phi}+B_1) q_1 , q_2)
                                          \nonumber\\
-\frak I_\tau(q_1 ,{\mathcal C}_0 q_2 )\nonumber\\
+\frak I_\tau(q_1 ,T_{B_1}^*(( {\mathcal C}_2(\partial_{\overline z}
-\tau\partial_{\bar z}{\overline \Phi}+A_2) +{\mathcal B}_2)q_2))\nonumber\\
-\frak I_\tau(P_{A_2}^*{\mathcal B}_1^* q_1 ,q_2 )
+ o\left(\frac 1\tau\right)\quad \mbox{as}\quad \tau\rightarrow +\infty.
\end{eqnarray}
The proof of the proposition is complete.
$\blacksquare$

Similarly to Proposition \ref{zanoza1}, we prove Proposition \ref{zanoza2}.

\begin{proposition}\label{zanoza2}
Let the functions $\widetilde U_j,
\widetilde V_j$ be defined in  Proposition \ref{mursilka0} and
$\widetilde q_j$ be given by (\ref{mishka66}).
The following asymptotic formulae hold true:
\begin{equation}
(\widetilde U_0e^{\tau\widetilde\Phi},
{\bf H}(x,\partial_z,\partial_{\overline z})
(\widetilde V_0e^{-\tau\overline {\widetilde \Phi}}))_{L^2(\Omega)}
= o\left(\frac 1\tau\right)\quad\mbox{as}\,\,\tau\rightarrow+\infty,
\end{equation}
\begin{equation}
(\widetilde U_0e^{\tau\widetilde \Phi}, {\bf H}
(x,\partial_z,\partial_{\overline z})(\widetilde V_1
e^{-\tau\overline{\widetilde \Phi}}))_{L^2(\Omega)}
=\frak I_\tau(P_{A_2}^*{\bf H}(x,\partial_z+\tau\partial_{z}\widetilde \Phi,
\partial_{\overline z})^*\widetilde U_0, \widetilde q_2)+o(\frac {1}{\tau}),
\end{equation}
\begin{equation}
(\widetilde U_1e^{\tau\widetilde\Phi},{\bf H}(x,\partial_z,
\partial_{\overline z})(\widetilde V_0e^{-\tau\overline {\widetilde \Phi}}))
_{L^2(\Omega)}
= \frak I_\tau (\widetilde q_1, T_{B_1}^*{\bf H}
(x,\partial_z,\partial_{\overline z}-\tau\partial_{\overline z}\overline{\widetilde \Phi})
\widetilde V_0)+o(\frac {1}{\tau}),
\end{equation}
\begin{equation}\label{zem5}
(\widetilde U_2 e^{\tau{\widetilde\Phi}},
{\bf  H}(x,\partial_z,\partial_{\overline z}) (\widetilde V_0
e^{-\tau\overline{\widetilde\Phi}}))_{L^2(\Omega)}
= -\frak I_\tau(\widetilde q_1,T_{B_1}^*Q_2(1)^* T_{B_1}^*({\mathcal C}_0
\partial_z+{\mathcal C}_2(\partial_{\overline z}
-\tau\partial_{\overline z}\overline{\widetilde \Phi})+b)V_0)+o\left(\frac 1\tau\right),
\end{equation}
\begin{equation}
(\widetilde U_0 e^{\tau\widetilde \Phi},
{\bf  H}(x,\partial_z,\partial_{\overline z}) (\widetilde V_2
e^{-\tau\overline{\widetilde \Phi}}))_{L^2(\Omega)}
= -\frak I_\tau ((P_{A_2}^*Q_1(2)^*P_{A_2}^*({\mathcal C}_1^*
(\partial_{z}+\tau\partial_{ z}\widetilde\Phi)+{\mathcal C}_0^*(x)\partial_{\overline z}
+ \widetilde b)U_0 ,\widetilde q_2)+o\left(\frac 1\tau\right),
\end{equation}
\begin{eqnarray}
(\widetilde U_1 e^{\tau\widetilde\Phi}, {\bf H}(x,\partial_z,\partial
_{\overline z}) (\widetilde V_1 e^{-\tau\overline{\widetilde\Phi}}))
_{L^2(\Omega)}
= \frak I_\tau(P_{A_2}^*({\mathcal C}_1^*(\partial_z+\tau\partial_{ z} \widetilde\Phi+B_1)
-{\mathcal B}_1^*) \widetilde q_1 ,\widetilde q_2)           \nonumber\\
- \frak I_\tau (\widetilde q_1 ,{\mathcal C}_0\widetilde q_2)
+ \frak I_\tau (\widetilde q_1 ,T_{B_1}^*(( {\mathcal C}_2
(\partial_{\overline z}-\tau\partial_{\overline z}\overline{\widetilde \Phi}+A_2)
+ {\mathcal B}_2)\widetilde q_2 ))+o\left(\frac 1\tau\right),
\end{eqnarray}
as $\tau \to +\infty$.
Here we recall (see Proposition \ref{zanoza1}) that
$b={\bf B_2}-\partial_{\overline z}\mathcal C_2+A_1^*\mathcal C_2$,
$\widetilde b={\bf B}_1+\mathcal C_1B_2-\partial_z{\mathcal C}_1$,
$\mathcal B_1=-2\partial_{z}\mathcal C_1 +2B_1^*\mathcal C_1-\mathcal C_0A_2+{\bf B}_1$
and $\mathcal B_2=B_1^*\mathcal C_0-\partial_z{\mathcal C}_0
-2\mathcal C_2A_2+{\bf B}_2$.
\end{proposition}

{\bf Proof.} Since the coefficients of the operator ${\bf H}$ have compact
supports, from  Proposition \ref{osel} we have
\begin{eqnarray}
(\widetilde U_0e^{\tau{\widetilde \Phi}}, {\bf H}(x,\partial_z,\partial_{\overline z})
(\widetilde V_0e^{-\tau\overline {\widetilde \Phi}}))_{L^2(\Omega)}            \nonumber\\
= (\widetilde U_0e^{\tau({\widetilde \Phi}-\overline {\widetilde \Phi})},
{\bf H}(x,\partial_z,\partial_{\overline z}-\tau\partial_{\overline z}\overline{\widetilde \Phi})
\widetilde V_0)_{L^2(\Omega)}=o\left(\frac 1\tau\right)
\quad\mbox{as}\,\,\tau\rightarrow +\infty.
\end{eqnarray}

Applying the stationary phase argument, we obtain
\begin{eqnarray}
(\widetilde U_1e^{\tau{\widetilde \Phi}}, {\bf H}(x,\partial_z,\partial_{\overline z})
(\widetilde V_0e^{-\tau\overline {\widetilde \Phi}}))_{L^2(\Omega)}
= (\widetilde U_1e^{\tau({\widetilde \Phi}-\overline {\widetilde \Phi})},
{\bf H}(x,\partial_z,\partial_{\overline z}-\tau\partial_{\overline z}\overline{\widetilde \Phi})
\widetilde V_0)_{L^2(\Omega)}                                        \nonumber\\
= (e^{\tau({\widetilde \Phi}-\overline {\widetilde \Phi})}\widetilde q_1,
T_{B_1}^*{\bf H}(x,\partial_z,\partial_{\overline z}
- \tau\partial_{\bar z}\overline{\widetilde \Phi})\widetilde V_0)_{L^2(\Omega)}     \nonumber\\
= \frak I_\tau(\widetilde q_1, T_{B_1}^*{\bf H}(x,\partial_z,
\partial_{\overline z}-\tau\partial_{\overline z}\overline{\widetilde \Phi})\widetilde V_0)
+o(\frac {1}{\tau})\quad\mbox{as}\,\,\tau\rightarrow +\infty,
\end{eqnarray}
and
\begin{eqnarray}
(\widetilde U_0e^{\tau{\widetilde \Phi}}, {\bf H}(x,\partial_z,\partial_{\overline z})
(\widetilde V_1e^{-\tau\overline{\widetilde \Phi}}))_{L^2(\Omega)}=({\bf H}(x,\partial_z,\partial
_{\overline z})^* \widetilde U_0e^{\tau\widetilde\Phi},\widetilde V_1e^{-\tau\overline{\widetilde \Phi}})_{L^2(\Omega)}
                                                                 \nonumber\\
=({\bf H}(x,\partial_z+\tau{\partial_z\widetilde\Phi},\partial_{\overline z})^*\widetilde U_0, \widetilde V_1
e^{\tau(\Phi-\overline \Phi)})_{L^2(\Omega)}
= (P_{A_2}^*{\bf H}(x,\partial_z+\tau{\partial_z\widetilde\Phi},\partial_{\overline z})^*\widetilde U_0,
\widetilde q_2)_{L^2(\Omega)}                             \nonumber\\
= \frak I_\tau (P_{A_2}^*{\bf H}(x,\partial_z+\tau\partial_z\widetilde\Phi,
\partial_{\overline z})^*\widetilde U_0,\widetilde q_2)
+ o(\frac {1}{\tau})\quad\mbox{as}\,\,\tau\rightarrow +\infty.
\end{eqnarray}

In order to prove the asymptotics (\ref{zem5}) it is convenient to represent
the operator ${\bf H}$ in the form
$$
{\bf H}(x,\partial_z,\partial_{\overline z})={\mathcal C}_1(x)\partial_{zz}^2
+ (\partial_{\overline z}-A_1^*){\mathcal C}_0(x)\partial_z
+(\partial_{\overline z}-A_1^*){\mathcal C}_2(x)\partial_{\overline z}
+ B(x,D),
$$
where $B(x,D)=(\partial_{\overline z}-A_1^*)b(x)+b_1(x)\partial_z+b_0(x)$
is the same as in the proof of Proposition \ref{zanoza1}.
Then Proposition \ref{elkazelenaja} yields
\begin{eqnarray}
(\widetilde U_2 e^{\tau{\widetilde \Phi}},{\bf  H}(x,\partial_z,\partial_{\overline z})
(\widetilde V_0 e^{-\tau\overline{\widetilde \Phi}}))_{L^2(\Omega)}   \\
= (\widetilde U_2 e^{\tau\widetilde\Phi},  ({\mathcal C}_1\partial_{zz}^2+(\partial_{\overline z}
-A_1^*){\mathcal C}_0\partial_z +(\partial_{\overline z}-A_1^*){\mathcal C}_2
\partial_{\overline z} +B(x,D)) (\widetilde V_0 e^{-\tau\overline{\widetilde \Phi}}))
_{L^2(\Omega)}                                         \nonumber\\
= (\widetilde U_2 e^{\tau\widetilde\Phi},  (\partial_{\overline z}-A_1^*){\mathcal C}_0
\partial_z +(\partial_{\overline z}-A_1^*){\mathcal C}_2\partial
_{\overline z}+(\partial_{\overline z}-A_1^*)b)
(\widetilde V_0 e^{-\tau\overline{\widetilde\Phi}}))_{L^2(\Omega)}
+ o\left(\frac 1\tau\right).\nonumber
\end{eqnarray}

Then using Proposition \ref{osel}, we have
\begin{eqnarray}
(\widetilde U_2 e^{\tau{\widetilde \Phi}},{\bf  H}(x,\partial_z,\partial_{\overline z})
(\widetilde V_0 e^{-\tau\overline{\widetilde \Phi}}))_{L^2(\Omega)}           \\
= -((\partial_{\overline z}+A_1)(\widetilde U_2 e^{\tau{\widetilde \Phi}}),{\mathcal C}_0
\partial_z+{\mathcal C}_2\partial_{\overline z}
+ b(\widetilde V_0 e^{-\tau\overline{\widetilde \Phi}}))_{L^2(\Omega)}
+ o\left(\frac 1\tau\right)                                 \nonumber\\
= -(T_{B_1}(Q_2(1) e^{\tau(\widetilde\Phi-\overline{\widetilde \Phi})}\widetilde U_1),
({\mathcal C}_0\partial_z+{\mathcal C}_2(\partial_{\overline z}
- \tau\partial_{\bar z}\overline{\widetilde \Phi})+b)\widetilde V_0 )_{L^2(\Omega)}
+ o\left(\frac 1\tau\right)                                \nonumber\\
= -(T_{B_1}(Q_2(1) T_{B_1}(e^{\tau(\widetilde\Phi-\overline{\widetilde\Phi})}\widetilde q_1)),
({\mathcal C}_0\partial_z+{\mathcal C}_2(\partial_{\overline z}
- \tau\partial_{\bar z}\overline{\widetilde \Phi})+b)\widetilde V_0 )_{L^2(\Omega)}
+ o\left(\frac 1\tau\right)                            \nonumber\\
= -(e^{\tau({\widetilde \Phi}-\overline{\widetilde \Phi})}
\widetilde q_1,T_{B_1}^*Q_2(1)^* T_{B_1}^*({\mathcal C}_0\partial_z
+{\mathcal C}_2(\partial_{\overline z}-\tau\partial_{\bar z}\overline {\widetilde \Phi})+b)
\widetilde V_0 )_{L^2(\Omega)}+o\left(\frac 1\tau\right)
\nonumber\\
= -\frak I_\tau(\widetilde q_1,T_{B_1}^*Q_2(1)^* T_{B_1}^*({\mathcal C}_0
\partial_z+{\mathcal C}_2(\partial_{\overline z}
-\tau\partial_{\bar z}\overline{\widetilde \Phi})+b)\widetilde V_0)+o\left(\frac 1\tau\right)\nonumber
\end{eqnarray}
and
\begin{eqnarray}
(\widetilde U_0 e^{\tau{\widetilde \Phi}},{\bf  H}(x,\partial_z,\partial_{\overline z})
(\widetilde V_2 e^{-\tau\overline{\widetilde \Phi}}))_{L^2(\Omega)}
=({\bf  H}^*(x,\partial_z,\partial_{\overline z})
(\widetilde U_0 e^{\tau{\widetilde \Phi}}), \widetilde V_2 e^{-\tau\overline{\widetilde \Phi}})
_{L^2(\Omega)}                                                   \\
=(((\partial_z-B_2^*){\mathcal C}_1^*\partial_{z}+(\partial_{ z}-B_2^*)
{\mathcal C}_0^*\partial_{\overline z} +{\mathcal C}_2^*(x)\partial
_{\overline z\overline z}^2 +(\partial_z-B_2^*) \widetilde b)
(\widetilde U_0 e^{\tau{\widetilde \Phi}}),\widetilde V_2 e^{-\tau\overline{\widetilde \Phi}})
_{L^2(\Omega)}+o\left(\frac 1\tau\right)                           \nonumber\\
= -(({\mathcal C}_1^*\partial_{z}+{\mathcal C}_0^*\partial_{\overline z}
+ \widetilde b)(\widetilde U_0 e^{\tau\widetilde\Phi}),e^{-\tau{\widetilde \Phi}}
(P_{A_2}(Q_1(2)\widetilde V_1) e^{\tau({\widetilde \Phi}-\overline{\widetilde \Phi})}))
_{L^2(\Omega)}+o\left(\frac 1\tau\right)                       \nonumber\\
= -(({\mathcal C}_1^*(\partial_{z}+\tau\partial_z{\widetilde \Phi})+{\mathcal C}_0^*
\partial_{\overline z}  + \widetilde b)\widetilde U_0 ,P_{A_2}(Q_1(2)\widetilde V_1
e^{\tau({\widetilde \Phi}-\overline{\widetilde \Phi})}))_{L^2(\Omega)}
+o\left(\frac 1\tau\right)\nonumber\\
= -(({\mathcal C}_1^*(\partial_{z}+\tau\partial_z\widetilde\Phi)+{\mathcal C}_0^*
\partial_{\overline z} + \widetilde b)\widetilde U_0 ,P_{A_2}(Q_1(2)P_{A_2}
(\widetilde q_2e^{\tau({\widetilde \Phi}-\overline{\widetilde \Phi})})))
_{L^2(\Omega)}+o\left(\frac 1\tau\right)                     \nonumber\\
= -(P_{A_2}^*Q_1(2)^*P_{A_2}^*({\mathcal C}_1^*(\partial_{z}
+ \tau\partial_z{\widetilde \Phi})+{\mathcal C}_0^*\partial_{\overline z}
+ \widetilde b)\widetilde U_0 ,\widetilde q_2e^{\tau({\widetilde \Phi}
-\overline{\widetilde \Phi})})_{L^2(\Omega)}+o\left(\frac 1\tau\right)
                                                              \nonumber\\
= -\frak I_\tau((P_{A_2}^*Q_1(2)^*P_{A_2}^*({\mathcal C}_1^*(\partial_{z}
+\tau\partial_z{\widetilde \Phi})+{\mathcal C}^*_0\partial_{\overline z}
+ \widetilde b)\widetilde U_0 ,\widetilde q_2)+o\left(\frac 1\tau\right).\nonumber
\end{eqnarray}
Integrating by parts, using Propositions \ref{elkazelenaja} and \ref{osel},
we obtain
\begin{eqnarray}
(\widetilde U_1 e^{\tau{\widetilde \Phi}}, {\bf H}(x,\partial_z,\partial_{\overline z})
(\widetilde V_1 e^{-\tau\overline{\widetilde \Phi}}))_{L^2(\Omega)}
= (\widetilde U_1 e^{\tau{\widetilde\Phi}}, ((-\partial_z+B_1^*)^2{\mathcal C}_1
+(\partial_z-B_1^*){\mathcal C}_0(\partial_{\overline z}+A_2)\nonumber\\
+ {\mathcal C}_2(\partial_{\overline z}+A_2)^2+(\partial_z-B_1^*)
{\mathcal B}_1 +{\mathcal B}_2(\partial_{\overline z}+A_2))
(\widetilde V_1 e^{-\tau\overline{\widetilde \Phi}}))_{L^2(\Omega)}+o\left(\frac 1\tau\right)
                                                                \nonumber\\
=((\partial_z+B_1)^2\widetilde U_1 e^{\tau{\widetilde \Phi}}, {\mathcal C}_1\widetilde V_1
e^{-\tau\overline{\widetilde \Phi}})_{L^2(\Omega)}-((\partial_z+B_1)\widetilde U_1
e^{\tau{\widetilde \Phi}},{\mathcal C}_0(\partial_{\overline z}+A_2)
(\widetilde V_1 e^{-\tau\overline{\widetilde\Phi}}))_{L^2(\Omega)}                  \nonumber\\
+ (\widetilde U_1 e^{\tau{\widetilde \Phi}}, ({\mathcal C}_2
(\partial_{\overline z}+A_2)^2 +{\mathcal B}_2(\partial_{\overline z}+A_2))
(\widetilde V_1 e^{-\tau\overline{\widetilde \Phi}}))_{L^2(\Omega)}\nonumber\\
-((\partial_z+B_1)(\widetilde U_1 e^{\tau{\widetilde \Phi}}),{\mathcal B}_1\widetilde V_1
e^{-\tau\overline{\widetilde\Phi}})_{L^2(\Omega)}+o\left(\frac 1\tau\right)
                        \nonumber\\
= ((\partial_z+B_1) \widetilde q_1 e^{\tau{\widetilde \Phi}},
{\mathcal C}_1\widetilde V_1 e^{-\tau\overline{\widetilde \Phi}})_{L^2(\Omega)}
-(\widetilde q_1 e^{\tau{\widetilde \Phi}},{\mathcal C}_0
\widetilde q_2 e^{-\tau\overline{\widetilde \Phi}})_{L^2(\Omega)}   \nonumber\\
+ (\widetilde U_1 e^{\tau{\widetilde\Phi}}, ({\mathcal C}_2(\partial_{\overline z}+A_2)
+{\mathcal B}_2)(\widetilde q_2 e^{-\tau\overline{\widetilde\Phi}}))
_{L^2(\Omega)}                              \nonumber\\
-( \widetilde q_1 e^{\tau{\widetilde \Phi}},{\mathcal B}_1V_1
e^{-\tau\overline{\widetilde \Phi}})_{L^2(\Omega)}
+ o\left(\frac 1\tau\right)                            \nonumber\\
= ((\partial_z+\tau \partial_z{\widetilde \Phi}+B_1)\widetilde  q_1,
{\mathcal C}_1\widetilde V_1 e^{\tau({\widetilde \Phi}-\overline{\widetilde\Phi})})
_{L^2(\Omega)}-(\widetilde q_1 e^{\tau\widetilde\Phi},
{\mathcal C}_0 \widetilde q_2 e^{-\tau\overline{\widetilde \Phi}})
_{L^2(\Omega)}                                               \nonumber\\
+ (\widetilde U_1 e^{\tau(\widetilde\Phi-\overline{\widetilde\Phi})}, ({\mathcal C}_2(\partial_{\overline z}
-\tau\partial_{\bar z}\overline {\widetilde \Phi}+A_2) +{\mathcal B}_2)\widetilde q_2)
_{L^2(\Omega)}\nonumber\\
-(\widetilde q_1 ,{\mathcal B}_1\widetilde V_1 e^{\tau({\widetilde \Phi}
- \overline{\widetilde \Phi})})_{L^2(\Omega)} + o\left(\frac 1\tau\right)
                                                           \nonumber\\
= ((\partial_z+\tau \partial_z{\widetilde \Phi}+B_1)\widetilde q_1 ,
{\mathcal C}_1P_{A_2}(\widetilde q_2 e^{\tau({\widetilde \Phi}
-\overline{\widetilde \Phi})}))_{L^2(\Omega)}            \nonumber\\
- (\widetilde q_1e^{\tau({\widetilde \Phi}-\overline{\widetilde \Phi})},
{\mathcal C}_0
\widetilde q_2 )_{L^2(\Omega)}\nonumber\\
+ (T_{B_1}(\widetilde q_1 e^{\tau({\widetilde \Phi}
- \overline{\widetilde \Phi})}), ({\mathcal C}_2(\partial_{\overline z}
-\tau\partial_{\bar z}\overline {\widetilde \Phi}+A_2) +{\mathcal B}_2) \widetilde q_2)
_{L^2(\Omega)}                                    \nonumber\\
- (\widetilde q_1 ,{\mathcal B}_1\widetilde V_1 e^{\tau({\widetilde \Phi}
- \overline{\widetilde \Phi})})_{L^2(\Omega)}+o\left(\frac 1\tau\right)
                                                             \nonumber\\
= (P_{A_2}^*{\mathcal C}_1^*(\partial_z+\tau \partial_z{\widetilde \Phi}+B_1)
\widetilde q_1 , \widetilde q_2 e^{\tau({\widetilde \Phi}
-\overline{\widetilde \Phi} )})_{L^2(\Omega)}\nonumber\\-(\widetilde q_1
e^{\tau({\widetilde \Phi}-\overline{\widetilde \Phi})},{\mathcal C}_0
\widetilde q_2)_{L^2(\Omega)}\nonumber\\
+ (\widetilde q_1 e^{\tau({\widetilde \Phi}-\overline{\widetilde \Phi})},
T_{B_1}^*( {\mathcal C}_2(\partial_{\overline z}
-\tau\partial_{\bar z}\overline{\widetilde \Phi}+A_2) +{\mathcal B}_2)\widetilde q_2 ))
_{L^2(\Omega)}\nonumber\\
- (\widetilde q_1,{\mathcal B}_1P_{A_2}(\widetilde q_2
e^{\tau({\widetilde \Phi}-\overline{\widetilde \Phi})}))_{L^2(\Omega)}
                                                               \nonumber\\
= +\frak I_\tau(P_{A_2}^*{\mathcal C}_1^*(\partial_z+\tau \partial_z{\widetilde \Phi}
+ B_1)\widetilde q_1 ,\widetilde q_2)                    \nonumber\\
- \frak I_\tau(\widetilde q_1 ,{\mathcal C}_0\widetilde q_2 )\nonumber\\
+\frak I_\tau(\widetilde q_1 ,T_{B_1}^*( {\mathcal C}_2(\partial_{\overline z}
- \tau\partial_{\bar z}\overline {\widetilde \Phi}+A_2) +{\mathcal B}_2)\widetilde q_2))
                                                           \nonumber\\
- \frak I_\tau(P_{A_2}^*{\mathcal B}_1^* \widetilde q_1 ,\widetilde q_2 )
+o\left(\frac 1\tau\right)\quad\mbox{as}\,\,\tau\rightarrow +\infty.
\end{eqnarray}
The proof of the proposition is complete.
$\blacksquare$

We have
\begin{proposition}\label{zanoza}
Let $k, \ell\in\Bbb N_+$ and $k+\ell\ge 3$ and ${\bf H}(x,\partial_z,\partial_{\overline z})$ be a second-order differential operator with compactly supported  smooth coefficients.  Let the functions $U$ and $V$ be given by
Proposition \ref{mursilka0} and $q_1,q_2$ be given by (\ref{mishka}).
Assume that the restriction of the function $\psi$ on $\partial\Omega$ has
a finite number of critical points and all these points are nondegenerate.
Then, if $q_1(\widetilde x)=q_2(\widetilde x)=0$  and
\begin{equation}
T_{B_1}^*({\partial_{\bar z}\bar  \Phi}{\mathcal C}_2V_0)={\partial_{\bar z}\bar  \Phi}T_{B_1}^*
({\mathcal C}_2V_0)\quad \mbox{and}\quad
P_{A_2}^*({\partial_z\Phi}{\mathcal C}_1^*U_0)={\partial_z\Phi}P_{A_2}^*({\mathcal C}_1^*U_0),
\end{equation}
then we have
\begin{equation}\label{novka1}
(U_k e^{\tau\Phi}, {\bf H}(x,\partial_z,\partial_{\overline z})
(V_\ell e^{\tau\overline\Phi}))_{L^2(\Omega)}=o\left(\frac 1\tau\right)
\quad \mbox{as}\,\,\tau\rightarrow +\infty.
\end{equation}
\end{proposition}

{\bf Proof.} By our assumption
\begin{equation}\label{novka}
q_1(\widetilde x)=q_2(\widetilde x)=0.
\end{equation}
We consider several cases separately.
Let $\ell=0$ and $k=3$. We compute the asymptotics
\begin{eqnarray}\label{ruba1}
\mathcal K_1(\tau)=\vert (U_3 e^{\tau\overline\Phi}, {\bf H}(x,\partial_z,
\partial_{\overline z}) (V_0 e^{-\tau\overline\Phi}))_{L^2(\Omega)}\vert                                                 \nonumber\\
=\vert (U_3 e^{\tau\overline\Phi}, ({\mathcal C}_0\partial^2_{z\overline z}
+{\mathcal C}_2\partial^2_{\overline z\overline z}+{\bf B}_2\partial_{\overline z})
(V_0 e^{-\tau\overline\Phi}))_{L^2(\Omega)}\vert+o\left(\frac 1\tau\right)                                         \nonumber\\
=\vert (\partial_{\overline z}({\mathcal C}_0^*U_3 e^{\tau\overline\Phi}),
\partial_{z}(V_0 e^{-\tau\overline\Phi}))_{L^2(\Omega)}
+ (\partial_{\overline z}
({\mathcal C}_2^*U_3 e^{\tau\overline\Phi}),\partial_{\overline z}
(V_0 e^{\tau\overline\Phi}))_{L^2(\Omega)} +
(\partial_{\overline z}({\bf B}_2^*U_3 e^{\tau\overline\Phi}),
V_0 e^{-\tau\overline\Phi})_{L^2(\Omega)}\vert\nonumber\\
+ o(\frac {1}{\tau}) =\vert ({\mathcal C}_0^*\partial_{\overline z}
(U_3 e^{\tau\overline\Phi}), \partial_{z}(V_0 e^{-\tau\overline\Phi}))
_{L^2(\Omega)}+
({\mathcal C}_2^*\partial_{\overline z}(U_3 e^{\tau\overline\Phi}),
\partial_{\overline z}(V_0 e^{-\tau\overline\Phi}))_{L^2(\Omega)}\nonumber\\
+({\bf B}_2^*\partial_{\overline z}(U_3 e^{\tau\overline\Phi}),
V_0 e^{-\tau\overline\Phi})_{L^2(\Omega)}\vert
+ o\left(\frac{ 1}{\tau}\right)                                 \nonumber\\
 = \frac{1}{2}\vert ({\mathcal C}_0^*T_{B_1}(Q_1(1)U_2) e^{\tau\overline\Phi},
\partial_{z}(V_0 e^{-\tau\overline\Phi}))_{L^2(\Omega)}+ ({\mathcal C}_2^*
T_{B_1}(Q_1(1)U_2) e^{\tau{\overline\Phi}},\partial_{\overline z}
(V_0 e^{-\tau\overline\Phi}))_{L^2(\Omega)}
                                                        \nonumber\\
+ ({\bf B}_2^*T_{B_1}(Q_1(1)U_2) e^{\tau\overline\Phi},
V_0 e^{-\tau\overline\Phi})_{L^2(\Omega)}\vert
+  o\left(\frac{1}{\tau}\right)                \nonumber\\
=\frac{1}{2}\vert (e^{\tau\overline \Phi}{\mathcal C}_2^*T_{B_1}(Q_1(1)U_2 ),
\partial_{\overline z}(V_0 e^{-\tau\overline\Phi}))_{L^2(\Omega)} \vert
+ o\left(\frac{1}{\tau}\right)\nonumber\\
= \frac{\tau}{2}\vert  (U_2 ,Q_1(1)^*T_{B_1}^*
({\mathcal C}_2{\partial_{\bar z}\bar  \Phi}V_0) )
_{L^2(\Omega)} \vert+ o\left(\frac{ 1}{\tau}\right).
\end{eqnarray}
By the assumption, we have
$T_{B_1}^*({\partial_{\bar z}\bar  \Phi}{\mathcal C}_2V_0)={\partial_{\bar z}
\bar  \Phi}T_{B_1}^*
({\mathcal C}_2V_0)$.
Therefore
\begin{eqnarray}
\mathcal K_1(\tau)=\frac \tau2(U_2 ,Q_1(1)^*{\partial_{\bar z}
\bar  \Phi}T_{B_1}^*
({\mathcal C}_2V_0) )_{L^2(\Omega)}
+ o\left(\frac{ 1}{\tau}\right)\\
= (e^{\tau(\overline\Phi-\Phi)}T_{B_1}(Q_2(1)e^{\tau(\Phi-\overline\Phi)}U_1),
P_{A_1}^*(e^{\tau(\Phi-\overline\Phi)}Q_1(1)^*
{\partial_{\bar z}\overline\Phi}T_{B_1}^*({\mathcal C}_2V_0)))_{L^2(\Omega)}
+ o\left(\frac{ 1}{\tau}\right)\quad
\mbox{as}\quad\tau\rightarrow +\infty.\nonumber
\end{eqnarray}
By Proposition 8 of \cite{IY6}, we see
\begin{equation}\label{ii!}
e^{\tau(\overline\Phi-\Phi)}P_{A_1}^*(e^{\tau(\Phi-\overline\Phi)}Q_1(1)^*
{\partial_{\bar z}\bar  \Phi}T_{B_1}^*({\mathcal C}_2V_0))=\frac {Q_1(1)^*}{2\tau}T_{B_1}^*
({\mathcal C}_2V_0) + o_{L^2(\Omega)}\left(\frac{1}{\tau}\right)\quad\mbox{as}
\quad\tau\rightarrow +\infty
\end{equation}
and
\begin{equation}\label{metla}
U_1=\frac{q_1}{2\tau{\partial_z\Phi}}+o_{L^2(\Omega)}\left(\frac{1}{\tau}\right)
\quad\mbox{as}\quad\tau\rightarrow +\infty.
\end{equation}
Using Proposition 2.4 of \cite{IUY}, by (\ref{ii!}) and (\ref{metla}),
we obtain from  (\ref{ruba1})
\begin{equation}
\mathcal K_1(\tau)=\frac {1}{4\tau}\int_\Omega e^{\tau(\overline\Phi-\Phi)}
\frac{(q_1,Q_1(1)^*T_{B_1}^*({\mathcal C}_2V_0))}{2{\partial_z\Phi}}dx+o\left(\frac 1\tau\right)
=o\left(\frac 1\tau\right)\quad\mbox{as}\,\tau\rightarrow +\infty.
\end{equation}

The proof for the cases of $\ell=3$ and $k=0$  is similar.

Next we consider the case $k=2$ and $\ell=1.$
\begin{eqnarray}\label{ruba1X}
\mathcal K_2(\tau)= (U_2 e^{\tau\overline\Phi},
{\bf H}(x,\partial_z,\partial_{\overline z}) (V_1 e^{-\tau\overline\Phi}))
_{L^2(\Omega)}                            \\
= (U_2 e^{\tau\overline\Phi}, ({\mathcal C}_0\partial^2_{z\overline z}
+ {\mathcal C}_2\partial^2_{\overline z\overline z}+{\mathcal C}_1\partial^2_{ zz}
+ {\bf B}_1\partial_{ z}+{\bf B}_2\partial_{\overline z}) (V_1 e^{-\tau\overline\Phi}))_{L^2(\Omega)}+o\left(\frac 1\tau\right) \nonumber \\
= \mathcal K_{2,1}(\tau) + \mathcal K_{2,2}(\tau)+o\left(\frac 1\tau\right),
\nonumber
\end{eqnarray}
where
$$
\mathcal K_{2,1}(\tau)= (U_2 e^{\tau\overline\Phi}, ({\mathcal C}_0
\partial^2_{z\overline z} +{\mathcal C}_2\partial^2_{\overline z\overline z}+{\bf B}_2
\partial_{\overline z}) (V_1 e^{-\tau\overline\Phi}))_{L^2(\Omega)}
$$
and
$$
\mathcal K_{2,2}(\tau)= (U_2 e^{\tau\overline\Phi},
({\mathcal C}_1\partial^2_{ zz}+{\bf B}_1\partial_{ z})
(V_1 e^{-\tau\overline\Phi}))_{L^2(\Omega)}.
$$
We start with the computation of the asymptotics of $\mathcal K_{2,2}(\tau):$
\begin{eqnarray}\mathcal K_{2,2}(\tau)=(U_2 e^{\tau\overline\Phi},
({\mathcal C}_1\partial_{ z}+{\bf B}_1) (e^{-\tau\Phi} [\partial_z,P_{A_2}]
(q_2 e^{\tau(\Phi-\overline\Phi)})))_{L^2(\Omega)}\\
\nonumber+\tau(U_2 e^{\tau\overline\Phi}, ({\mathcal C}_1\partial_{ z}
+{\bf B}_1) (e^{-\tau\Phi} P_{A_2}((\tau({\partial_\zeta\Phi}(\zeta)-{\partial_z\Phi}(z))q_2+\partial_z q_2)
e^{\tau(\Phi-\overline\Phi)})))_{L^2(\Omega)}.\nonumber
\end{eqnarray}
By Proposition \ref{elkazelenaja}  we obtain
\begin{eqnarray}\label{Yukpos4}
\mathcal K_{2,2}(\tau)
=(U_2 e^{\tau\overline\Phi}, {\mathcal C}_1\partial_{ z} (e^{-\tau\Phi}
[\partial_z,P_{A_2}] (q_2 e^{\tau(\Phi-\overline\Phi)}))_{L^2(\Omega)}                                                      \\
+(U_2 e^{\tau\overline\Phi}, {\mathcal C}_1\partial_{ z}
(e^{-\tau\Phi} P_{A_2}((\tau({\partial_\zeta\Phi}(\zeta)-{\partial_z\Phi}(z))q_2+\partial_z q_2)
e^{\tau(\Phi-\overline\Phi)}))_{L^2(\Omega)}
+ o\left(\frac 1\tau\right)\quad\mbox{as}\quad\tau\rightarrow +\infty.\nonumber
\end{eqnarray}

By (\ref{victory}) there exists a constants $C_j(\tau)$ such that
\begin{equation}\label{victory1}
P_{A_2}(({\partial_z\Phi}(z)-{\partial_\zeta\Phi}(\zeta))q_2 e^{\tau(\Phi-\overline\Phi)})
=\Pi C_1(\tau)+\widetilde \Pi C_2(\tau),
\end{equation}
where the matrices $\Pi,\widetilde \Pi$ are introduced in the previous section.
The stationary phase argument yields
\begin{equation}\label{victory2}
C_j(\tau)=O(\frac{1}{\tau^\frac 32})\quad\mbox{as}\quad\tau\rightarrow
+\infty\quad \forall j\in\{1,2\}.
\end{equation}
By (\ref{victory1}), (\ref{victory2}) and Proposition \ref{elkazelenaja}, we
have
\begin{equation}\label{Yukos1}
\tau(U_2 e^{\tau\overline\Phi}, {\mathcal C}_1\partial_{ z} (e^{-\tau\Phi}
P_{A_2}(({\partial_z\Phi}(z)-{\partial_\zeta\Phi}(\zeta))q_2 e^{\tau(\Phi-\overline\Phi)}))
_{L^2(\Omega)}=o(\frac{1}{\tau})\quad\mbox{as}\quad\tau\rightarrow +\infty.
\end{equation}

Next we compute the following asymptotics
\begin{eqnarray}\label{zmei}
(U_2e^{\tau\bar\Phi},\mathcal C_1\partial_z[e^{-\tau\Phi}P_{A_1}(\partial_z q_2e^{\tau(\Phi-\bar\Phi)})])_{L^2(\Omega)}=(U_2e^{\tau\bar\Phi},\tau\mathcal C_1e^{-\tau\Phi}P_{A_1}((\Phi(\zeta)-\Phi(z))\partial_z q_2e^{\tau(\Phi-\bar\Phi)}))_{L^2(\Omega)}\nonumber\\+(U_2e^{\tau\bar\Phi},\tau\mathcal C_1e^{-\tau\Phi}[\partial_z,P_{A_1}](\partial_z q_2e^{\tau(\Phi-\bar\Phi)}))_{L^2(\Omega)}\\+(U_2e^{\tau\bar\Phi},\tau\mathcal C_1e^{-\tau\Phi}[\partial_z,P_{A_1}](\partial^2_{zz} q_2e^{\tau(\Phi-\bar\Phi)}))_{L^2(\Omega)}.\nonumber
\end{eqnarray}
By Proposition \ref{elkazelenaja} and Proposition \ref{pen} the last two terms in the above  formula are of the order $o(\frac {1}{\tau}).$

By (\ref{victory}) there exists a constants $\widetilde C_j(\tau)$ such that
\begin{equation}\label{victory11}
P_{A_2}(({\partial_z\Phi}(z)-{\partial_\zeta\Phi}(\zeta))q_2 e^{\tau(\Phi-\overline\Phi)})
=\Pi\widetilde  C_1(\tau)+\widetilde \Pi\widetilde  C_2(\tau),
\end{equation}
where the matrices $\Pi,\widetilde \Pi$ are introduced in the previous section.
The stationary phase argument yields
\begin{equation}\label{victory21}
\widetilde C_j(\tau)=O(\frac{1}{\tau})\quad\mbox{as}\quad\tau\rightarrow
+\infty\quad \forall j\in\{1,2\}.
\end{equation}
By (\ref{victory11}), (\ref{victory21}) and Proposition \ref{elkazelenaja}, we
have
$$
(U_2e^{\tau\bar\Phi},\tau\mathcal C_1e^{-\tau\Phi}P_{A_1}((\Phi(\zeta)-\Phi(z))\partial_z q_2e^{\tau(\Phi-\bar\Phi)}))_{L^2(\Omega)}=o(\frac{1}{\tau})\quad\mbox{as}\,\,\tau\rightarrow +\infty.
$$
Hence the each of the three terms of the right hand side of the formula (\ref{zmei}) equal $o(\frac 1\tau).$

Next we compute the asymptotics:
\begin{eqnarray*}
&&(U_2 e^{\tau\overline\Phi}, {\mathcal C}_1\partial_{ z}
(e^{-\tau\Phi} [\partial_z,P_{A_2}] (q_2 e^{\tau(\Phi-\overline\Phi)}))
_{L^2(\Omega)}\\
&& =(U_2 e^{\tau\overline\Phi}, {\mathcal C}_1 (-\tau{\partial_z\Phi}e^{-\tau\Phi}
[\partial_z,P_{A_2}] (q_2 e^{\tau(\Phi-\overline\Phi)})+ [\partial_z,P_{A_2}]
(q_2\tau {\partial_z\Phi} e^{\tau(\Phi-\overline\Phi)})\nonumber\\
&&+[\partial_z,P_{A_2}] (\partial_zq_2 e^{\tau(\Phi-\overline\Phi)}))_{L^2(\Omega)} \\
&& + (U_2 e^{\tau\overline\Phi}, {\mathcal C}_1 e^{-\tau\Phi}[\partial_z,
[\partial_z,P_{A_2}]] (q_2 e^{\tau(\Phi-\overline\Phi)}))_{L^2(\Omega)}.
\end{eqnarray*}

By Propositions \ref{elkazelenaja} and \ref{zanuda} we have
\begin{equation}
(U_2 e^{\tau\overline\Phi},[\partial_z,P_{A_2}]
(\partial_zq_2 e^{\tau(\Phi-\overline\Phi)}))_{L^2(\Omega)}
= o\left(\frac 1\tau\right)\quad\mbox{as}\quad\tau\rightarrow +\infty.
\end{equation}

By Propositions \ref{elkazelenaja} and \ref{zanuda1} we have
\begin{equation}\label{Yukos2}
(U_2 e^{\tau\overline\Phi}, {\mathcal C}_1 (-\tau{\partial_z\Phi}e^{-\tau\Phi}
[\partial_z,P_{A_2}] (q_2 e^{\tau(\Phi-\overline\Phi)})\nonumber+ [\partial_z,P_{A_2}] (q_2\tau {\partial_z\Phi} e^{\tau(\Phi-\overline\Phi)}))_{L^2(\Omega)}
=o\left(\frac 1\tau\right)
\end{equation}
and
\begin{equation}\label{Yukos}
(U_2 e^{\tau\overline\Phi}, {\mathcal C}_1 e^{-\tau\Phi}[\partial_z,
[\partial_z,P_{A_2}]] (q_2 e^{\tau(\Phi-\overline\Phi)}))_{L^2(\Omega)}
=o\left(\frac 1\tau\right)
\end{equation}
as $\tau\rightarrow +\infty$.

Combining (\ref{Yukpos4}) and (\ref{Yukos1})-(\ref{Yukos}), we prove that
\begin{equation}
\mathcal K_{2,2}(\tau)=o\left(\frac 1\tau\right)\quad\mbox{as}
\quad\tau\rightarrow +\infty.
\end{equation}
Next we compute the asymptotics for $\mathcal K_{2,1}(\tau)$.
By (\ref{mursilka1}) we have
\begin{equation}
\mathcal K_{2,1}(\tau)= (U_2 e^{\tau\overline\Phi}, ({\mathcal C}_0
\partial_{z} +{\mathcal C}_2\partial_{\overline z}+{\bf B}_2)
(-A_2V_1 +q_2) e^{-\tau\overline\Phi})_{L^2(\Omega)}.
\end{equation}
Proposition \ref{elkazelenaja} yields
\begin{equation}
(U_2 e^{\tau\overline\Phi},{\bf B}_2 (-A_2V_1 +q_2) e^{-\tau\overline\Phi})
_{L^2(\Omega)}=o\left(\frac 1\tau\right)\quad\mbox{as}\quad\tau\rightarrow
+\infty.
\end{equation}
By (\ref{mursilka1}) and  Proposition \ref{elkazelenaja} we have
\begin{eqnarray}\label{makaka1}
(U_2 e^{\tau\Phi}, {\mathcal C}_0\partial_{ z}((-A_2V_1 +q_2)
e^{-\tau\overline\Phi}))_{L^2(\Omega)}=
(U_2 e^{\tau\Phi}, {\mathcal C}_0(-(\partial_{ z}A_2)V_1
-A_2\partial_zV_1
+ \partial_zq_2) e^{-\tau\overline\Phi}))_{L^2(\Omega)}            \nonumber\\
= -(U_2 , {\mathcal C}_0A_1\partial_{ z}V_1  )_{L^2(\Omega)}
+ o\left(\frac 1\tau\right)
= -(U_2 e^{\tau(\overline\Phi-\Phi)}, {\mathcal C}_0A_1[\partial_{ z},
P_{A_2}] (q_2e^{\tau(\Phi-\overline\Phi)}))_{L^2(\Omega)}     \nonumber\\
- (U_2 e^{\tau(\overline\Phi-\Phi)}, {\mathcal C}_0A_1P_{A_2} (\partial_zq_2
e^{\tau(\Phi-\overline\Phi)}))_{L^2(\Omega)}\nonumber\\
+(U_2 e^{\tau(\overline\Phi-\Phi)}, \tau{\mathcal C}_0A_1P_{A_2}
(({\partial_z\Phi}(z)-{\partial_\zeta\Phi}(\zeta))q_2e^{\tau(\Phi-\overline\Phi)}))_{L^2(\Omega)}
+ o\left(\frac 1\tau\right).
\end{eqnarray}
By the stationary phase argument, there exists a constant $C$ independent of
$\tau$ such that
$$
\tau\Vert P_{A_2}( ({\partial_z\Phi}(z)-{\partial_\zeta\Phi}(\zeta))q_2e^{\tau(\Phi-\overline\Phi)})
\Vert_{L^2(\Omega)}\le C.
$$
Hence by Proposition \ref{elkazelenaja}, we have
\begin{equation}
(U_2 e^{\tau(\overline\Phi-\Phi)}, \tau{\mathcal C}_0A_1P_{A_2}
(({\partial_z\Phi}(z)-{\partial_\zeta\Phi}(\zeta))q_2e^{\tau(\Phi-\overline\Phi)}))_{L^2(\Omega)}
=o\left(\frac 1\tau\right)\quad\mbox{as}\quad\tau\rightarrow +\infty.
\end{equation}
Applying Proposition \ref{elkazelenaja} we have
\begin{equation}
(U_2 e^{\tau(\overline\Phi-\Phi)}, {\mathcal C}_0A_1P_{A_2} (\partial_zq_2
e^{\tau(\Phi-\overline\Phi)}))_{L^2(\Omega)}=o\left(\frac 1\tau\right)
\quad\mbox{as}\quad\tau\rightarrow +\infty.
\end{equation}
By Propositions \ref{elkazelenaja} and \ref{zanuda}, we have
\begin{eqnarray}\label{makaka}
\vert(U_2 e^{\tau(\overline\Phi-\Phi)}, {\mathcal C}_0A_1[\partial_{ z},
P_{A_2}] (q_2e^{\tau(\Phi-\overline\Phi)}))_{L^2(\Omega)}\vert\\
\le C \Vert U_2\Vert_{L^2(\Omega)}\Vert [\partial_{ z},P_{A_2}]\Vert
_{\mathcal L(C^0(\overline\Omega), L^2(\Omega))}\Vert q_2\Vert
_{C^0(\overline\Omega)}=
o\left(\frac 1\tau\right)\quad\mbox{as}\quad\tau\rightarrow +\infty.\nonumber
\end{eqnarray}
Then  from (\ref{makaka1})-(\ref{makaka}) we have
\begin{eqnarray}\label{idiot}
\mathcal K_{2,1}(\tau)=(U_2 e^{\tau\overline\Phi}, {\mathcal C}_2
\partial_{\overline z}((-A_2V_1 +q_2) e^{-\tau\overline\Phi}))_{L^2(\Omega)}
= (U_2, {\mathcal C}_2\tau {\partial_{\bar z}\bar  \Phi}q_2 )_{L^2(\Omega)}       \nonumber\\
+(U_2 e^{\tau\overline\Phi}, {\mathcal C}_2((-\partial_{\overline z}A_2V_1
+ \frac 12 A_2(-A_2V_1+q_2) +\partial_{\overline z}q_2)
e^{\tau\overline\Phi}))_{L^2(\Omega)}\\
= \tau(e^{\tau(\overline\Phi-\Phi)}T_{B_1}(Q_2(1)e^{\tau(\Phi-\overline\Phi)}
U_1),P_{A_1}^*(e^{\tau(\Phi-\overline\Phi)}
{\partial_{\bar z}\bar  \Phi}{\mathcal C}_2q_2)_{L^2(\Omega)}+o\left(\frac 1\tau\right).
                                        \nonumber
\end{eqnarray}
By Proposition 8 of \cite{IY6} we obtain
\begin{equation}\label{metla1}
e^{\tau(\overline\Phi-\Phi)}P_{A_1}^*(e^{\tau(\Phi-\overline\Phi)}
{\partial_{\bar z}\bar  \Phi}{\mathcal C}_2q_2)=\frac {1}{2\tau}{\mathcal C}_2q_2
+ o_{L^2(\Omega)}\left(\frac{1}{\tau}\right)
\quad\mbox{as}\quad\tau\rightarrow +\infty.
\end{equation}
By (\ref{metla}) and (\ref{metla1}) we have
$$
\mathcal K_{2,1}(\tau)=\frac{1}{4\tau}\int_\Omega
e^{\tau(\overline\Phi-\Phi)} \frac{({\mathcal C}_2q_2,q_1)}{{\partial_z\Phi}}dx
+ o\left(\frac 1\tau\right)\quad\mbox{as}\quad\tau\rightarrow +\infty.
$$
Therefore from  Proposition 2.4 of \cite{IUY} we have
\begin{equation}
\mathcal K_{2,1}(\tau)=o\left(\frac 1\tau\right)\quad
\mbox{as}\quad\tau\rightarrow +\infty.
\end{equation}
The proof in the case $\ell=2$ and $k=1$  is similar.

Finally we consider the case when $k+\ell\ge 4, k,\ell\in {\Bbb N}_+.$
We set $\Psi=\Phi$ and $\Psi=\overline{\Phi}$ for $k=0,1$ and
$\Psi=\overline \Phi$ and $\Psi_1=\Phi$ otherwise.
Then
\begin{eqnarray}\label{ux}
(U_k e^{\tau\Psi}, {\bf H}(x,\partial_z,\partial_{\overline z})
(V_\ell e^{-\tau\Psi_1}))_{L^2(\Omega)}           \\
=(U_k e^{\tau\Psi}, ({\mathcal C}_0\partial^2_{z\overline z}+{\mathcal C}_1
\partial^2_{z z} +{\mathcal C}_2\partial^2_{\overline z\overline z}+{\bf B}_1\partial_{ z}
+{\bf B}_2\partial_{\overline z}) (V_\ell e^{-\tau\Psi_1}))_{L^2(\Omega)}
+o\left(\frac 1\tau\right)                     \nonumber\\
= (\partial_{\overline z}({\mathcal C}_0^*U_k e^{\tau\Psi}), \partial_{z}
(V_\ell e^{-\tau\Psi_1}))_{L^2(\Omega)}
+ (\partial_{ z}({\mathcal C}_1^*U_k e^{\tau\Psi}), \partial_{z}
(V_\ell e^{-\tau\Psi_1}))_{L^2(\Omega)}\nonumber\\ +(\partial_{\overline z}({\mathcal C}_2^*U_k e^{\tau\overline\Phi}),\partial_{\overline z}
(V_\ell e^{-\tau\Psi_1}))_{L^2(\Omega)} \nonumber\\
+ (U_k e^{\tau\Psi},{\bf B}_1\partial_{ z}+{\bf B}_2 \partial_{\overline z}
(V_\ell e^{-\tau\Psi_1}))_{L^2(\Omega)}+o\left(\frac 1\tau\right)
\quad\mbox{as}\quad\tau\rightarrow +\infty.\nonumber
\end{eqnarray}

To complete the proof of the proposition we show that each of the four integrals in the right-hand side of (\ref{ux}) are equal to $o(\frac 1\tau).$
In order to prove this fact it is suffices to show that  for any regular matrix $\mathcal C(x)$
\begin{equation}
(\partial_1 e^{\tau\Psi},\mathcal C \partial_2
(V_\ell e^{-\tau\Psi_1}))_{L^2(\Omega)}=o(\frac 1\tau)\quad\mbox{as}\quad\tau\rightarrow +\infty \quad \partial_j\in\{\partial_z,\partial_{\bar z}\}.
\end{equation}
In order to prove the above formula consider several cases.
If $k\ge 1$, then we have two cases:

{\bf Case A.} Either  $k=1$ and $\partial=\partial_z$ or $k\ge 2$ and
$\partial=\partial_{\overline z}.$
Then by (\ref{mursilka}) there exist matrices $A_k$ and functions
$F_k( \tau,\cdot)$ such that
$$
\partial (U_k e^{\tau\Psi})=A_kU_k e^{\tau\Psi}+ F_k e^{\tau\Psi},
$$
and for some constant $C$ independent of $\tau$, we have an estimate:
$$
\Vert F_k\Vert_{L^2(\Omega)}\le C\Vert U_{k-1}\Vert_{L^2(\Omega)}.
$$
By Proposition \ref{elkazelenaja} for any $\epsilon\in (0,1)$ there exists
\begin{equation}\label{Y1}
\Vert e^{-\tau\varphi}\partial (U_k e^{\tau\Psi})\Vert_{L^2(\Omega)}
\le C(\epsilon)/\tau^{\kappa},\quad \kappa=\max \{0, k-1-\epsilon\}.
\end{equation}

{\bf Case B.} Let $k=1$ and $\partial=\partial_{\overline z}$ or $k\ge 2$
and $\partial=\partial_{ z}.$

Then by Proposition \ref{zanuda1}, for any $\epsilon\in (0,1)$ there exists
a constant $C(\epsilon)$ such that
\begin{equation}\label{Y2}
\Vert e^{-\tau\varphi}\partial (U_k e^{\tau\Psi})\Vert_{L^2(\Omega)}
\le C(\epsilon)\tau^\epsilon \Vert U_{k-1}\Vert_{L^2(\Omega)} \le C(\epsilon)/\tau^{\kappa}, \quad\kappa=\max \{0, k-1-2\epsilon\}.
\end{equation}

If $\ell\ge 1$ we also have two cases.

{\bf Case C.} Let $\ell=1$ and $\partial=\partial_{\overline z}$ or
$\ell\ge 2$ and $\partial=\partial_{z}.$

Then by (\ref{mursilka1}) there exist matrices $\widetilde A_\ell$ and
functions $\widetilde F_\ell$ such that
$$
\partial (V_\ell e^{-\tau\Psi_1})=\widetilde A_\ell V_\ell e^{-\tau\Psi_1}
+ \widetilde F_\ell e^{-\tau\Psi_1},
$$
where
$$
\Vert F_\ell\Vert_{L^2(\Omega)}\le C\Vert V_{\ell-1}\Vert_{L^2(\Omega)}.
$$
By Proposition  \ref{elkazelenaja} we obtain
\begin{equation}\label{XY}
\Vert e^{\tau\varphi}\partial (V_\ell e^{-\tau\Psi})\Vert_{L^2(\Omega)}
\le C(\epsilon)/\tau^{\widetilde \kappa}, \quad \widetilde\kappa
=\max \{0, \ell-1-\epsilon\}.
\end{equation}

{\bf Case D.} Let $\ell=1$ and $\partial=\partial_{ z}$ or $\ell\ge 2$
and $\partial=\partial_{\overline  z}.$

Then by Proposition \ref{zanuda1} for any $\epsilon\in (0,1)$,
there exists a constant $C(\epsilon)$ such that
\begin{equation}\label{XY1}
\Vert e^{\tau\varphi}\partial (V_\ell e^{-\tau\Psi_1})\Vert_{L^2(\Omega)}
\le C(\epsilon)\tau^\epsilon \Vert V_{\ell-1}\Vert_{L^2(\Omega)}
\le C(\epsilon)/\tau^{\widetilde \kappa},\quad\widetilde \kappa
=\max \{0, \ell-1-2\epsilon\}.
\end{equation}

Hence if  $k\ge 1,$ $\ell\ge 1$ and $k+\ell\ge 4$, then by (\ref{Y1}), (\ref{Y2}), (\ref{XY}), (\ref{XY1})
\begin{equation}\label{uxx}
\vert(\partial_1U_k e^{\tau\Psi}, \mathcal C
\partial_2(V_\ell e^{-\tau\Psi_1}))_{L^2(\Omega)}\vert
\le \frac{C}{\tau^{\kappa+\widetilde \kappa}}=\frac{C(\epsilon)}
{\tau^{k+\ell-2-4\epsilon}}=o(\frac{1}{\tau})
\quad\mbox{as}\,\,\tau\rightarrow +\infty.
\end{equation}

Now the only case we need to consider is the situation when
either  $\ell=0, k\ge 4$ or $\ell\ge 4, k=0.$ Let for example $\ell=0.$ Then
using the Proposition  \ref{elkazelenaja} we obtain
\begin{equation}\label{uxx}
\vert(\partial_1U_k e^{\tau\Psi}, \mathcal C
\partial_2(V_\ell e^{-\tau\Psi_1}))_{L^2(\Omega)}\vert
\le \frac{C}{\tau^{\kappa+\widetilde \kappa}}
\le \frac{C(\epsilon)}{\tau^{\kappa-2-2\epsilon}}
=o\left(\frac{1}{\tau}\right)\quad\mbox{as}\,\,\tau\rightarrow +\infty.
\end{equation}
The proof of proposition is complete.

 $\blacksquare$

In a way similar to Proposition \ref{zanoza}, we prove

\begin{proposition}\label{zanoza!}Let $k, \ell\in\Bbb N_+$ and
$k+\ell\ge 3$ and ${\bf H}(x,\partial_z,\partial_{\bar z})$ be a second-order differential operator
with compactly supported smooth coefficients.
Let the functions $\widetilde U_k$ and $\widetilde V_\ell$ be given by
Proposition \ref{mursilka00} and $\widetilde q_1, \widetilde q_2$ be
given by (\ref{mishka66}). Then
$$
(\widetilde U_k e^{\tau\widetilde \Phi}, {\bf H}(x,\partial_z,\partial_{\overline z})
(\widetilde V_\ell e^{\tau\overline{\widetilde\Phi}}))_{L^2(\Omega)}
= o\left(\frac 1\tau\right)\quad \mbox{as}\,\,\tau\rightarrow +\infty.
$$
\end{proposition}

Henceforth we assume that the domain $\Omega$ is a ball
centered at origin.
Denote
$$
\widetilde {\frak C}_1(x)=P_{A_2}^*( 2\partial_z\widetilde \Phi\partial_{ z}
(\mathcal C_1^* \widetilde U_0)+2\mathcal C_1^*\widetilde U_0
+ \partial_z{\widetilde \Phi}\partial_{\overline z}(\mathcal C_0^*{\widetilde U}_0)
+ \partial_z\widetilde \Phi{\bf B}_1^*\widetilde U_0+ Q_1(2)^*P_{A_2}^*(\mathcal C_1^*
\partial_z\widetilde \Phi\widetilde U_0)+\mathcal C_1^*\partial_z\widetilde \Phi\widetilde q_1),
$$
$$
\widetilde {\frak C}_2(x)=-T_{B_1}^*(\mathcal C_2 (2\partial_{\bar z}\overline{\widetilde \Phi}
\partial_{\overline z}\widetilde V_0+2\widetilde V_0)+\mathcal C_0
\partial_{\bar z}\overline{\widetilde \Phi}\partial_{z} V_0
+\partial_{\bar z} \overline {\widetilde \Phi}{\bf B}_2\widetilde V_0+Q_2(1)^*T_{B_1}^*
(\mathcal C_2\partial_{\bar z}\overline{\widetilde  \Phi}\widetilde V_0)
+ \mathcal C_2\partial_{\bar z}\overline{\widetilde \Phi}\widetilde q_2).
$$
We have
\begin{proposition} \label{garmoshka}
Let ${\bf H}(x,\partial_z,\partial_{\overline z})$ be a second-order
differential operator given by (\ref{poker}) with compactly supported smooth
coefficients in $\Omega$.
Let $\widetilde U, \widetilde V$ be the solutions constructed in
Proposition \ref{mursilka00}.   Suppose that
\begin{equation}\label{lida1}
(\widetilde U,{\bf H}(x,\partial_z,\partial_{\overline z})\widetilde V)
_{L^2(\Omega)}=o\left(\frac 1\tau\right)\quad \mbox{as}\quad \tau\rightarrow
+\infty.
\end{equation}
Then the following equalities hold true:
\begin{equation}\label{4.177}
P_{A_2}^*[z^k\mathcal C_1^*U_0]\vert_{\partial\Omega}
= T_{B_1}^*[z^k\mathcal C_2V_0]\vert_{\partial\Omega}=0
\quad k\in \{0,1,2\},
\end{equation}
\begin{equation}\label{4.178}
P_{A_2}^*( 2z^k\partial_{ z}(\mathcal C_1^*  U_0)+2z^k\mathcal C_1^*U_0
+z^k\partial_{\overline z}(\mathcal C_0^*{ U}_0) + z^k{\bf B}_1^* U_0
+ Q_1(2)P_{A_2}^*(\mathcal C_1^*z^k U_0)+\mathcal C_1^* z^kq_1)
=0\quad\mbox{on}\quad\partial\Omega
\end{equation}
and
\begin{equation}\label{4.179}
-T_{B_1}^*(\mathcal C_2 (2\overline z^k\partial_{\overline z}\widetilde V_0
+2\overline z^k V_0)+\mathcal C_0\overline z^k\partial_{z} V_0
+\overline  z^k {\bf B}_2 V_0+Q_2(1)T_{B_1}^*(\mathcal C_2\overline z^k V_0)
+ \mathcal C_2\overline z^kq_2)=0\quad\mbox{on}\quad\partial\Omega.
\end{equation}
\end{proposition}

 {\bf Proof.} From Proposition \ref{zanoza!} we obtain
\begin{eqnarray}\label{elka}
(\widetilde U,{\bf H}(x,\partial_z,\partial_{\overline z})\widetilde V)
_{L^2(\Omega)}
= ((\widetilde U_0-\widetilde U_1+\widetilde U_2)e^{\tau\widetilde \Phi},
{\bf H}(x,\partial_z,\partial_{\overline z})((\widetilde V_0-\widetilde V_1
+\widetilde V_2)e^{-\tau\widetilde \Phi}))_{L^2(\Omega)}\nonumber\\
+ o\left(\frac 1\tau\right)\quad \mbox{as}\quad \tau\rightarrow +\infty.
\end{eqnarray}
Then applying  Proposition \ref{zanoza2}, we have
\begin{eqnarray}\label{999A}
(\widetilde U,{\bf H}(x,\partial_z,\partial_{\overline z})\widetilde V)
_{L^2(\Omega)}=-(\widetilde U_0e^{\tau\widetilde\Phi}, {\bf H}(x,\partial_z,
\partial_{\overline z})(\widetilde V_1e^{-\tau\overline {\widetilde \Phi}}))
_{L^2(\Omega)}\\
- (\widetilde U_1e^{\tau\widetilde\Phi}, {\bf H}(x,\partial_z,\partial
_{\overline z})(\widetilde V_0e^{-\tau\overline {\widetilde \Phi}}))
_{L^2(\Omega)}
+ (\widetilde U_0 e^{\tau\widetilde \Phi},{\bf  H}(x,\partial_z,\partial
_{\overline z}) (\widetilde V_2 e^{-\tau\overline{\widetilde \Phi}}))
_{L^2(\Omega)}\nonumber\\
+ (\widetilde U_2 e^{\tau\widetilde\Phi},{\bf  H}(x,\partial_z,\partial
_{\overline z}) (\widetilde V_0 e^{-\tau\overline{\widetilde\Phi}}))
_{L^2(\Omega)}+(\widetilde U_1 e^{\tau{\widetilde\Phi}},
{\bf  H}(x,\partial_z,\partial_{\overline z})
(\widetilde V_1 e^{-\tau\overline{\widetilde\Phi}}))_{L^2(\Omega)}\nonumber\\
+ o\left(\frac{1}{\root\of{\tau}}\right)\quad\mbox{as}\quad\tau\rightarrow
+\infty.                        \nonumber
\end{eqnarray}
The function $\widetilde\psi$ on $\partial\Omega$ has only two critical points,
that is, the point of minimum $x_{min}$ and the point of maximum $x_{max}.$
Moreover after appropriate choice of the function $\widetilde \Phi$
we can assume that  $x_{max}=\hat x$ where $\hat x$ is an arbitrary fixed
point from $\partial\Omega.$  The second tangential derivatives of the
function $\widetilde \psi$ at points $x_{max}$ and $x_{min}$ are not equal
to zero. Hence the first asymptotic term on the right-hand side of
(\ref{999A}) is of the order $\root\of{\tau}$. Therefore we can remove
some terms and write (\ref{999A}) as
\begin{eqnarray}
(\widetilde U,{\bf H}(x,\partial_z,\partial_{\overline z})\widetilde V)
_{L^2(\Omega)}                             \nonumber\\
= \frak J_\tau (P_{A_2}^*{\bf H}(x,\partial_z+\tau\partial_z\widetilde\Phi,
\partial_{\overline z})^*\widetilde U_0, \widetilde q_2)
+ \frak J_\tau(\widetilde q_1, T_{B_1}^*{\bf H}(x,\partial_z,\partial
_{\overline z}-\tau\partial_{\bar z}\overline{\widetilde \Phi})V_0)
+o(\frac{1}{\root\of{\tau}})            \nonumber\\
=\int_{\partial\Omega} (P_{A_2}^*{\bf H}(x,\partial_z+\tau\partial_z\widetilde\Phi,
\partial_{\overline z})^*\widetilde
U_0,\widetilde q_2)\frac{(\nu_1-i\nu_2)}{2\tau \partial_z\widetilde\Phi}
e^{\tau(\widetilde\Phi-\overline{\widetilde \Phi})}d\sigma\nonumber\\
+ \int_{\partial\Omega} (\widetilde q_1, T_{B_1}^*{\bf H}
(x,\partial_z,\partial_{\overline z}-\tau\partial_{\bar z}\overline{\widetilde \Phi})V_0)
\frac{(\nu_1-i\nu_2)}{2\tau \partial_z\widetilde\Phi}e^{\tau(\widetilde\Phi
-\overline {\widetilde\Phi})}d\sigma
+ o\left(\frac{1}{\root\of{\tau}}\right)\quad\mbox{as}
\quad\tau\rightarrow +\infty.\nonumber
\end{eqnarray}

Then applying the stationary phase argument, we obtain that there exists
 a function $\kappa$, not equal to zero at any point of $\partial\Omega$
such that
\begin{eqnarray}\label{999AA}
(\widetilde U,{\bf H}(x,\partial_z,\partial_{\overline z})\widetilde V)
_{L^2(\Omega)}                        \\
= \frac{(\nu_1-i\nu_2)\root\of{\tau} \kappa}{2 \partial_z\widetilde\Phi}
(x_{max})\left\{ (\widetilde q_1, T_{B_1}^*[z^k\mathcal C_2V_0])\right.
+\left.(P_{A_2}^*[z^k\mathcal C_1^*U_0],\widetilde  q_2)\right\}(x_{max})
e^{2\tau\widetilde\psi(x_{max})}                          \nonumber\\
+ \frac{(\nu_1-i\nu_2)\root\of{\tau} \kappa}{2\tau \partial_z\widetilde\Phi}
(x_{min})\left\{ (\widetilde q_1, T_{B_1}^*[z^k\mathcal C_2V_0])\right.
\left.+(P_{A_2}^*[z^k\mathcal C_1^*U_0],\widetilde  q_2)\right\}
(x_{min})e^{2\tau\widetilde\psi(x_{min})}+o({\root\of{\tau}}).\nonumber
\end{eqnarray}
Since $\psi(x_{min})\ne \psi(x_{max})$, the above equality implies
$$
\left \{(\widetilde q_1, T_{B_1}^*[z^k\mathcal C_2V_0])
+(P_{A_2}^*[z^k\mathcal C_1^*U_0],\widetilde q_2)\right\}(x_{max})=0.
$$
This equality and Proposition \ref{nikita}  imply (\ref{4.177}).
Since the matrices $\mathcal C_j$ have compact supports, the equalities
(\ref{4.177}) imply
\begin{equation}
\frac{\partial^j}{\partial\nu^j} T_{B_1}^*[z^k\mathcal C_2V_0]\vert
_{\partial \Omega}
= \frac{\partial^j}{\partial\nu^j} P_{A_2}^*[z^k\mathcal C_1^*U_0]\vert
_{\partial \Omega}=0\quad \forall j\in \{0,1,\dots, 4\}.
\end{equation}
This equality implies in particular that
\begin{equation}\label{lada}
\frak J_\tau (P_{A_2}^*[z^k\mathcal C_1^*U_0],\widetilde q_2)
+ \frak J_\tau(\widetilde q_1, T_{B_1}^*[z^k\mathcal C_2V_0])
= o\left(\frac{1}{\root\of{\tau}}\right)
\quad\mbox{as}\quad\tau\rightarrow +\infty.
\end{equation}

By (\ref{lada}), we have
\begin{equation}\label{999AAA}
(\widetilde U,{\bf H}(x,\partial_z,\partial_{\overline z})\widetilde V)
_{L^2(\Omega)}
\end{equation}
\begin{eqnarray*}
&& =\frak J_\tau (P_{A_2}^*{\bf H}(x,\partial_z+\tau\partial_z\widetilde\Phi,
\partial_{\overline z})^*\widetilde U_0,\widetilde q_2)
+ \frak J_\tau(\widetilde q_1, T_{B_1}^*{\bf H}(x,\partial_z,
\partial_{\overline z}-\tau\partial_{\bar z}\overline{\widetilde \Phi})V_0) \\
&&+ o\left(\frac{1}{\root\of{\tau}}\right)\\
&& =\int_{\partial\Omega} (\widetilde{\frak C_1}, \widetilde q_2)
\frac{e^{2i\theta}(\nu_1-i\nu_2)}{2 \partial_z\widetilde\Phi}
e^{\tau(\widetilde\Phi-\overline{\widetilde \Phi})}d\sigma\\
&& +\int_{\partial\Omega} (\widetilde q_1,\widetilde {\frak C}_2)
\frac{e^{2i\theta}(\nu_1-i\nu_2)}{2 \partial_z\widetilde\Phi}
e^{\tau(\widetilde\Phi-\overline {\widetilde\Phi})}d\sigma
+ o\left(\frac{1}{\root\of{\tau}}\right)\\
&& \qquad \mbox{as}\quad\tau\rightarrow +\infty.
\end{eqnarray*}
Applying the stationary phase argument to the right-hand side of
(\ref{999AAA}), we obtain
$$
(\widetilde{\frak C_1},\widetilde q_2)(x_{max})+(\widetilde q_1,
\widetilde {\frak C}_2)(x_{max})=0.
$$
Then using Proposition \ref{nikita}, we obtain (\ref{4.178}) and (\ref{4.179}).
$\blacksquare$

Denote
\begin{eqnarray*}
&& \frak H_{U_0,V_0}(x)=\{(U_0,{\mathcal C}_1\partial^2_{zz}V_0)
-(\partial_{\overline z}(\mathcal C_0^* U_0),\partial_zV_0)\\
&& + (\partial^2_{\overline z\overline z}(\mathcal C_2^*U_0),V_0)
+ (U_0,{\bf B}_1\partial_zV_0)
-(\partial_{\overline z}({\bf B}_2^*U_0),V_0)+(U_0,{\bf B}_0V_0) \}(x).
\end{eqnarray*}

The following proposition plays the crucial part in the proof of the
uniqueness of the determination of coefficients for the Navier-Stokes
equations and the Lam\'e system.

\begin{proposition} \label{gavnuk}
Let ${\bf H}(x,\partial_z,\partial_{\overline z})$ be a second-order
differential operator given by (\ref{poker}) with smooth coefficients which
have compact supports in $\Omega$. Assume that the restriction of the function $\psi$ on $\partial\Omega$ has
a finite number of critical points and all these points are nondegenerate.
Let $U, V$ be the solutions constructed in Proposition \ref{mursilka0},
$q_1(\widetilde x)=q_2(\widetilde x)=0$ and (\ref{4.177})-(\ref{4.179}) hold true.  Suppose that
\begin{equation}\label{lida1}
(U,{\bf H}(x,\partial_z,\partial_{\overline z})V)_{L^2(\Omega)}
=o\left(\frac 1\tau\right)\quad \mbox{as}\quad \tau\rightarrow +\infty.
\end{equation}
Then the following equality holds true:
\begin{eqnarray}\label{kaput3}
(P_{A_2}^*[\mathcal C_1^*U_0],\frac 12\partial_{z}[Q_2(2)V_0]
- \frac 14 B_2Q_2(2)V_0)(\widetilde x)                     \nonumber\\
+ (T^*_{B_1}[\mathcal C_2V_0],\frac 12\partial_{\overline z}
[Q_1(1)U_0]-\frac 14A_1Q_1(1)U_0)(\widetilde x)
\nonumber\\
+2(\partial_zP_{A_2}^*[\mathcal C_1^*U_0]
,\frac 12Q_2(2)V_0)(\widetilde x)+2(\partial_{\overline z}T^*_{B_1}[\mathcal C_2V_0],\frac 12 Q_1(1)U_0)(\widetilde x)\nonumber\\-(P_{A_2}^*(2\partial_z(\mathcal C_1^*U_0)+\partial_{\overline z}
(\mathcal C_0^*U_0)+{\bf B}_1^*U_0+Q_1(2)^*P_{A_2}^*(\mathcal C_1^*U_0)
+ \mathcal C_1^*q_1),\frac 12Q_2(2)V_0)(\widetilde x)\nonumber\\
- (T_{B_1}^*(2\mathcal C_2\partial_{\overline z}V_0+\mathcal C_0\partial_z V_0
+ {\bf B}_2V_0 +2Q_2(1)^*T_{B_1}^* (\mathcal C_2V_0)+\mathcal C_2q_2)
,\frac 12 Q_1(1)U_0)(\widetilde x)\nonumber\\
+ \frak H_{U_0,V_0} (\widetilde x)=0.
\end{eqnarray}
\end{proposition}
{\bf Proof.}
First, using (\ref{poker}), integrating by parts and applying
the stationary phase argument, since the coefficients of the operator
${\bf H}$ are compactly supported,  we obtain the asymptotic formula
\begin{eqnarray}\label{lida}\nonumber
(U_0e^{\tau\Phi}, {\bf H}(x,\partial_z,\partial_{\overline z})
(V_0e^{-\tau\overline \Phi}))_{L^2(\Omega)}                   \nonumber\\
= \int_\Omega \frac 1\tau\left\{(U_0,{\mathcal C}_1\partial^2_{zz}V_0)
-(\partial_{\overline z}(\mathcal C_0^* U_0),\partial_zV_0)
+ (\partial^2_{\overline z\overline z}(\mathcal C_2^*U_0),V_0)
+ (U_0,{\bf B_1}\partial_zV_0)\right.\nonumber\\
\left. -(\partial_{\overline z}({\bf B}_2^*U_0),V_0)
+ (U_0,{\bf B}_0V_0)\right \}e^{\tau(\Phi-\overline\Phi)}dx
= \frak F_{\widetilde x,\tau}(U_0, {\bf H}(x,\partial_z,\partial_{\overline z}
- \tau{\partial_{\bar z}\bar  \Phi})V_0)+o\left(\frac 1\tau\right)\nonumber\\
= \frac{\pi}{2\tau}\frak H_{U_0,V_0} (\widetilde x)+o\left(\frac 1\tau\right)
\quad\mbox{as}\,\,\tau\rightarrow +\infty.
\end{eqnarray}
From Proposition \ref{zanoza} we obtain
\begin{eqnarray}\label{elka}
(U,{\bf H}(x,\partial_z,\partial_{\overline z})V)_{L^2(\Omega)}
= (((U_0-U_1)e^{\tau\Phi}+U_2e^{\tau\overline\Phi}),{\bf H}
(x,\partial_z,\partial_{\overline z})((V_0-V_1)e^{-\tau\overline\Phi}
+ V_2e^{-\tau\Phi}))_{L^2(\Omega)}\nonumber\\+o(\frac 1\tau)\quad \mbox{as}\quad
\tau\rightarrow +\infty.
\end{eqnarray}
Using this equality, (\ref{poker}), (\ref{lida}) and Proposition \ref{zanoza1},
we have
\begin{eqnarray}\label{999}
(U,{\bf H}(x,\partial_z,\partial_{\overline z})V)_{L^2(\Omega)}
= -(U_0e^{\tau\Phi}, {\bf H}(x,\partial_z,\partial_{\overline z})
(V_1e^{-\tau\overline \Phi}))_{L^2(\Omega)}\\-(U_1e^{\tau\Phi},
{\bf H}(x,\partial_z,\partial_{\overline z})(V_0e^{-\tau\overline \Phi}))
_{L^2(\Omega)}+o(\tau)                      \nonumber\\
= -\frak F_{\widetilde x,\tau} (P_{A_2}^*{\bf H}(x,\partial_z
+\tau{\partial_z\Phi},\partial_{\overline z})^*U_0, q_2)
- \frak F_{\widetilde x,\tau}(q_1, T_{B_1}^*{\bf H}(x,\partial_z,\partial
_{\overline z}-\tau{\partial_{\bar z}\bar  \Phi})V_0)\nonumber\\
- \frak J_\tau (P_{A_2}^*{\bf H}(x,\partial_z+\tau{\partial_z\Phi},\partial_{\overline z})
^*U_0, q_2)-\frak J_{\tau}(q_1, T_{B_1}^*{\bf H}(x,\partial_z,\partial
_{\overline z}-\tau{\partial_{\bar z}\bar  \Phi})V_0)+o(\tau)            \nonumber\\
=-4\tau(P_{A_2}^*[\mathcal C_1^*(z-\widetilde z)^2U_0](\widetilde x),
q_2(\widetilde x))-4\tau(T_{B_1}^*
[\mathcal C_2(\overline z-\overline{\widetilde z})^2V_0](\widetilde x),
q_1(\widetilde x))+o(\tau).\nonumber
\end{eqnarray}
Here, in order to obtain the last equality we used the fact that
$\frak J_\tau (P_{A_2}^*{\bf H}(x,\partial_z+\tau{\partial_z\Phi},
\partial_{\overline z})^*U_0, q_2)=o(\tau).$ This follows from (\ref{4.177}).

Next we claim that
\begin{equation}\label{longmon}
P_{A_2}^*[z^k{\mathcal C}_1^*U_0]=z^kP_{A_2}^*[{\mathcal C}_1^*U_0]
\quad\mbox{and}\quad T_{B_1}^*[\overline z^k{\mathcal C}_2V_0]
=\overline z^kT_{B_1}^*[{\mathcal C}_2V_0]\quad \forall k\in \{0,1,2\}.
\end{equation}

Indeed, using the notations $r_{1,k}= P_{A_2}^*[\mathcal C_1^*z^kU_0],
r_{2,k}=z^kP_{A_2}^*[\mathcal C_1^*U_0]$, we observe that
$$
(-\partial_{\overline z} +A_2^*)r_{j,k}=z^k\mathcal C_1^*U_0\quad\mbox{in}
\,\,\Omega.
$$
Hence using (\ref{4.177}) we obtain
$$
(-\partial_{\overline z} +A_2^*)(r_{1,k}-r_{2,k})=0
\quad\mbox{in}\,\,\Omega,\quad (r_{1,k}-r_{2,k})\vert_{\partial\Omega}=0.
$$
By uniqueness of the Cauchy problem for the $\partial_z$ equation,
we obtain that $r_{1,k}-r_{2,k}\equiv 0$. The proof of the second equality
in (\ref{longmon}) is the same.

We introduce the following notations
$$
m_2(x)=P_{A_2}^*[{\mathcal C}_1^*(\partial_z\Phi )^2U_0],\quad
m_1(x)=T_{B_1}^*[{\mathcal C}_2(\partial_{\bar z}\overline{\Phi})^2V_0],
$$
$$
\frak C_1(x)=P_{A_2}^*( 2{\partial_z\Phi}\partial_{ z}(\mathcal C_1^* U_0)
+ 2\mathcal C_1^*U_0+{\partial_z\Phi}\partial_{\overline z}(\mathcal C_0^*U_0)
+ {\partial_z\Phi}{\bf B}_1^*U_0+ Q_1(2)^*P_{A_2}^*(\mathcal C_1^*{\partial_z\Phi}U_0)
+ \mathcal C_1^*{\partial_z\Phi}q_1),
$$
$$
\frak C_2(x)=-T_{B_1}^*(\mathcal C_2 (2{\partial_{\bar z}\bar  \Phi}\partial
_{\overline z}V_0+2V_0)
+ \mathcal C_0{\partial_{\bar z}\bar  \Phi}\partial_{z}V_0
+\partial_{\bar z} \overline \Phi{\bf B}_2V_0+Q_2(1)^*T_{B_1}^*(\mathcal C_2\partial_{\bar z}\overline \Phi V_0)
+ \mathcal C_2{\partial_{\bar z}\bar  \Phi}q_2).
$$

In particular (\ref{longmon}) implies
\begin{equation}\label{longmont1}
m_2(x)=(\partial_z\Phi )^2P_{A_2}^*[{\mathcal C}_1^*U_0]\quad\mbox{and}\quad m_1(x)
=(\partial_{\bar z}\overline{\Phi})^2T_{B_1}^*[{\mathcal C}_2V_0]\quad\mbox{in}\,\,\Omega.
\end{equation}



The formula (\ref{longmont1}) implies the following equalities
\begin{equation}\label{lobster1}
m_2(\widetilde x)=\partial_z m_2(\widetilde x)=0,\quad \partial^2_{zz}
m_2(\widetilde x)= 8P_{A_2^*}[\mathcal C_1^* U_0](\widetilde x),\quad
\partial^3_{zzz} m_2(\widetilde x)
= 24\partial_zP_{A_2^*}[\mathcal C_1^* U_0](\widetilde x),
\end{equation}
\begin{equation}\label{lobster2}
m_1(\widetilde x)
= \partial_{\overline z} m_1(\widetilde x)=0,\quad \partial^2_{\overline z
\overline z} m_1(\widetilde x)= 8T_{B_1^*}[\mathcal C_2 V_0](\widetilde x),
\quad \partial^3_{\overline z\overline z\overline z} m_1(\widetilde x)
= 24\partial_{\overline z}T_{B_1^*}[\mathcal C_2 V_0](\widetilde x).
\end{equation}
Moreover  thanks to our assumption that $q_1(\widetilde x)
= q_2(\widetilde x)=0$, formulae (\ref{lobster1}) and (\ref{lobster2})  implies the
following equalities:
\begin{eqnarray}\label{lobster3}
\frac{1}{32}\partial^4_{zzzz}(m_2,q_2)(\widetilde x)
=\frac{1}{32}\{ 4(\partial^3_{zzz}m_2,\partial_z q_2)
+6(\partial^2_{zz}m_2,\partial^2_{zz} q_2)\}(\widetilde x)         \nonumber\\
= (\partial_zP_{A_2^*}[\mathcal C_1^* U_0],
\partial_zq_2)(\widetilde x)
+ \frac 32(P_{A_2^*}[\mathcal C_1^* U_0],
\partial^2_{zz} q_2)(\widetilde x)
\end{eqnarray}
and
\begin{eqnarray}\label{lobster4}
\frac{1}{32}\partial_{\overline z\overline z\overline z\overline z}^4(m_1,q_1)(\widetilde x)
=\frac{1}{32}\{ 4(\partial^3_{\overline z\overline z\overline z}m_1,\partial_{\overline z} q_1)
+6(\partial^2_{\overline z\overline z}m_1,\partial^2_{\overline z\overline z} q_1)\}(\widetilde x)
                                                   \nonumber\\
= (\partial_{\overline z}T_{B_1^*}[\mathcal C_2 V_0],
\partial_{\overline z}q_1)(\widetilde x)+\frac 32(T_{B_1^*}
[\mathcal C_2 V_0],\partial^2_{\overline z\overline z} q_1)(\widetilde x).
\end{eqnarray}

In addition, since the matrices $\mathcal C_j$ are compactly supported and
the functions $m_j$ satisfy the equations
$$
-(\partial_{\overline z}+A_1^*)m_1=\mathcal C_1^*U_0\quad\mbox{in}\,\,\Omega,
\quad (-\partial_{ z}+B_2^*)m_2=\mathcal C_2V_0
\quad\mbox{in}\,\,\Omega,\quad m_j\vert_{\partial\Omega}=0 \quad
\forall j\in\{1,2\},
$$
we have
\begin{equation}
\frac{\partial^k m_j}{\partial\nu^k}\vert_{\partial\Omega}
= 0\quad\forall k\in\{0,1,2,3\}\,\thinspace j\in \{1,2\}.
\end{equation}
This implies that
\begin{equation}\label{puel1}
\frak J_\tau(m_k)
= o\left(\frac{1}{\tau^3}\right)\quad\mbox{as}\,\,\tau\rightarrow
+\infty,\quad  \forall j\in \{1,2\}.
\end{equation}

Repeating the above arguments  and using (\ref{4.178}) and (\ref{4.179}),
we obtain
\begin{equation}\label{X1}
\frak C_1(x)={\partial_z\Phi}P_{A_2}^*( 2\partial_{ z}(\mathcal C_1^* U_0)
+ \partial_{\overline z}(\mathcal C_0^*U_0) + {\bf B}_1^*U_0+ Q_1(2)^*P_{A_2}^*
(\mathcal C_1^*U_0)+\mathcal C_1^*q_1)+2P_{A_2}^*(\mathcal C_1^*U_0)
\end{equation}
and
\begin{equation}\label{X2}
\frak C_2(x)=-{\partial_{\bar z}\overline\Phi}T_{B_1}^*
(\mathcal C_2 (2\partial_{\overline z}V_0)
+ \mathcal C_0\partial_{z}V_0 +{\bf B}_2V_0+Q_2(1)^*T_{B_1}^*(\mathcal C_2V_0)
+ \mathcal C_2q_2)-2T^*_{B_1}(\mathcal C_2V_0).
\end{equation}
The formulae (\ref{X1}) and (\ref{X2}) imply that
\begin{equation}\label{sneg}
\frak C_1(\widetilde x)=2P_{A_2}^*(\mathcal C_1^*U_0)(\widetilde x),
\quad \frak C_2(\widetilde x)=-2T_{B_1}^*
(\mathcal C_2V_0)(\widetilde x).
\end{equation}
Using the formulae
$$
{\bf H}(x,\partial_z,\partial_{\overline z})(V_0e^{-\tau\overline \Phi})
=e^{-\tau\overline\Phi}( \tau^2{\mathcal C_2}(\partial_{\bar z}{\overline \Phi})^2V_0
-\tau (\mathcal C_2 (2{\partial_{\bar z}\bar  \Phi}\partial_{\overline z}V_0+2V_0)
+ \mathcal C_0{\partial_{\bar z}\bar  \Phi}\partial_{z}V_0 +\partial_{\bar z}{\overline \Phi}{\bf B}_2V_0)
$$
$$
+ {\bf H}(x,\partial_z,\partial_{\overline z})V_0)
$$
and
\begin{eqnarray}
{\bf H}(x,\partial_z,\partial_{\overline z})^*(U_0e^{\tau\Phi})
= e^{\tau\Phi}( \tau^2{\mathcal C_1^*}({\partial_z\Phi})^2U_0 +\tau ( 2{\partial_z\Phi}\partial
_{ z}(\mathcal C_1^*U_0+2\mathcal C_1^*U_0)
+ {\partial_z\Phi}\partial_{\overline z}(\mathcal C_0^*U_0) + {\partial_z\Phi}{\bf B}^*_1U_0)
                                               \nonumber\\
+ {\bf H}(x,\partial_z,\partial_{\overline z})^*U_0),\nonumber
\end{eqnarray}  formula (\ref{elka}) and asymptotics (\ref{lida})
we write the left-hand side of (\ref{lida1}) as
\begin{eqnarray}\label{K999I}
(U,{\bf H}(x,\partial_z,\partial_{\overline z})V)_{L^2(\Omega)}
= \int_\Omega \{ \tau^2(q_1,T_{B_1}^*((\partial_{\bar z}{\overline \Phi})^2\mathcal C_2 V_0))
+ \tau (q_1,\frak C_2)+(q_2,\mathcal M_1)\} e^{\tau(\Phi-\overline\Phi)}dx
                                                         \nonumber\\
+ \int_\Omega \{ \tau^2(q_2,P_{A_2}^*(\mathcal C_1^*({\partial_z\Phi})^2 U_0))
+ \tau (q_2,\frak C_1)+(q_1,\mathcal M_2)\} e^{\tau(\Phi-\overline\Phi)}dx
                                                        \nonumber\\
+  \frac{\pi}{2\tau}\frak H_{U_0,V_0} (\widetilde x)
+ o\left(\frac 1\tau\right)\quad\mbox{as}\,\,\tau\rightarrow +\infty,
\end{eqnarray}
where
$$
\mathcal M_1=P_{A_2}^*({\bf H}(x,\partial_z,\partial_{\overline z})^*U_0)
-P_{A_2}^*Q_1(2)^*P_{A_2}^*((\mathcal C_1^*\partial_z
+\mathcal C_0^*\partial_{\overline z}+\widetilde b)U_0)+P_{A_2}^*
((\mathcal C_1^*(\partial_z+B_1)-\mathcal B_1^*)q_1)
$$
and
$$
\mathcal M_2=T_{B_1}^*({\bf H}(x,\partial_z,\partial_{\overline z})V_0)
-T^*_{B_1}Q_2(1)^*T_{B_1}^*((\mathcal C_0\partial_z
+ \mathcal C_2\partial_{\overline z}+b_1)V_0)
-\mathcal C_0 q_2+T_{B_1}^*((\mathcal C_2(\partial_{\overline z}
+A_2)+{\mathcal B_2})q_2).
$$
By the stationary phase argument
$$
\int_\Omega ((q_1,\mathcal M_1)+(q_2,\mathcal M_2)) e^{\tau(\Phi-\overline\Phi)}dx=o(\frac{1}{\tau})\quad\mbox{as}\quad\tau\rightarrow +\infty.
$$

Then we rewrite (\ref{K999I}) as
 \begin{eqnarray}\label{999I}
(U,{\bf H}(x,\partial_z,\partial_{\overline z})V)_{L^2(\Omega)}
= \int_\Omega \{ \tau^2(q_1,T_{B_1}^*((\partial_{\bar z}{\overline \Phi})^2\mathcal C_2 V_0))
+ \tau (q_1,\frak C_2)\} e^{\tau(\Phi-\overline\Phi)}dx
                                                         \nonumber\\
+ \int_\Omega \{ \tau^2(q_2,P_{A_2}^*(\mathcal C_1^*({\partial_z\Phi})^2 U_0))
+ \tau (q_2,\frak C_1)\} e^{\tau(\Phi-\overline\Phi)}dx
                                                        \nonumber\\
+  \frac{\pi}{2\tau}\frak H_{U_0,V_0} (\widetilde x)
+ o\left(\frac 1\tau\right)\quad\mbox{as}\,\,\tau\rightarrow +\infty.
\end{eqnarray}

Computing the next term in the asymptotics (\ref{lida1}) using the  representation (\ref{999I}), we obtain
\begin{eqnarray}\label{gnom-1}
(U,{\bf H}(x,\partial_z,\partial_{\overline z})V)_{L^2(\Omega)}
=\mathcal I_1(\widetilde x)+\mathcal I_2(\widetilde x)
+ \tau(\frak J_\tau(\frak C_1)+\frak J_\tau({\frak C}_2))+o(1)
\quad\mbox{as}\quad\tau\rightarrow +\infty,\nonumber
\end{eqnarray}
where
\begin{eqnarray}\label{gnom}
\mathcal I_1(x)=
\frac 14\{-\partial^2_{zz}(P_{A_2}^*({\mathcal C}_1^*({\partial_z\Phi})^2U_0), q_2)(x)
+\partial^2_{\overline z\overline z}(P_{A_2}^*({\mathcal C}_1^*
({\partial_z\Phi})^2U_0), q_2)(x)\}+ (\frak C_1, q_2)(x)\nonumber
\end{eqnarray}
and
\begin{eqnarray}\label{gnom0}
\mathcal I_2( x)
= \frac 14\{-\partial^2_{z  z}(T_{B_1}^*
({\mathcal C}_2(\partial_{\bar z}{\overline \Phi})^2V_0), q_1)(x)
+\partial^2_{\bar z \bar z}(T_{B_1}^*({\mathcal C}_2(\partial_{\bar z}{\overline \Phi})^2V_0), q_1)(x)\}
+(\frak C_2, q_1)( x).\nonumber
\end{eqnarray}
By (\ref{4.178}), we have
\begin{equation}\label{puel}
(\frak J_\tau(\frak C_1)+\frak J_\tau(\frak C_2))=o(\frac{1}{\tau^2})
\quad \mbox{as}\,\tau\rightarrow +\infty.
\end{equation}
We claim that $\mathcal I_1(\widetilde x)=\mathcal I_2(\widetilde x)=0$. Indeed,
by (\ref{longmon}), we obtain
$$
\partial^2_{\overline z\overline z}(P_{A_2}^*({\mathcal C}_1^*(\partial_z\Phi)^2U_0), q_2)
(\widetilde x)=0\quad\mbox{and}\quad \partial^2_{zz}(-T_{B_1}^*({\mathcal C}_2
(\partial_{\bar z}{\overline \Phi})^2V_0), q_1)(\widetilde x)=0.
$$
Using the above equality, (\ref{mishka}), (\ref{longmont1}) and (\ref{X1}),
we can compute $\mathcal I_1(\widetilde x)$ as
\begin{eqnarray}\label{gnom1}
\mathcal I_1(\widetilde x)=-\frac 14\partial_{zz}(P_{A_2}^*
(\mathcal C_1^*(\partial_z{\Phi})^2U_0),q_2)(\widetilde x)\nonumber\\
+(P_{A_2}^*(2\mathcal C_1^*U_0), q_2)(\widetilde x)=-2(P_{A_2}^*
(2\mathcal C_1^*U_0), q_2)(\widetilde x)+(P_{A_2}^*
(2\mathcal C_1^*U_0), q_2)(\widetilde x)=0.
\nonumber
\end{eqnarray}
Similarly, using (\ref{longmont1}) and (\ref{X2}),
we compute $\mathcal I_2$  at point $\widetilde x :$
\begin{eqnarray}\label{pida3}
\mathcal I_2(\widetilde x)=\frac 14\partial^2_{\overline z\overline z}
(T_{B_1}^*({\mathcal C}_2(\partial_{\bar z}{\overline \Phi})^2V_0), q_1)(\widetilde x)
-(q_1 ,T_{B_1}^*( {\mathcal C}_2\partial_{\bar z}{\overline \Phi}q_2 ))(\widetilde x)
=\nonumber\\ \frac 14\partial^2_{\overline z\overline z}(T_{B_1}^*({\mathcal C}
_2(\partial_{\bar z}{\overline \Phi})^2V_0), q_1)(\widetilde x)
-2(q_1, T_{B_1}^*(\mathcal C_2V_0))(\widetilde x)=0.\nonumber
\end{eqnarray}

Finally we compute the term of order $\frac{1}{\tau}$ in the asymptotics
of the left-hand side of (\ref{lida1}).
We remind our assumption
\begin{equation}\label{zombi}
q_1(\widetilde x)=q_2(\widetilde x)=0.
\end{equation}
We introduce the couple operators
$$
 L\phi=\frac{1}{4}(-\partial^2_{zz}\phi+\partial^2_{\overline z\overline z}\phi)
(\widetilde x)\quad \mbox{and}\quad \frak P\phi
= (\frac{1}{32}\partial^4_{zzzz} \phi-\frac{1}{16}\partial^4
_{zz\overline z\overline z} \phi
+ \frac{1}{32}\partial^4_{\overline z\overline z\overline z\overline z}\phi)(\widetilde x).
$$
By (\ref{puel}) we see that
\begin{eqnarray}\label{YYY}
(U,{\bf H}(x,\partial_z,\partial_{\overline z})V)_{L^2(\Omega)}
=\frac 1\tau(\frak B(m_1,q_1)+\frak B(m_2,q_2)
 +L (\frak C_2,q_1)+L(\frak C_1,q_2))                     \\
+\frac{\pi}{2\tau}\frak H_{U_0,V_0} (\widetilde x)+\tau (\frak I_\tau(\frak C_1)+\frak I_\tau(\frak C_2))
=\frac 1\tau( \frak B(m_1,q_1)+\frak B(m_2,q_2)\nonumber\\
 +L (\frak C_2,q_1)+L(\frak C_1,q_2))+\frac{\pi}{2\tau}\frak H_{U_0,V_0} (\widetilde x)
+o\left(\frac 1\tau\right)\quad\mbox{as}\,\,\tau\rightarrow +\infty.
                                          \nonumber
\end{eqnarray}

First we compute  $\frak B(m_2,q_2).$
Observe that
\begin{equation}
m_2(\widetilde x)=\partial_{\overline z}m_2(\widetilde x)
= \partial^2_{\overline z\overline z} m_2(\widetilde x)
=\partial^3_{\overline z\overline z\overline z}m_2(\widetilde x)
=\partial^4_{\overline z\overline z\overline z\overline z}m_2(\widetilde x)
=0.
\end{equation}

Then
\begin{equation}\label{Y}
\partial^4_{\overline z\overline z\overline z\overline z}(m_2,q_2)(\widetilde x)=0.
\end{equation}

Using (\ref{zombi}) and (\ref{lobster1}), we obtain
\begin{eqnarray}\label{op7}
\partial^4_{zz\overline z\overline z}(m_2,q_2)(\widetilde x)
=\{(\partial^4_{zz\overline z\overline z}m_2,q_2)
+2(\partial^3_{z\overline z\overline z}m_2,\partial_zq_2)
+ 2(\partial^3_{zz\overline z}m_2,\partial_{\overline z}q_2)
+ (\partial^2_{zz}m_2,\partial^2_{\overline z \overline z}q_2)\nonumber\\
+ (\partial^2_{\overline z\overline z}m_2,\partial^2_{z  z}q_2)
+ 4(\partial^2_{ z\overline z}m_2,\partial^2_{z  \overline z}q_2)
+ 2(\partial_{z}m_2,\partial^2_{z  \overline z\overline z}q_2)
+ 2(\partial_{\overline z}m_2,\partial^2_{z  z\overline z}q_2)
+ (m_2,\partial^4_{zz\overline z\overline z}q_2)\}(\widetilde x)
                                               \nonumber\\
= 2(\partial^3_{zz\overline z}m_2,\partial_{\overline z}q_2)(\widetilde x)
+(\partial^2_{zz}m_2,\partial^2_{\overline z \overline z}q_2)(\widetilde x)
= 16(\partial_{\overline z}P_{A_2}^*[\mathcal C_1^*U_0],\partial_{\overline z}
q_2)(\widetilde x)+(8P_{A_2}^*[\mathcal C_1^*U_0],
\partial^2_{\overline z \overline z}q_2)(\widetilde x).\nonumber\\
\end{eqnarray}

By (\ref{lobster4}), (\ref{Y}), (\ref{op7})

\begin{eqnarray}\label{cloun1}
\frak B(m_2,q_2)=-(\partial_{\overline z}P_{A_2}^*[\mathcal C_1^*U_0],\partial_{\overline z}
q_2)(\widetilde x)-\frac 12(P_{A_2}^*[\mathcal C_1^*U_0],
\partial^2_{\overline z \overline z}q_2)(\widetilde x)\nonumber\\+(\partial_zP_{A_2^*}[\mathcal C_1^* U_0],
\partial_zq_2)(\widetilde x)
+ \frac 32(P_{A_2^*}[\mathcal C_1^* U_0],
\partial^2_{zz} q_2)(\widetilde x).
\end{eqnarray}

Next, in similar way, we compute   $\frak B(m_1,q_1).$
Observe that
\begin{equation}\label{regnum}
m_1(\widetilde x)=\partial_{ z}m_1(\widetilde x)
= \partial^2_{ zz} m_1(\widetilde x)
= \partial^3_{ zzz}m_1(\widetilde x)
= \partial^4_{zz zz}m_1(\widetilde x)=0.
\end{equation}

Then, by (\ref{regnum})
\begin{equation}\label{z}
\partial^4_{ zzzz}(m_1,q_1)(\widetilde x)=0.
\end{equation}
Short computations, (\ref{zombi}), (\ref{lobster2}) provide the formulae
\begin{eqnarray}\label{op3}
\partial^4_{zz\overline z\overline z}(m_1,q_1)(\widetilde x)
= \{(\partial^4_{zz\overline z\overline z}m_1,q_1)
+ 2(\partial^3_{z\overline z\overline z}m_1,\partial_zq_1)
+ 2(\partial^3_{zz\overline z}m_1,\partial_{\overline z}q_1)
+ (\partial^2_{zz}m_1,\partial^2_{\overline z \overline z}q_1)\nonumber\\
+ (\partial^2_{\overline z\overline z}m_1,\partial^2_{z  z}q_1)
+ 4(\partial^2_{ z\overline z}m_1,\partial^2_{z  \overline z}q_1)
+ 2(\partial_{z}m_1,\partial^2_{z  \overline z\overline z}q_1)
+ 2(\partial_{\overline z}m_1,\partial^2_{z  z\overline z}q_1)
+ (m_1,\partial^4_{zz\overline z\overline z}q_1)\}(\widetilde x)\nonumber\\
= 2(\partial^3_{z\overline z\overline z}m_1,\partial_zq_1)(\widetilde x)
+(\partial^2_{\overline z\overline z}m_1,\partial^2_{z  z}q_1)(\widetilde x).
\end{eqnarray}

By (\ref{lobster2}), (\ref{z}), (\ref{op3}) we obtain
\begin{eqnarray}\label{cloun2}
\frak B(m_1,q_1)=-2(\partial^3_{z\overline z\overline z}m_1,\partial_zq_1)(\widetilde x)
-(\partial^2_{\overline z\overline z}m_1,\partial^2_{z  z}q_1)(\widetilde x)\nonumber\\
 +(\partial_{\overline z}T_{B_1^*}[\mathcal C_2 V_0],
\partial_{\overline z}q_1)(\widetilde x)+\frac 32(T_{B_1^*}
[\mathcal C_2 V_0],\partial^2_{\overline z\overline z} q_1) (\widetilde x).
\end{eqnarray}
Using (\ref{zombi}) we have
\begin{eqnarray}\label{op}
L(q_1,\frak C_2)=\frac 14\{-(\partial_{zz}^2 q_1,\frak C_2)(\widetilde x)
-2(\partial_z q_1,\partial_z\frak C_2)(\widetilde x)
+(\partial_{\overline z\overline z}^2 q_1,\frak C_2)(\widetilde x)
+2(\partial_{\overline z} q_1,\partial_{\overline z}\frak C_2)(\widetilde x)\}
                                       \nonumber\\
\end{eqnarray}
and
\begin{eqnarray}\label{op1}
L(q_2,\frak C_1)=\frac 14\{-(\partial_{zz}^2 q_2,\frak C_1)(\widetilde x)
-2(\partial_z q_2,\partial_z\frak C_1)(\widetilde x)
+(\partial_{\overline z\overline z}^2 q_2,\frak C_1)(\widetilde x)
+ 2(\partial_{\overline z} q_2,\partial_{\overline z}\frak C_1)
(\widetilde x)\}.\nonumber\\
\end{eqnarray}

By (\ref{op}), (\ref{op1}), (\ref{cloun1}), (\ref{cloun2}), (\ref{lobster4}) we obtain from (\ref{YYY})
\begin{eqnarray}\label{kaput}
\frac 14(\frak C_1-\frac{1}{4}(\partial^2_{zz}m_2 ,\partial^2
_{\overline z\overline z}q_2))(\widetilde x)\\
+\frac 14 (-\frak C_2-\frac{1}{4}(\partial^2_{\overline z\overline z}m_1,
\partial^2_{ zz}q_1))(\widetilde x)\nonumber\\
+ (\frac 12\partial_{\overline z}\frak C_1-\frac{1}{8}(\partial^3
_{z z\overline  z}m_2 ,\partial_{\overline z}q_2))(\widetilde x)
                               \nonumber\\
+ (-\frac 12\partial_{ z}\frak C_2-\frac{1}{8}(\partial^3
_{z \overline z\overline z}m_1,\partial_{ z}q_1))(\widetilde x)\nonumber\\
+ 3(\partial_zP_{A_2}^*[\mathcal C_1^*U_0],\partial_zq_2)(\widetilde x)+ \frac 32(P_{A_2}^*
[\mathcal C_1^*U_0],\partial^2_{zz}q_2)(\widetilde x)\nonumber\\
+ 3(\partial_{\overline z}T^*_{B_1}[\mathcal C_2V_0]
,\partial_{\overline z}q_1)
(\widetilde x)+ \frac 32(T^*_{B_1}[\mathcal C_2V_0],\partial^2_{\overline z\overline z}q_1)
(\widetilde x)\nonumber\\
+ \frac14 (-\partial^2_{zz}(\frak C_1,q_2)+\partial^2_{\overline z\overline z}
(\frak C_2,q_1))(\widetilde x)\nonumber
+\frak H_{U_0,V_0}(\widetilde x)=0.
\end{eqnarray}
By (\ref{X1}) and (\ref{lobster1}), we have
\begin{equation}\label{kauk1}
\frak C_1(\widetilde x)-\frac{1}{4}\partial^2_{zz}m_2 (\widetilde x)=0.
\end{equation}
By (\ref{X2}) and (\ref{lobster2}), we see
\begin{equation}\label{kauk2}
\frak C_2(\widetilde x)+\frac{1}{4}\partial^2_{\overline z\overline z}
m_1(\widetilde x)=0.
\end{equation}
Applying (\ref{X1}) and (\ref{lobster1}), we obtain
\begin{equation}\label{kauk3}
\left(\frac 12\partial_{\overline z}\frak C_1-\frac{1}{8}\partial^3
_{z z\overline  z}m_2 \right)(\widetilde x)=0.
\end{equation}
Applying (\ref{X2}) and (\ref{lobster2}), we obtain
\begin{equation}\label{kauk4}
\left(-\frac 12\partial_{ z}\frak C_2-\frac{1}{8}
\partial^3_{z \overline z\overline z}m_1\right)(\widetilde x)=0.
\end{equation}
By (\ref{kauk1})-(\ref{kauk4}), we rewrite (\ref{kaput}) as
\begin{eqnarray}\label{kaput1}
\{3(\partial_zP_{A_2}^*[\mathcal C_1^*U_0]
,\partial_zq_2)+ \frac 32(P_{A_2}^*
[\mathcal C_1^*U_0],\partial^2_{zz}q_2)\}(\widetilde x)          \\
  +\{3(\partial_{\overline z}T^*_{B_1}[\mathcal C_2V_0],\partial_{\overline z}q_1)
+ \frac 32(T^*_{B_1}[\mathcal C_2V_0],\partial^2_{\overline z\overline z}q_1)\}(\widetilde x)
\nonumber\\
+ \frac14 (-\partial^2_{zz}(\frak C_1,q_2)+\partial^2_{\overline z\overline z}
(\frak C_2,q_1))(\widetilde x)\nonumber
+ \frak H_{U_0,V_0}(\widetilde x)=0.
\end{eqnarray}
Now we compute the last term in (\ref{kauk4}):
\begin{eqnarray}\label{kaput2}
\frac14 (-\partial^2_{zz}(\frak C_1,q_2)+\partial^2_{\overline z\overline z}
(\frak C_2,q_1))(\widetilde x)=\frac14 (-2(\partial_z\frak C_1,\partial_z q_2)
- (\frak C_1,\partial^2_{z z} q_2)\\+2(\partial_{\overline z}\frak C_2,
\partial_{\overline z}q_1)+(\frak C_2,\partial^2_{\overline z\overline z}q_1))
(\widetilde x)
= \frac14 (-2(\partial_z\frak C_1,\partial_z q_2)-2(P_{A_2}^*
[\mathcal C_1^*U_0],\partial^2_{ zz} q_2)\nonumber\\+2(\partial
_{\overline z}\frak C_2,\partial_{\overline z}q_1)-2(T^*_{B_1}
(\mathcal C_2 V_0),\partial^2_{\overline z\overline z}q_1))(\widetilde x)     \nonumber\\
= \frac14 (-4(\partial_zP_{A_2}^*[\mathcal C_1^* U_0],
\partial_zq_2)(\widetilde x)-2(P_{A_2}^*[\mathcal C_1^*U_0],\partial^2_{zz} q_2)
(\widetilde x)\nonumber\\-4(\partial_{\overline z}T_{B_1}^*[\mathcal C_2 V_0]
, \partial_{\overline z}q_1)(\widetilde x)-2(T^*_{B_1}
(\mathcal C_2 V_0)(\widetilde x),\partial^2_{\overline z\overline z}q_1)(\widetilde x))
                            \nonumber\\
- 4(P_{A_2}^*(2\partial_z(\mathcal C_1^*U_0)+\partial_{\overline z}
(\mathcal C_0^*U_0)+{\bf B}_1^*U_0+Q_1(2)^*P_{A_2}^*(\mathcal C_1^*U_0)
+ \mathcal C_1^*q_1),\partial_zq_2)(\widetilde x)\nonumber\\
- 4(T_{B_1}^*(2\mathcal C_2\partial_{\overline z}V_0
+ \mathcal C_0\partial_z V_0+{\bf B}_2V_0
+ 2Q_2(1)^*T_{B_1}^* (\mathcal C_2V_0)+\mathcal C_2q_2),
\partial_{\overline z}q_1)(\widetilde x)).\nonumber
\end{eqnarray}
From (\ref{kaput1}) and (\ref{kaput2}), we obtain
\begin{eqnarray}\label{kaput33}
2(\partial_zP_{A_2}^*[\mathcal C_1^*U_0]
,\partial_zq_2)(\widetilde x)+2(\partial_{\overline z}T^*_{B_1}[\mathcal C_2V_0],\partial_{\overline z}q_1)(\widetilde x)\nonumber\\
+(P_{A_2}^*[\mathcal C_1^*U_0],\partial^2_{zz}q_2)(\widetilde x)
+ (T^*_{B_1}[\mathcal C_2V_0],\partial^2_{\overline z\overline z}q_1)(\widetilde x)\\
-(P_{A_2}^*(2\partial_z(\mathcal C_1^*U_0)+\partial_{\overline z}
(\mathcal C_0^*U_0)+{\bf B}_1^*U_0+Q_1(2)^*P_{A_2}^*(\mathcal C_1^*U_0)
+ \mathcal C_1^*q_1),\partial_z q_2)(\widetilde x)\nonumber\\
- (T_{B_1}^*(2\mathcal C_2\partial_{\overline z}V_0+\mathcal C_0\partial_z V_0
+ {\bf B}_2V_0 +2Q_2(1)^*T_{B_1}^* (\mathcal C_2V_0)+\mathcal C_2q_2),
\partial_{\overline z} q_1)(\widetilde x) +
\frak H_{U_0,V_0}(\widetilde x)=0.\nonumber
\end{eqnarray}

Observe that by (\ref{mishka})
$$
2\partial_{\overline z} q_1+A_1q_1=Q_1(1)U_0\quad\mbox{and}\quad
2\partial_z q_2+B_2q_2=Q_2(2)V_0.
$$
Then using (\ref{zombi}), we have
\begin{equation}\label{elka}
\quad \partial_{\overline z} q_1(\widetilde x)
= \frac 12 Q_1(1)U_0(\widetilde x)\quad
\mbox{and}\quad \partial_z q_2(\widetilde x)=\frac 12Q_2(2)V_0(\widetilde x).
\end{equation}
Taking into account that
$$
2\partial^2_{\overline z\overline z} q_1+\partial_{\overline z} A_1q_1
+ A_1 \partial_{\overline z} q_1=\partial_{\overline z}[Q_1(1)U_0],
\quad 2\partial^2_{zz} q_2+\partial_z B_2q_2 +B_2 \partial_z q_2=\partial_{z}
[Q_2(2)V_0],
$$
we have
\begin{equation}\label{elka1}
\quad \partial^2_{zz} q_2(\widetilde x)
= \frac 12\partial_{z}[Q_2(2)V_0](\widetilde x)
-\frac 14 B_2Q_2(2)V_0(\widetilde x),
\quad
\partial^2_{\overline z\overline z} q_1(\widetilde x)=\frac 12\partial_{\overline z
}[Q_1(1)U_0](\widetilde x)-\frac 14A_1Q_1(1)U_0(\widetilde x).
\end{equation}
Using (\ref{elka}) and (\ref{elka1}), we rewrite (\ref{kaput33}) as
\begin{eqnarray}\label{kaput4}
(P_{A_2}^*[\mathcal C_1^*U_0],\frac 12\partial_{z}[Q_2(2)V_0]
- \frac 14 B_2Q_2(2)V_0)(\widetilde x)                     \nonumber\\
+ (T^*_{B_1}[\mathcal C_2V_0],\frac 12\partial_{\overline z}
[Q_1(1)U_0]-\frac 14A_1Q_1(1)U_0)(\widetilde x)
\nonumber\\
+2(\partial_zP_{A_2}^*[\mathcal C_1^*U_0]
,\frac 12Q_2(2)V_0)(\widetilde x)+2(\partial_{\overline z}T^*_{B_1}[\mathcal C_2V_0],\frac 12 Q_1(1)U_0)(\widetilde x)\nonumber\\-(P_{A_2}^*(2\partial_z(\mathcal C_1^*U_0)+\partial_{\overline z}
(\mathcal C_0^*U_0)+{\bf B}_1^*U_0+Q_1(2)^*P_{A_2}^*(\mathcal C_1^*U_0)
+ \mathcal C_1^*q_1),\frac 12Q_2(2)V_0)(\widetilde x)\nonumber\\
- (T_{B_1}^*(2\mathcal C_2\partial_{\overline z}V_0+\mathcal C_0\partial_z V_0
+ {\bf B}_2V_0 +2Q_2(1)^*T_{B_1}^* (\mathcal C_2V_0)+\mathcal C_2q_2)
,\frac 12 Q_1(1)U_0)(\widetilde x)\nonumber\\
+ \frak H_{U_0,V_0} (\widetilde x)=0.
\end{eqnarray}
The proof of the proposition is complete.
$\blacksquare$

{\bf Remark 4.1.}{\it The equation (\ref{kaput3}) at each point $\widetilde x$
depends on the choice of the functions $q_1$ and $q_2$ since they are
supposed to satisfy the condition
$q_1(\widetilde x)=q_2(\widetilde x)=0.$  Of course the choice of such functions
for any fixed point $\widetilde x$ is unique only modulo a function $q_{1,*}\in
\mbox{Ker} \,(2\partial_{\bar z}+A_1) \, $ and
$q_{2,*}\in \mbox{Ker}\,(2\partial_{ z}+B_2)$. On the other hand,
any function $q_{1,*}\in \mbox{Ker}\,(2\partial_{\bar z}+A_1)$
can be represented in the form $q_{1,*}=zq_{11}, q_{11}\in\mbox{Ker}
\,(2\partial_{\bar z}+A_1)$ and any function $q_{2,*}\in\mbox{Ker}
\,(2\partial_{z}+B_2)$ can be represented in the form
$q_{2,*}=\bar zq_{22}, q_{22}\in\mbox{Ker} \,(2\partial_{z}+B_2).$
Therefore by (\ref{longmont1}) as long as the point $\widetilde x$ is fixed the
choice of the functions $q_j$ does not affect the equation (\ref{kaput3}). }

We complete this section, presenting one of many possible choices
of functions $q_j$ which will later in Sections 5 and 6.
Let us fix some point $x^0$ in $\Omega$ and consider a ball $B(x^0,\delta)$
centered at $x^0$ of the small positive radius $\delta.$
By Proposition \ref{nikita}, there exist regular functions $q_{k,j}$
such that
$$
q_{k,j}(x^0)=\vec e_j, \quad
q_{1,j}\in \mbox{Ker}\,(2\partial_{\bar z} +A_1), \quad q_{2,j}\in
\mbox{Ker}\,(2\partial_{ z} +B_2),  \quad k=1,2, \thinspace j\in\{1,2,3\}.
$$
Then, provided that $\delta>0$ is sufficiently small, there exist functions
$r_{k,j}(\widetilde x)$ such that
$$
\sum_{j=1}^3r_{k,j}(\widetilde x)q_{k,j}(\widetilde{x})
=-q_k^0(\widetilde x) \quad \forall \widetilde x\in B(x^0,\delta),\,\,\mbox{where}
\,\,\quad q_2^0=T_{B_2}(Q_2(2)V_0),\quad q_1^0=P_{A_1}(Q_1(1)U_0).
$$

Then we set
\begin{equation}\label{victoryz}
q_k(x) := q_k(x,\widetilde x)
= q_k^0(x)+\sum_{j=1}^3r_{k,j}(\widetilde x) q_{k,j}(x), \quad k=1,2.
\end{equation}

By (\ref{4.178}), (\ref{4.179}) we note that
\begin{equation}\label{boundary}
P_{A_2}^*(\mathcal C_1^*q_{1,j})\vert_{\partial\Omega}=T_{B_1}^*(\mathcal C_2q_{2,j})\vert_{\partial\Omega} =0\quad\forall j\in\{1,2,3\}.
\end{equation}
Using the boundary conditions (\ref{boundary}) and using Proposition
\ref{inga} we obtain
\begin{equation}\label{!zanoza1!}
\left\Vert s\,\root\of{\phi_s} P_{A_2}^*
\left(\mathcal C_1^*\sum_{j=1}^3r_{1,j}(\widetilde x)
q_{1,j}\right)e^{s\phi_s}\right\Vert_{L^2(\Omega)}
\le C\Vert \mathcal C_2^*e^{s\phi_s}\Vert_{L^2(\Omega)}
\end{equation}
and
\begin{equation}\label{!zanoza2!}
\left\Vert s\,\root\of{\phi_s} T_{B_1}^*
\left(\mathcal C_2\sum_{j=1}^3r_{2,j}
(\widetilde x)q_{2,j}\right)e^{s\phi_s}\right\Vert_{L^2(\Omega)}
\le C\Vert \mathcal C_2^*e^{s\phi_s}\Vert_{L^2(\Omega)}.
\end{equation}

\section{Construction of complex geometric optics solutions for
the Navier-Stokes equations}\label{sec2}

Let ${\bf u}=(u_1,u_2),p$ be a solution to the Stokes equations
\begin{equation}\label{stokes}
L_\mu(x,D)(\mbox{\bf u},
p)=(\sum_{j=1}^2\partial_{x_j}(\mu(\partial_{x_j}
u_1+\partial_{x_1}u_j))+\partial_{x_1}p,\sum_{j=1}^2\partial_{x_j}
(\mu(\partial_{x_j}
u_2+\partial_{x_2}u_j))+\partial_{x_2}p)=0.
\end{equation}

We construct the complex geometric optics  solution to the Stokes equations.

As the first step of such a construction we reduce the Stokes
equations to a decoupled elliptic system.

\begin{proposition}\label{gigi} Let functions  ${\bf w}=(w_1,w_2),f $
be some  solutions to the elliptic system
\begin{equation}\label{-11}
\Delta f=\mbox{div} \,{\bf w}
\quad\mbox{in}\,\,\Omega,
\end{equation}
\begin{equation}
\mu\Delta w_1+2\mu_{x_1}\partial_{x_1}w_1+\partial_{x_2}\mu(\partial_{x_2}
w_1+\partial_{x_1}w_2)+\partial_{x_1}\mu\mbox{div}\, {\bf w}+2\partial^2_{x_1
x_1}\mu\partial_{x_1}f+2\partial^2_{x_2x_1}\mu\partial_{x_2}
f=0\quad\mbox{in}\,\,\Omega,
\end{equation}
\begin{equation}\label{-12}
\mu\Delta w_2+2\partial_{x_2}\mu\partial_{x_2}w_2+\partial_{x_1}
\mu(\partial_{x_2}
w_1+\partial_{x_1}w_2)+\partial_{x_2}\mu\mbox{div}\,
{\bf w}+2\partial^2_{x_2x_2}\mu\partial_{x_2}f+2\partial^2_{x_1x_2}
\mu\partial_{x_1}
f=0\quad\mbox{in}\,\,\Omega.
\end{equation}
Then the pair
$({\bf u},p)=({\bf w}-\nabla f,\frac 12 \mu\Delta
f+\partial_{x_1}\mu \partial_{x_1}f+\partial_{x_2}\mu\partial_{x_2}f)$
solves the Stokes equations (\ref{stokes}).
\end{proposition}

{\bf Proof.}
The equation (\ref{-11}) implies $\mbox{div}\, {\bf u}=\mbox{div}\,
({\bf w}-\nabla f)=\mbox{div}\, {\bf w}-\Delta  f=0.$
Short computations provide
$$
\mathcal J_k:=\frac 12\sum_{j=1}^2\partial_{x_j}(\mu(\partial_{x_j}
(w_k-\partial_{x_k}f)+\partial_{x_k}(w_j-\partial_{x_j}f)))+\partial_{x_k}p
$$
$$
= \frac 12\sum_{j=1}^2\partial_{x_j}(\mu(\partial_{x_j}w_k
-\partial^2_{x_kx_j}f+\partial_{x_k}w_j-\partial^2_{x_jx_k}f))+\partial_{x_k}p
$$
$$
= \frac 12\sum_{j=1}^2\partial_{x_j}\mu(\partial_{x_j}w_k
-\partial^2_{x_kx_j}f+\partial_{x_k}w_j-\partial^2_{x_jx_k}f)
$$
$$
+ \frac 12\{\mu(\partial^2_{x_jx_j}w_k-\partial^3_{x_kx_jx_j}f
+ \partial_{x_kx_j}w_j-\partial^3_{x_jx_jx_k}f)\}+\partial_{x_k}p
$$
$$
= \frac 12\sum_{j=1}^2\partial_{x_j}\mu(\partial_{x_j}w_k
 -\partial^2_{x_kx_j}f+\partial_{x_k}w_j-\partial^2_{x_jx_k}f)
$$
$$
+ \frac 12\{\mu \Delta w_k-\mu\partial_{x_k}\Delta f
+\mu\partial_{x_k}\mbox{div}\,{\bf w}-\mu\partial_{x_k}\Delta f\}
+ \partial_{x_k}p.
$$
The equation (\ref{-11}) yields
$$
-\mu\partial_{x_k}\Delta f+\mu\partial_{x_k}\mbox{div}\,{\bf w}\equiv 0.
$$
Hence
$$
\mathcal J_k=\frac 12\sum_{j=1}^2\partial_{x_j}\mu(\partial_{x_j}w_k-\partial^2_{x_kx_j}f+\partial_{x_k}w_j-\partial^2_{x_jx_k}f)
$$
$$
+ \frac 12\{\mu \Delta w_k-\mu\partial_{x_k}\Delta f\}+\partial_{x_k}p
$$
$$
= \frac 12\sum_{j=1}^2\partial_{x_j}\mu(\partial_{x_j}w_k-\partial^2_{x_kx_j}f
+ \partial_{x_k}w_j-\partial^2_{x_jx_k}f)
$$
$$
+ \frac 12\{\mu \Delta w_k-\partial_{x_k}(\mu\Delta f)
+ \partial_{x_k}\mu\Delta f\}+\partial_{x_k}p
$$
$$
= \frac 12\sum_{j=1}^2\partial_{x_j}\mu(\partial_{x_j}w_k
-\partial^2_{x_kx_j}f+\partial_{x_k}w_j-\partial^2_{x_jx_k}f)
$$
$$
+ \frac 12\{\mu \Delta w_k-\partial_{x_k}(\mu\Delta f)
+\partial_{x_k}\mu\mbox{div}\, {\bf w}\}+\partial_{x_k}p
$$
$$
= \frac 12\sum_{j=1}^2(\partial_{x_k}(-2\partial_{x_j}\mu\partial_{x_j}f)
+2\partial_{x_jx_k}\mu\partial_{x_j}f)
$$
$$
+ \frac 12\{\mu \Delta w_k-\partial_{x_k}(\mu\Delta f)
+\partial_{x_k}\mu\mbox{div}\, {\bf w}+\sum_{j=1}^2\partial_{x_j}
\mu(\partial_{x_j}w_k+\partial_{x_k}w_j)\}+\partial_{x_k}p=0.
$$
Thus the proof of the proposition is complete.
$\blacksquare$

We set
\begin{equation}\label{ii}
{\bf H}(x,\partial_z,\partial_{\overline z})
\left (\begin{matrix}{\bf v}\\g\end{matrix}\right )
\end{equation}
$$
=\left(
\begin{matrix}
\frac{\mu}{\mu_2}\mathcal M_2(x,D) ({\bf v},g)- 2(\nabla\mu)
\mbox{div}\,{\bf  v}-\left( \begin{matrix} (\nabla \mu,\nabla v_1)\\
(\nabla \mu,\nabla v_2) \end{matrix} \right)
-2\left( \begin{matrix} (\nabla \partial_{x_1}\mu,\nabla g)\\
(\nabla \partial_{x_2}\mu,\nabla g) \end{matrix} \right)\\
\mbox{div}\, \left (\frac{\mu}{\mu_2}\mathcal M_2 (x,D)({\bf v},g)
- 2(\nabla\mu)\mbox{div}\,{\bf  v}
-\left( \begin{matrix} (\nabla \mu,\nabla v_1)\\
(\nabla \mu,\nabla v_2) \end{matrix} \right )-2\left( \begin{matrix}
(\nabla \partial_{x_1}\mu,\nabla g)\\(\nabla \partial_{x_2}\mu,\nabla g)
\end{matrix} \right) \right )\end{matrix}\right),
$$
where the operator $\mathcal M_2(x,D)$ is given by
\begin{equation}\label{ii1}
\mathcal M_2(x,D) ({\bf v},g)
\end{equation}
$$
=\left (\begin{matrix}
2\partial_{x_1}\mu_{2}\partial_{x_1}v_1+\partial_{x_2}\mu_{2}(\partial_{x_2}
v_1+\partial_{x_1}v_2)+\partial_{x_1}\mu_{2}\mbox{div}\, {\bf v}
+2\partial^2_{x_1x_1}\mu_{2}\partial_{x_1}g
+2\partial^2_{x_1x_2}\mu_{2}\partial_{x_2}g\\
2\partial_{x_2}\mu_{2}\partial_{x_2}v_2+\partial_{x_1}\mu_{2}(\partial_{x_2}
v_1+\partial_{x_1}v_2)+\partial_{x_2}\mu_{2}\mbox{div}\,
{\bf v}+2\partial^2_{x_2x_2}\mu_{2}\partial_{x_2}g+2\partial^2_{x_1x_2}
\mu_{2}\partial_{x_1}g\end{matrix}\right ).
$$

We note that the coefficients of ${\bf H}(x,\partial_z,\partial_{\overline z})$
are given by the derivatives of $\mu$.

\begin{proposition}\label{vo2}
Let $({\bf w} ,f), ({\bf v},g)$ be some regular  solutions to the
system (\ref{-11})-(\ref{-12}) with coefficients $\mu_1$ and $\mu_2$
respectively.   Provided that the function $\mu=\mu_1-\mu_2$ has a compact support in domain $\Omega$ we have
\begin{equation}\label{po}
({\bf w}-\nabla f, L_\mu(x,D)({\bf v}-\nabla g,0))_{L^2(\Omega)}
=(({\bf w},f),{\bf H}(x,\partial_z,\partial_{\overline z}) ({\bf v},g))
_{L^2(\Omega)}.
\end{equation}
\end{proposition}
{\bf Proof.} The short computations imply
$$
I = \Biggl({\bf w}-\nabla f,\bigg [ -\sum_{j=1}^2\partial_{x_j}(\mu(\partial_{x_j}
(v_1-\partial_{x_1}g)+\partial_{x_1}(v_j-\partial_{x_j} g))),
$$
$$
- \sum_{j=1}^2\partial_{x_j}(\mu(\partial_{x_j}(v_2-\partial_{x_2}g)
+\partial_{x_2}( v_j-\partial_{x_j} g)))\bigg ]\Biggr)_{L^2(\Omega)}
:= I_1+I_2,
$$
where we set
$$
I_1 = \left({\bf w}-\nabla f,\bigg [
-\sum_{j=1}^2\partial_{x_j}(\mu(\partial_{x_j}v_1+\partial_{x_1} v_j)),
-\sum_{j=1}^2\partial_{x_j}(\mu(\partial_{x_j}v_2+\partial_{x_2} v_j))\bigg ]
\right)_{L^2(\Omega)}
$$
and
$$
I_2 = \left({\bf w}-\nabla f,\bigg [
\sum_{j=1}^2\partial_{x_j}(\mu(\partial_{x_jx_1}g+\partial_{x_1x_j} g)),
\sum_{j=1}^2\partial_{x_j}(\mu(\partial_{x_jx_2}g+\partial_{x_2x_j} g))\bigg ]
\right)_{L^2(\Omega)}.
$$
Observe that
$$
I_2
= \left({\bf w}-\nabla f,2\nabla\sum_{j=1}^2\frac{\partial}{\partial x_j}
(\mu\partial_{x_j}g)\right)_{L^2(\Omega)}
$$
$$
- ({\bf w}-\nabla f,(2(\nabla \partial_{x_1}\mu,\nabla g)+2\partial_{x_1}\mu
\Delta g,(2\nabla \partial_{x_2}\mu,\nabla g)+2\partial_{x_2}\mu \Delta g))
_{L^2(\Omega)}.
$$
Since the vector field ${\bf w}-\nabla f$ is divergence free and
$\sum_{j=1}^2\frac{\partial}{\partial x_j}(\mu\partial_{x_j}g)\vert
_{\partial\Omega}=0$, we have
$$
I_2=-({\bf w}-\nabla f,(2(\nabla \partial_{x_1}\mu,\nabla g)+2\partial_{x_1}
\mu \Delta g,2(\nabla \partial_{x_2}\mu,\nabla g)+2\partial_{x_2}\mu \Delta g))
_{L^2(\Omega)}.
$$
Using the equation (\ref{-11}), we obtain
\begin{equation}\label{begemot}
I_2=-({\bf w}-\nabla f,(2(\nabla \partial_{x_1}\mu,\nabla g)+2\partial_{x_1}\mu \mbox{div}\, {\bf v},2(\nabla \partial_{x_2}\mu,\nabla g)*
+2\partial_{x_2}\mu \mbox{div}\,{\bf v} ))_{L^2(\Omega)}.
\end{equation}

Since the vector field ${\bf w}-\nabla f$ is divergence free and
$\mu\mbox{div}\,{\bf v}\vert_{\partial\Omega}=0$, we have
$$
I_1=({\bf w}-\nabla f,(-\mu  \Delta v_1, -\mu\Delta v_2)-((\nabla \mu,\nabla v_1),(\nabla \mu,\nabla v_2)))_{L^2(\Omega)}
$$
$$
+({\bf w}-\nabla f,\nabla (\mu \mbox{div}\, {\bf v}))_{L^2(\Omega)}
-({\bf w}-\nabla f,\nabla \mu \mbox{div}\, {\bf v})_{L^2(\Omega)}
$$
$$
= ({\bf w}-\nabla f,(-\mu  \Delta v_1,- \mu\Delta v_2)
-((\nabla \mu,\nabla v_1),(\nabla \mu,\nabla v_2)))_{L^2(\Omega)}
-({\bf w}-\nabla f,\nabla \mu \mbox{div}\,{\bf  v})_{L^2(\Omega)}
$$
$$
= ({\bf w}-\nabla f,\mathcal M_2 (x,D)({\bf v},g))_{L^2(\Omega)}
-({\bf w}-\nabla f,((\nabla \mu,\nabla v_1),(\nabla \mu,\nabla v_2)))
_{L^2(\Omega)}
$$
$$
- ({\bf w}-\nabla f,\nabla \mu \mbox{div}\, {\bf v})_{L^2(\Omega)}.
$$
Combining this equality with (\ref{begemot}), we obtain (\ref{po}).
$\blacksquare$

 We introduce the following matrices:
\begin{eqnarray}\label{yokumoku1}
{\bf C}_1=\frac 12\left( \begin{matrix}  \frac{3\partial_{x_1}\mu_1}{\mu_1} &
\frac{\partial_{x_2}\mu_1}{\mu_1} & \frac{2\partial^2_{x_1 x_1}\mu_1}{\mu_1}
\\
\frac{\partial_{x_2}\mu_1}{\mu_1}&\frac{\partial_{x_1}\mu_1}{\mu_1}&\frac{2\partial^2_{x_1 x_2}\mu_1}{\mu_1}\\
1 & 0 & 0
\end{matrix}\right),{\bf C}_2=\frac 12\left( \begin{matrix}
\frac{\partial_{x_2}\mu_1}{\mu_1} & \frac{\partial_{x_1}\mu_1}{\mu_1} & \frac{2\partial^2_{x_1 x_2}\mu_1}{\mu_1} \\
\frac{\partial_{x_1}\mu_1}{\mu_1}&\frac{3\partial_{x_2}\mu_1}{\mu_1}&
\frac{2\partial^2_{x_2 x_2}\mu_1}{\mu_1}
\\
0 & 1 & 0
\end{matrix}\right),\nonumber\\\quad  A_1=({\bf C}_1+i{\bf C}_2),\quad
B_1=({\bf C}_1-i{\bf C}_2).
\end{eqnarray}

By Proposition \ref{gigi}, the system (\ref{-11})-(\ref{-12})
can be written in the form (\ref{-10}) with $A_1, B_1$ and $C_1=0$.
Therefore, by Proposition \ref{gigi} we can construct the
complex geometric optics solution $u_\tau$ to this system in the form
\begin{equation}\label{vo1}
\mbox{ u}_\tau=\left(\begin{matrix}w_1\\w_2\end{matrix}\right )-\nabla f,
\quad U=\left (\begin{matrix} w_1\\w_2\\f\end{matrix}\right),\quad
U=e^{\tau\Phi}( U_0-U_1)+\sum_{j=2}^\infty(-1)^jU_j e^{\tau\overline \Phi},
\end{equation}
where the matrices $A_1,B_1$ are defined in (\ref{yokumoku1}) and
the function $U_j$ determined in Proposition \ref{mursilka0}.

Next we construct the complex geometric optics solution for the
stationary Navier-Stokes equations.

The stationary Navier-Stokes equations can be written in the form
\begin{equation}\label{viking}
G_\mu(\mbox{\bf u},{ p})=L_\mu(x,D)(\mbox{\bf u},
p)+(\mbox{\bf u},\nabla)\mbox{\bf u}=0.
\end{equation}

We construct the complex geometric optics solution
for the Navier-Stokes equations using the
Newton-Kantorovich iteration scheme. More precisely we use Theorem 6
(1.XVIII) from \cite{KA} p.708.

We recall that by (\ref{nol})
\begin{equation}\label{inequality}
\varphi(x)<- 1 <0\quad\forall x\in\Omega.
\end{equation}
It is convenient for us to change the unknown function $({\bf u}, p)=({\bf r},
\mbox{q})+ (\mbox{u}_\tau, \mbox{p}_\tau).$
By (\ref{viking}) the pair $({\bf r}, q)$ satisfies the  equation
\begin{equation}\label{gemoroi}L_\mu(x,D)(\mbox{\bf r},\mbox{
q})+(\mbox{\bf r}+\mbox{u}_\tau,\nabla)(\mbox{\bf r}+\mbox{u}_\tau)=0.
\end{equation}

Denote $X=\{({\bf u},p)\in W_2^{2,\tau}(\Omega)\times
L^2(\Omega)\vert \thinspace \mbox {div}\,{\bf u}=0\}$ with the norm\newline
$\Vert ({\bf u},p)\Vert_X=\Vert e^{-\tau\varphi} ({\bf u},p)\Vert
_{W_2^{2,\tau}(\Omega)\times L^2(\Omega)}.$  Let
$L_\mu(x,D)^{-1}$ be the operator from $L^2_{e^{-\tau\varphi}}(\Omega)$
into orthogonal complement of $\mbox{Ker}\,L_\mu(x,D)$ in $X.$
Applying to both sides of equation (\ref{gemoroi}) the operator
$L_\mu(x,D)^{-1}$, we obtain
$$
P(\mbox{\bf r},\mbox{q}) : = (\mbox{\bf r},\mbox{
q})+L_\mu(x,D)^{-1}(\mbox{\bf r}+\mbox
{u}_\tau,\nabla)(\mbox{\bf r}
+\mbox {u}_\tau)=0.
$$
The operator $P$ is twice continuously differentiable as the
mapping from $X$ into $X.$ By Proposition \ref{vanka}, we have
\begin{equation}\label{i02}
\Vert \Gamma_0\Vert_{\mathcal L(X;X)}\le C\tau^2\quad\forall \tau\ge 1.
\end{equation}

We set ${\bf x}_0=(0,0)$ and $\Gamma_0=[P'({\bf x}_0)]^{-1}.$ From this
inequality and (\ref{inequality}) we have
$$
\mbox{sup}_{x\in\Omega_0}\Vert \Gamma_0 P({\bf x}_0)\Vert_X
\le C\tau^2e^{-2\tau}=\eta(\tau)\quad\forall \tau\ge 1.
$$
Here $P'({\bf x}_0)$ denotes the Fr\'echet derivative at ${\bf x}_0$.
We set $\Omega_0=\{{\bf x}\vert\Vert {\bf x}-{\bf x}_0\Vert_X\le r_0\}.$
By (\ref{i02}), we have
$$
\Vert \Gamma_0P''({\bf x})\Vert_X\le C\tau^2=K(\tau).
$$
Then $h=K\eta\le \tau^4e^{-2\tau}$ and
$r_0(\tau)=\frac{1-\root\of{1-2h}}{h}\eta\le 2\tau^2e^{-2\tau}<\frac
12 $ for all sufficiently large $\tau.$ Then there exists a solution
${\bf x}_*$ to the equation  $P({\bf x})=0$ such that
$\Vert{\bf  x}_*\Vert_X \le r_0(\tau).$
The construction of the complex geometric optics solution for the
Navier-Stokes equations completed.

Then a solution to the Navier-Stokes can be represented in the form
\begin{equation}\label{voz}
\mbox{\bf u}=\mbox{u}_\tau+u_{cor}, \quad \mbox{where} \quad
\Vert u_{cor}\Vert_{W_2^{2,\tau}(\Omega)}=o\left(\frac 1\tau\right)
\quad \mbox{as}\,\,\tau\rightarrow +\infty.
\end{equation}

\section{Completion of the proof of Theorem \ref{vokal}}\label{sec4}

Suppose that for positive smooth functions $\mu_1, \mu_2$
the Dirichlet-to-Neumann maps (\ref{z1}) are the same.
Then
\begin{equation}\label{boleslav}
\frac{\partial^\ell \mu_1}{\partial \nu^\ell}=\frac{\partial^\ell \mu_2}
{\partial \nu^\ell}\quad\mbox{on}\quad\partial\Omega, \quad\forall \vert
\ell\vert\le 10.
\end{equation}
Since the domain $\Omega$ is assumed to be bounded, there exists
a ball $B(0,r)$ such that $\overline \Omega\subset B(0,r).$
Thanks to (\ref{boleslav}), we can extend the coefficients $\mu_j$
into $B(0,r)\setminus\Omega$ in such a way that
$$
\mu_j\in C^{10}(\overline{B(0,r)}), \quad \mu_1=\mu_2\quad \mbox{in}
\quad B(0,r) \setminus\Omega
$$
and the functions $\mu_j$ are constant in some neighborhood of $S(0,r).$
Therefore, without loss of generality, we can assume that $\Omega=B(0,r)$,
(\ref{boleslav}) holds true
and
$$
\mu_1(x)=\mu_2(x)=const \quad \forall x\,\,\mbox{from some
neighborhood of }\,\,\partial\Omega.
$$

Let $(\mbox{\bf u}_1,p_1)$ be the
complex geometric optics solution for the operator $G_{\mu_1}$
given by formulae (\ref{voz}). Then there exists a pair
$(\mbox{\bf u}_2,p_2)$ such that
$$
G_{\mu_2}(x,D)(\mbox{\bf u}_2,
p_2)=0\quad\mbox{in}\,\,\Omega, \quad
(\mbox{\bf u}_1 -\mbox{\bf u}_2)\vert_{\partial\Omega}
= \left(\frac{\partial \mbox{\bf u}_1}{\partial\nu}-\frac{\partial
\mbox{\bf u}_2}{\partial\nu}\right)\vert_{\partial\Omega}
= (p_1-p_2)\vert_{\partial\Omega}=0.
$$

We set $ \mbox{\bf u}=\mbox{\bf u}_1-\mbox{\bf u}_2$ and $p=p_1-p_2$.
Then
\begin{equation}\label{Bob}
L_{\mu_2}(x,D)(\mbox{\bf u}, p)+(\mbox{\bf u},\nabla)\mbox{\bf u}
= L_{\mu_2-\mu_1}(x,D)(\mbox{\bf
u}_1,p_1)\quad\mbox{in}\,\,\Omega, \,\, \mbox{\bf
u}\vert_{\partial\Omega}=\frac{\partial {\bf u}}{\partial\nu}
\vert_{\partial\Omega}=p\vert_{\partial\Omega}=0.
\end{equation}
Let $(\widetilde {\mbox{\bf v}},\widetilde p)$ be the complex geometric
optics solution to the Stokes equations
\begin{equation}
L_{\mu_2}(x,D)(\widetilde{\mbox{\bf v}},\widetilde p)=0
\quad\mbox{in}\,\,\Omega
\end{equation}
which is given by formula
$$
\widetilde {\mbox{\bf v}}={\bf v}-\nabla g,
 \quad V=\left (\begin{matrix} v_1\\v_2\\g\end{matrix}\right), \quad
V=e^{-\tau\overline \Phi}( V_0-V_1)
+ \sum_{j=2}^\infty(-1)^jV_j e^{-\tau \Phi},
$$ 
where the function $V$ is constructed in Proposition \ref{mursilka0}.
Taking the scalar product in $L^2(\Omega)$ of equation (\ref{Bob})
and the function  $\widetilde {\mbox{\bf v}}$, we obtain:
\begin{eqnarray}\label{magnit}
(L_{\mu_2}(x,D)(\mbox{\bf u}, p)+(\mbox{\bf u}_1,\nabla)
\mbox{\bf u}_1, \widetilde {\mbox{\bf v}})_{L^2(\Omega)}
= \int_{\partial\Omega}((\nu,\widetilde{\mbox{\bf v}})
- 2\sum_{i,j=1}^2\mu_2\nu_i\widetilde v_i\epsilon_{ij}(\mbox{\bf u}))
d\sigma\\
+ \int_{\partial\Omega}2\sum_{i,j=1}^2\mu_2\nu_i u_i\epsilon_{ij}
(\widetilde{\mbox{\bf v}})d\sigma
+ \int_\Omega p\mbox{div}\, \widetilde {\mbox{\bf v}}
+ (\mbox{\bf u}_1,L_{\mu_2}(x,D)(\widetilde{\mbox{\bf v}},\widetilde p)
-\nabla\widetilde p)dx+((\mbox{\bf u}_1,\nabla)\mbox{\bf u}_1, \widetilde
{\mbox{\bf v}})_{L^2(\Omega)}\nonumber\\=((\mbox{\bf u}_1,\nabla)
\mbox{\bf u}_1, \widetilde {\mbox{\bf v}})_{L^2(\Omega)}
= (L_{\mu_2-\mu_1}(x,D)(\mbox{\bf u}_1,p_1),\widetilde{\bf v})_{L^2(\Omega)}
= (\mbox{\bf u}_1,L_{\mu_2-\mu_1}(x,D)(\widetilde{\bf v},0))_{L^2(\Omega)}.
\nonumber
\end{eqnarray}

By (\ref{voz}) there exists a constant $C$ independent of $\tau$ such that
for all $\tau\ge \tau_0$
$$
\vert ((\mbox{\bf u}_1,\nabla)\mbox{\bf u}_1, \widetilde {\mbox{\bf v}})
_{L^2(\Omega)}\vert
\le \Vert e^{-\tau\varphi}(\mbox{\bf u}_1,\nabla)\mbox{\bf u}_1\Vert
_{L^2(\Omega)}\Vert e^{\tau \varphi}\widetilde {\mbox{\bf v}}\Vert
_{L^2(\Omega)}
$$
$$
\le \Vert e^{-\tau\varphi}\mbox{\bf u}_1\Vert_{L^4(\Omega)}\Vert
e^{- \tau-\tau\varphi}\nabla {\bf u}_1\Vert_{W_2^1(\Omega)}
\le C e^{-\tau/2}.
$$
Then this inequality and (\ref{magnit}) imply
$$
(\mbox{\bf u}_1,L_{\mu_2-\mu_1}(x,D)(\widetilde{\bf v},0))_{L^2(\Omega)}
= o\left(\frac 1\tau\right)\quad \mbox{as}\,\,\tau\rightarrow +\infty.
$$
Using (\ref{voz}) we rewrite the above equality as
$$
(\mbox{
u}_\tau,L_{\mu_2-\mu_1}(x,D)(\widetilde{\bf v},0))_{L^2(\Omega)}
= o\left(\frac 1\tau\right)\quad \mbox{as}\,\,\tau\rightarrow +\infty.
$$
Then thanks to (\ref{vo1}), we can apply Proposition \ref{vo2} to
transform the left-hand side of the above equality as
$$
(({\bf w},f),{\bf H}(x,\partial_z,\partial_{\overline z}) ({\bf v},g))
_{L^2(\Omega)}
= o\left(\frac 1\tau\right)\quad \mbox{as}\,\,\tau\rightarrow +\infty.
$$
We recall that the operator ${\bf H}(x,\partial_z,\partial_{\overline z})$
is given by (\ref{ii}) and (\ref{ii1}).
Since the coefficients $\mu_j$ are constants near $\partial\Omega$,
the operator ${\bf H}$ are compactly supported in $\Omega.$ Moreover since the  domain $\Omega$ is the ball, all the conditions of the Proposition
\ref{garmoshka} and Proposition \ref{zanoza} hold true.
This in turn implies that all the conditions on Proposition \ref{gavnuk}
hold true.
By Proposition \ref{gavnuk} the equality (\ref{kaput3}) holds true.

Denote
\begin{eqnarray}
\frak Q_{V_0,U_0,q_1,q_2}(\widetilde x)=
(P_{A_2}^*[\mathcal C_1^*U_0],\frac 12\partial_{z}[Q_2(2)V_0]
- \frac 14 B_2Q_2(2)V_0)(\widetilde x)                     \nonumber\\
+ (T^*_{B_1}[\mathcal C_2V_0](\widetilde x),\frac 12\partial_{\overline z}
[Q_1(1)U_0]-\frac 14A_1Q_1(1)U_0)(\widetilde x)
\nonumber\\
+2(\partial_zP_{A_2}^*[\mathcal C_1^*U_0]
,\frac 12Q_2(2)V_0)(\widetilde x)
+ 2(\partial_{\overline z}T^*_{B_1}[\mathcal C_2V_0],\frac 12 Q_1(1)U_0)
(\widetilde x)                              \nonumber\\
-(P_{A_2}^*(2\partial_z(\mathcal C_1^*U_0)+\partial_{\overline z}
(\mathcal C_0^*U_0)+{\bf B}_1^*U_0+Q_1(2)^*P_{A_2}^*(\mathcal C_1^*U_0)
+ \mathcal C_1^*q_1),\frac 12Q_2(2)V_0)(\widetilde x)\nonumber\\
- (T_{B_1}^*(2\mathcal C_2\partial_{\overline z}V_0
+ \mathcal C_0\partial_z V_0
+ {\bf B}_2V_0 +2Q_2(1)^*T_{B_1}^* (\mathcal C_2V_0)+\mathcal C_2q_2),
\frac 12 Q_1(1)U_0)(\widetilde x)\nonumber .
\end{eqnarray}
Then the equality (\ref{kaput3}) can be written in the form
\begin{equation}\label{pobeda}
\frak H_{V_0,U_0}(\widetilde x)
+ \frak Q_{V_0,U_0,q_1,q_2}(\widetilde x)=0.
\end{equation}
For the matrix differential operator ${\bf H}(x,\partial_z,\partial
_{\overline z})$, we denote the $(i,j)$-entry by ${\bf H}_{ij}
(x,\partial_z,\partial_{\overline z})$. The  operators $H_{ij}$ with either
$j\ne 3$ or $i\ne 3$ are the second order differential operators with respect
to $\mu.$

Let $\mathcal{C}_{j,k\ell}$ denote the $(k,\ell)$-entry of the matrix
$\mathcal C_j$.
From (\ref{ii}) and (\ref{ii1}) we compute
\begin{equation}
{\bf H}_{33}(x,\partial_z,\partial_{\overline z})g
= -2\sum_{i,j=1}^2\partial^2_{x_ix_j}\mu\partial^2_{x_ix_j}g
- 2\sum_{k=1}^2(\nabla \partial^2_{x_k x_k}\mu,\nabla g)
+ \mbox{div}\,\left(\frac{\mu}{\mu_2}\mu_2''\nabla g\right).
\end{equation}
Here we recall that $\mu_2''$ is the Hessian matrix.
Then by (\ref{poker}) we have
$$
\mathcal C_{1,33}=-2(\partial^2_{x_1x_1}\mu+2i\partial^2_{x_1x_2}\mu
-\partial^2_{x_2x_2}\mu),
\mathcal C_{2,33}=-2(\partial^2_{x_1x_1}\mu-2i\partial^2_{x_1x_2}\mu
-\partial^2_{x_2x_2}\mu)=-8\partial_{zz}^2\mu,
$$
$$
\mathcal C_{0,33}=-4(\partial^2_{x_1x_1}\mu+\partial^2_{x_2x_2}\mu),
{\bf B}_{1,33}
= -2(\partial^3_{x_1x_1x_1}\mu+i\partial^3_{x_1x_1x_2}\mu
+\partial^3_{x_1x_2x_2}\mu+i\partial^3_{x_2x_2x_2}\mu),
$$
$$
{\bf B}_{2,33}=-2(\partial^3_{x_1x_1x_1}\mu
-i\partial^3_{x_1x_1x_2}\mu+\partial^3_{x_1x_2x_2}\mu
-i\partial^3_{x_2x_2x_2}\mu)=-4\partial_z\Delta\mu.
$$

Then we have
\begin{equation} \partial^2_{\overline z\overline z} \mathcal C^*_{2,33}
-\partial_{\overline z}{\bf B}^*_{2,33}=8\Delta^2 \mu\quad\mbox{in}\,\,\Omega.
\end{equation}
For the fixed functions $U,V$ the function $\frak H_{U,V}$ can be considered
as the fourth-order operator applied to the function $\mu$.
The principal part of this operator is $8V_{0,3}U_{0,3}\Delta^2 \mu.$
By Proposition \ref{nikita}, for each point $\widetilde x$ we can choose
functions $U_{0,\widetilde x}=(U_{0,1,\widetilde x}, U_{0,2,\widetilde x},
U_{0,3,\widetilde x})$ and $V_{0,\widetilde x}=(V_{0,1,\widetilde x},
V_{0,2,\widetilde x},V_{0,3,\widetilde x})$ satisfying (\ref{-5})
and (\ref{-55}) respectively and
$V_{0,3,\widetilde x}(\widetilde x)=U_{0,3,\widetilde x}(\widetilde x)=1.$
Therefore for any $\widetilde x$ there exists $\delta (\widetilde x)>0$ such
that on the ball $B(\widetilde x,\delta(\widetilde x))$ we have
\begin{equation}\label{pobeda1}
V_{0,3,\widetilde x}(x)U_{0,3,\widetilde x}(x)\ge \frac 12\quad \forall
x\in B(\widetilde x,\delta(\widetilde x)).
\end{equation}
Since $\overline \Omega$ is covered by $\cup_{\widetilde x\in \overline\Omega}
B(\widetilde x,\delta(\widetilde x))$, from such a covering one can choose
a finite subcovering $B(\widetilde x_j,\delta(\widetilde x_j))$,
$j\in \{1,\dots, J\}.$

From (\ref{pobeda}) and (\ref{pobeda1}), there exist functions
$c_\beta ,p_k\in L^\infty(\Omega)$ such that
\begin{equation}\label{lazy5}
8\Delta^2 \mu+\sum_{\vert \beta\vert\le 3}c_\beta\partial^\beta_x \mu
+ \sum_{k=1}^J p_k \frak Q_{V_{0,\widetilde x_k},U_{0,\widetilde x_k},q_1,q_2}
=0\quad \mbox{in}\,\,\Omega .
\end{equation}
From (\ref{kaput3}) and Propositions \ref{inga} and \ref{zanoza!},
there exists a constant $C $ such that
\begin{equation}\label{lazy2}
\Vert s\phi_s^\frac 12\frak Q_{V_0,U_0,q_1,q_2}e^{s\phi_s}\Vert_{L^2(\Omega)}
\le C\Vert \sum_{\vert\beta\vert\le 3}\partial^\beta_x \mu e^{s\phi_s}\Vert
_{L^2(\Omega)}\quad\forall s\ge s_0.
\end{equation}
Indeed in order to prove (\ref{lazy2}), observe that  by Proposition
\ref{inga} and Proposition \ref{zanoza!}
\begin{equation}\label{lenta1}
\Vert s\phi_s^\frac 12 P_{A_2}^*[\mathcal C_1^*U_0]e^{s\phi_s}\Vert
_{L^2(\Omega)}
\le C\Vert \mathcal C_1^*U_0e^{s\phi_s}\Vert_{L^2(\Omega)}
\le C\Vert \mathcal C_1^*e^{s\phi_s}\Vert_{L^2(\Omega)}
\end{equation}
and
\begin{equation}\label{lenta2}
\Vert s\phi_s^\frac 12 T_{B_1}^*[\mathcal C_2^*V_0]e^{s\phi_s}\Vert
_{L^2(\Omega)}
\le C\Vert \mathcal C_2^*V_0e^{s\phi_s}\Vert_{L^2(\Omega)}
\le C\Vert \mathcal C_2^*e^{s\phi_s}\Vert_{L^2(\Omega)}.
\end{equation}
Let the functions $q_j$ are given by (\ref{victoryz}).
By (\ref{lenta1}), (\ref{!zanoza1!}), Proposition \ref{inga} and Proposition \ref{zanoza!}, we have
\begin{eqnarray}\label{lenta3}
\Vert s\phi_s^\frac 12(P_{A_2}^*(2\partial_z(\mathcal C_1^*U_0)+\partial
_{\overline z}(\mathcal C_0^*U_0)+{\bf B}_1^*U_0+Q_1(2)^*P_{A_2}^*
(\mathcal C_1^*U_0)
+\mathcal C_1^*q_1
)e^{s\phi_s}\Vert_{L^2(\Omega)} \nonumber\\
\le \Vert s\phi_s^\frac 12(P_{A_2}^*(2\partial_z(\mathcal C_1^*U_0)+\partial
_{\overline z}(\mathcal C_0^*U_0)+{\bf B}_1^*U_0+Q_1(2)^*P_{A_2}^*
(\mathcal C_1^*U_0)
+\mathcal C_1^*q_1^0
)e^{s\phi_s}\Vert_{L^2(\Omega)}\nonumber\\+\Vert s\phi_s^\frac 12\sum_{j=1}^3r_{1,j} P_{A_2}^*(\mathcal C_1^*q_{1,j}
)e^{s\phi_s}\Vert_{L^2(\Omega)}\nonumber\\
\le C \Vert (2\partial_z(\mathcal C_1^*U_0)+\partial_{\overline z}
(\mathcal C_0^*U_0)+{\bf B}_1^*U_0+Q_1(2)^*P_{A_2}^*(\mathcal C_1^*U_0)
+\mathcal C_1^*q_1^0
)e^{s\phi_s}\Vert_{L^2(\Omega)}\nonumber\\+C\sum_{j=1}^3\Vert \mathcal C_1^*q_{1,j}
e^{s\phi_s}\Vert_{L^2(\Omega)}\nonumber\\
\le C \Vert (2\partial_z(\mathcal C_1^*U_0)+\partial_{\overline z}
(\mathcal C_0^*U_0)+{\bf B}_1^*U_0+\mathcal C_1^*)e^{s\phi_s}\Vert
_{L^2(\Omega)}
+C\Vert \mathcal C_1^*e^{s\phi_s}\Vert_{L^2(\Omega)}.
\end{eqnarray}
By (\ref{lenta3}), (\ref{!zanoza2!}), Propositions \ref{inga} and \ref{zanoza!}, we obtain
\begin{eqnarray}\label{lenta4}
\Vert  s\phi_s^\frac 12(T_{B_1}^*(2\mathcal C_2\partial_{\overline z}V_0
+\mathcal C_0\partial_z V_0+{\bf B}_2V_0 +2Q_2(1)^*T_{B_1}^*
(\mathcal C_2V_0)
+\mathcal C_2q_2
)e^{s\phi_s}\Vert_{L^2(\Omega)}\nonumber\\
\le\Vert  s\phi_s^\frac 12(T_{B_1}^*(2\mathcal C_2\partial_{\overline z}V_0
+\mathcal C_0\partial_z V_0+{\bf B}_2V_0 +2Q_2(1)^*T_{B_1}^*
(\mathcal C_2V_0)
+\mathcal C_2q_2^0
)e^{s\phi_s}\Vert_{L^2(\Omega)}\nonumber\\
+\Vert s\phi_s^\frac 12\sum_{j=1}^3r_{2,j} T_{B_1}^*(\mathcal C_2q_{2,j}
)e^{s\phi_s}\Vert_{L^2(\Omega)}\nonumber\\
 \le C\Vert
(2\mathcal C_2\partial_{\overline z}V_0+\mathcal C_0\partial_z V_0
+{\bf B}_2V_0 +\mathcal C_2q_2^0
)e^{s\phi_s}\Vert_{L^2(\Omega)}
+C\Vert \mathcal C_2e^{s\phi_s}\Vert_{L^2(\Omega)}
.
\end{eqnarray}

Finally, observing that
$$
\sum_{j=0}^2 (\vert \nabla {\mathcal C_j}(x)\vert+\vert  {\mathcal C_j}
(x)\vert) +\sum_{j=1}^2\vert {\bf B}_j(x)\vert\le\sum_{\vert\beta\vert\le 3}
\vert\partial^\beta_x \mu (x)\vert
\quad \forall x\in\Omega
$$
from (\ref{lenta1})-(\ref{lenta4}), we obtain (\ref{lazy2}).

Applying the Carleman estimate (\ref{lazy1}) to equation (\ref{lazy5})
and using (\ref{lazy2}), we obtain $\mu\equiv 0.$
Thus the proof of the theorem is complete. $\blacksquare$

\section{Construction of complex geometric optics solutions for the Lam\'e
system and proof of Theorem 1.2 }\label{sec5}

Denote
\begin{eqnarray*}
&&\mathcal Z_\alpha(x,D) {\bf v}
= \left(\sum_{j=1}^2\partial_{x_j}(\alpha(\partial_{x_j}v_1
+\partial_{x_1}v_j)),
\sum_{j=1}^2\partial_{x_j}(\alpha(\partial_{x_j}v_2+\partial_{x_2}v_j)\right)\\
&&= \alpha\Delta {\bf v} + \alpha\nabla\dd{\bf v}
+ ((\nabla\alpha,\nabla v_1), (\nabla\alpha,\nabla v_2))
+ ((\nabla\alpha,\partial_{x_1}{\bf v}), (\nabla\alpha,\partial_{x_2}{\bf v})).
\end{eqnarray*}
Then we note
$$
\mathcal L_{\mu,\lambda}(x,D){\bf w}
= \mathcal{Z}_{\mu}(x,D){\bf w} + \nabla(\lambda\dd{\bf w}).
$$
\begin{proposition}\label{lame}
Let ${\bf w}=(w_1,w_2) $ and $f$ satisfy the elliptic system
\begin{eqnarray}\label{-111}
(\lambda+2\mu)\Delta f+(\lambda+\mu)\mbox{div}\thinspace{\bf w}
+2(\nabla\mu,\nabla f)=0\quad\mbox{in}\,\,\Omega
\end{eqnarray}
and
\begin{equation}\label{-112}
\mu\Delta {\bf w}+{\bf M}(x,D)({\bf w},f)=0\quad\mbox{in}\,\,\Omega,
\end{equation}
where
\begin{eqnarray}\label{yabloko}
{\bf M}_2(x,D)({\bf w},f)= -(\nabla \mu)\dd{\bf w}
+ ((\nabla\mu,\nabla w_1),(\nabla\mu,\nabla w_2))
+ ((\nabla \mu, \partial_{x_1}{\bf w}), (\nabla\mu,\partial_{x_2}{\bf w}))
                                         \nonumber\\
+ 2\frac{\lambda+\mu}{\lambda+2\mu}(\nabla \mu)\dd{\bf w}
- 2((\partial_{x_1}\nabla\mu,\nabla f),(\partial_{x_2}\nabla\mu,\nabla f))
+ 4(\nabla\mu)\frac{(\nabla\mu,\nabla f)}{\lambda+2\mu}.
\end{eqnarray}
Then the function
 ${\bf u}={\bf
w}+\nabla f$ solves the Lam\'e system (\ref{lammas}).
\end{proposition}

{\bf Proof.}
Short computations provide
\begin{equation}\label{Xz0}
\nabla (\lambda\mbox{div}\,\nabla f)=
\nabla (\lambda\Delta f)
\end{equation}
and
\begin{eqnarray}\label{Xz1}
\mathcal Z_\mu(x,D)\nabla f=
\left(\sum_{j=1}^2\partial_{x_j}(2\mu\partial_{x_1x_j} f),
\sum_{j=1}^2\partial_{x_j}(2\mu\partial_{x_2x_j} f)\right)\\
= 2\mu\nabla \Delta f+2((\nabla\mu,\nabla \partial_{x_1}f),
(\nabla\mu,\nabla \partial_{x_2}f))                    \nonumber\\
= 2\nabla (\mu\Delta f) - 2(\nabla \mu) \Delta f+2\nabla (\nabla\mu,\nabla f)
-2(((\partial_{x_1}\nabla\mu,\nabla f),(\partial_{x_2}
\nabla\mu,\nabla f)).\nonumber
\end{eqnarray}
Next we compute
\begin{eqnarray}\label{Xz2}
{\mathcal L}_{\mu,\lambda}(x,D){\bf w}
= \left(\sum_{j=1}^2\partial_{x_j}(\mu(\partial_{x_j}w_1+\partial_{x_1}w_j)),
\sum_{j=1}^2\partial_{x_j}(\mu(\partial_{x_j}w_2+\partial_{x_2}w_j))\right)
                    \nonumber\\
+ \nabla(\lambda\mbox{div}\,{\bf w})           \nonumber\\
= \mu\Delta {\bf w}+\mu\nabla \mbox{div}\,{\bf w}+((\nabla\mu,\nabla w_1),
(\nabla\mu,\nabla w_2))+((\nabla \mu, \partial_{x_1}{\bf w}),
(\nabla\mu,\partial_{x_2}{\bf w}))\nonumber\\ +\nabla(\lambda
\mbox{div}\,{\bf w})               \nonumber\\
= \mu\Delta {\bf w}+\nabla (\mu\mbox{div}\,{\bf w})
- (\nabla \mu)\mbox{div}\,{\bf w}
+ ((\nabla\mu,\nabla w_1),(\nabla\mu,\nabla w_2))
                                              \nonumber\\
+ ((\nabla \mu, \partial_{x_1}{\bf w}), (\nabla\mu,\partial_{x_2}{\bf w}))
+\nabla(\lambda\mbox{div}\,{\bf w})
\end{eqnarray}
$$
= {\mathcal Z}_{\mu}(x,D){\bf w} + \nabla(\lambda\dd{\bf w}).
$$
Combining the formulae (\ref{Xz0}), (\ref{Xz1}) and (\ref{Xz2}) we obtain
\begin{eqnarray}\label{lam}
{\mathcal L}_{\mu,\lambda}(x,D)({\bf w}+\nabla f)=\mu\Delta {\bf w}
- (\nabla \mu)\mbox{div}\,{\bf w}+((\nabla\mu,\nabla w_1),
(\nabla\mu,\nabla w_2))                       \nonumber\\
+\nabla((\lambda+\mu)\mbox{div} \mbox\,{\bf w})
+((\nabla \mu, \partial_{x_1}{\bf w}), (\nabla\mu,\partial_{x_2}{\bf w}))
                           \nonumber\\
+ \nabla ((\lambda+2\mu)\Delta f) - 2(\nabla \mu)\Delta f
+ 2\nabla(\nabla\mu,\nabla f)\nonumber\\
- 2((\nabla\partial_{x_1}\mu,\nabla f),(\nabla\partial_{x_2}\mu,\nabla f)).
\end{eqnarray}
Using the equation (\ref{-111}) we have
\begin{eqnarray}\label{lam}
{\mathcal L}_{\mu,\lambda}(x,D)({\bf w}+\nabla f)
= \mu\Delta {\bf w} - (\nabla \mu)\mbox{div}\,{\bf w}
+((\nabla\mu,\nabla w_1),(\nabla\mu,\nabla w_2))\nonumber\\
+((\nabla \mu, \partial_{x_1}{\bf w}), (\nabla\mu,\partial_{x_2}{\bf w}))
                        \nonumber\\
+ 2\frac{\lambda+\mu}{\lambda+2\mu}(\nabla \mu)\mbox{div}\,{\bf w}
- 2(((\partial_{x_1}\nabla\mu,\nabla f),(\partial_{x_2}\nabla\mu,\nabla f))
+ 4\nabla\mu \frac{(\nabla\mu,\nabla f)}{\lambda+2\mu}.
\end{eqnarray}
By (\ref{-112}) and (\ref{yabloko}) the right-hand side of (\ref{lam}) is zero.
The proof of the proposition is complete.
$\blacksquare$

Let
$$
{\bf P}(x,D)({\bf v},g)=-2((\partial_{x_1}\nabla\alpha,\nabla g),
(\partial_{x_2}\nabla\alpha,\nabla g))
$$
$$
+\left\{\frac{2\lambda_2+2\mu_2}{\lambda_2+2\mu_2} \mbox{div}\, {\bf v}
+ \frac{4(\nabla\mu_2,\nabla g)}{\lambda_2+2\mu_2}\right \}\nabla \alpha
+ \mathcal K_\alpha(x,D){\bf v}
- \frac{\alpha}{\mu_2}{\bf M}_2(x,D)({\bf v},g),
$$
where
$$
\mathcal K_\alpha(x,D){\bf v}
= ((\nabla\alpha,\nabla v_1),(\nabla\alpha,\nabla v_2))
+ ((\nabla\alpha,\partial_{x_1}{\bf v}),(\nabla\alpha,\partial_{x_2}{\bf v}))
- (\nabla \alpha)\mbox {div}\, {\bf v}
$$
and
$$
P_{\beta}(x,D)({\bf v},g)
= -\left ( \begin{matrix}-\nabla\left\{ \frac{\beta\mu_1}{\lambda_1
+ 2\mu_1}\left (\frac{\mu_2}{(\lambda_2+2\mu_2)}\mbox{div}\,{\bf v}
-\frac{2(\nabla\mu_2,\nabla g)}{\lambda_2+2\mu_2}\right )\right\}\\
\mbox{div}\,\left \{(\nabla\mu_1)\frac{2\beta}{\lambda_1+2\mu_1}
\left(\frac{\mu_2}{(\lambda_2+2\mu_2)}\mbox{div}\,{\bf v}
-\frac{2(\nabla\mu_2,\nabla g)}{\lambda_2+2\mu_2}\right)\right\}
\end{matrix}
\right).
$$
We have
\begin{proposition}\label{ignalina}
Let $\alpha,\beta\in C^4_0(\Omega)$ and
let $({\bf w} ,f), ({\bf v},g)$ be some regular  solutions to the
system (\ref{-111})-(\ref{-112}) with the Lam\'e coefficients
$(\mu_1,\lambda_1)$ and $(\mu_2,\lambda_2)$ respectively.
Then
\begin{eqnarray}\label{pop}
({\bf w}+\nabla f,{\mathcal  L}_{\alpha,\beta}(x,D)({\bf v}+\nabla g))
_{L^2(\Omega)}
= (({\bf w},f),{\bf H}(x,\partial_z,\partial_{\overline z}) ({\bf v},g))
_{L^2(\Omega)},
\end{eqnarray}
where
\begin{equation}\label{neapol}
{\bf H}(x,\partial_z,\partial_{\overline z}) ({\bf v},g)
= \left( \begin{matrix}
-\nabla \left \{\alpha \left (\frac{\lambda_2}{\lambda_2+2\mu_2}
\mbox{div}\,{\bf v}+\frac{4(\nabla\mu_2,\nabla g)}{\lambda_2+2\mu_2}\right)
\right\}\\ 0\end{matrix}\right)+\left(\begin{matrix}{\bf P}(x,D)({\bf v},g)\\-\mbox{div}\,{\bf P}(x,D)({\bf v},g)\end{matrix}\right)
\end{equation}
$$
+ P_{\beta}(x,D)({\bf v},g)
$$
$$
+ \left( \begin{matrix} \nabla\left\{
\alpha\frac{\lambda_1+\mu_1}{\lambda_1+2\mu_1} \left \{\frac{\lambda_2}
{\lambda_2+2\mu_2} \mbox{div}\,{\bf v}+\frac{4(\nabla\mu_2,\nabla g)}
{\lambda_2+2\mu_2}\right\}\right\}\\
-\mbox{div}\,\left\{(\nabla\mu_1)\frac{2\alpha}{\lambda_1+2\mu_1} \left
\{\frac{\lambda_2}{\lambda_2+2\mu_2} \mbox{div}\,{\bf v}
+\frac{4(\nabla\mu_2,\nabla g)}{\lambda_2+2\mu_2}\right\}\right \}
\end{matrix}\right)
+ \left(
\begin{matrix}\nabla \left(\frac{\mu_1}{\lambda_1+2\mu_1}
(\nabla\alpha,\nabla g) \right)\\
-\mbox{div}\left( \frac{4\nabla\mu_1}{\lambda_1+2\mu_1}
(\nabla\alpha,\nabla g) \right)
\end{matrix}\right).
$$
Here the operator ${\bf M}_2(x,D)$ is given by (\ref{yabloko}) with the
functions $\lambda_2,\mu_2$ instead of $\lambda$ and $\mu .$
\end{proposition}
{\bf Proof.} Integrating by parts and using equation (\ref{-111}) we have
\begin{eqnarray}\label{neapol1}
\int_\Omega(({\bf w}+\nabla f)),\nabla (\beta\mbox{div} ({\bf v}
+\nabla g))dx
= -\int_\Omega\beta\mbox{div}({\bf w}+\nabla f) \mbox{div}({\bf v}
+\nabla g)dx\nonumber\\
= -\int_\Omega\beta(\mbox{div}\,{\bf w}+\Delta f) (
\mbox{div}\,{\bf v}+\Delta g)dx\nonumber\\
= -\int_\Omega\beta\left(\frac{\mu_1}{(\lambda_1+2\mu_1)}
\mbox{div}\,{\bf w}-\frac{2(\nabla\mu_1,\nabla f)}{\lambda_1+2\mu_1}\right)
\left(\frac{\mu_2}{(\lambda_2+2\mu_2)}\mbox{div}\,{\bf v}-\frac{2(\nabla\mu_2,
\nabla g)}{\lambda_2+2\mu_2}\right )dx               \nonumber\\
= (({\bf w},f), P_{\beta}(x,D)({\bf v},g))_{L^2(\Omega)}.
\end{eqnarray}

Observe that
\begin{eqnarray}
\mathcal Z_{\alpha}(x,D){\bf v}=\alpha \Delta{\bf v}
+ \nabla (\alpha \mbox{div}\, {\bf v})+\mathcal K_\alpha(x,D){\bf v}\nonumber\\
= -\frac{\alpha}{\mu_2}{\bf M}_2({\bf v},g)+\nabla (\alpha \mbox{div}\,{\bf v})
+ \mathcal K_\alpha(x,D){\bf v}\nonumber.
\end{eqnarray}

Hence, using (\ref{neapol1}) and (\ref{Xz1}), we have
\begin{eqnarray}\label{neapol2}
({\bf w}+\nabla f,{\mathcal  Z}_{\alpha}(x,D)({\bf v}+\nabla g))_{L^2(\Omega)}
                                             \nonumber\\
= ({\bf w}+\nabla f,
2\nabla(\alpha\Delta g) - 2(\nabla\alpha)\Delta g
+ 2((\nabla\alpha,\nabla \partial_{x_1}g),
(\nabla\alpha,\nabla \partial_{x_2}g))
+ \nabla (\alpha \mbox{div}\, {\bf v})                     \nonumber\\
+ \mathcal K_\alpha(x,D){\bf v}-\frac{\alpha}{\mu_2}{\bf M}_2(x,D)
({\bf v},g))_{L^2(\Omega)}                               \nonumber\\
= ({\bf w}+\nabla f,-\nabla\left \{\alpha\left (\frac{\lambda_2}
{\lambda_2+2\mu_2} \mbox{div}\,{\bf v}+\frac{4(\nabla\mu_2,\nabla g)}
{\lambda_2+2\mu_2}\right)\right\} +2((\nabla\alpha,\nabla \partial_{x_1}g),
(\nabla\alpha,\nabla \partial_{x_2}g))         \nonumber\\
+ \left\{\frac{2\lambda_2+2\mu_2}{\lambda_2+2\mu_2} \mbox{div}\, {\bf v}
+ \frac{4(\nabla\mu_2,\nabla g)}{\lambda_2+2\mu_2}\right \}\nabla \alpha
+ \mathcal K_\alpha(x,D){\bf v}-\frac{\alpha}{\mu_2}{\bf M}_2(x,D)({\bf v},g))
_{L^2(\Omega)}                   \nonumber\\
\end{eqnarray}
Here at the second equality of (\ref{neapol2}), we use:
\begin{eqnarray*}
&& 2\nabla(\alpha\Delta g) - 2(\nabla\alpha)\Delta g
+ \nabla(\alpha \dd{\bf v})
= \nabla(\alpha(2\Delta g + \dd{\bf v}))
- 2(\nabla\alpha)\Delta g\\
&&= -2\nabla\left(\alpha
\left(\frac{\lambda_2+\mu_2}{\lambda_2+2\mu_2}\dd{\bf v}
+ \frac{2(\nabla\mu_2, \nabla g)}{\lambda_2+2\mu_2}
- \frac{1}{2}\dd{\bf v}\right)\right)\\
&&+ 2\nabla\alpha
\left(\frac{\lambda_2+\mu_2}{\lambda_2+2\mu_2}\dd{\bf v}
+ \frac{2(\nabla\mu_2, \nabla g)}{\lambda_2+2\mu_2}\right)\\
&&= -2\nabla\left(\alpha
\left(\frac{\lambda_2}{2\lambda_2+4\mu_2}\dd{\bf v}
+ \frac{2(\nabla\mu_2, \nabla g)}{\lambda_2+2\mu_2}\right)\right)\\
&&+ 2\nabla\alpha
\left(\frac{\lambda_2+\mu_2}{\lambda_2+2\mu_2}\dd{\bf v}
+ \frac{2(\nabla\mu_2, \nabla g)}{\lambda_2+2\mu_2}\right)
\end{eqnarray*}
by
$$
\Delta g = - \frac{\lambda_2+\mu_2}{\lambda_2+2\mu_2}\dd {\bf v}
- \frac{2(\nabla\mu_2, \nabla g)}{\lambda_2+2\mu_2}.
$$

Integrating by parts and using the fact what the function $\alpha$ has a
compact support in $\Omega$ we obtain
\begin{eqnarray}
({\bf w}+\nabla f,2((\nabla\alpha,\nabla \partial_{x_1}g),
(\nabla\alpha,\nabla \partial_{x_2}g)))_{L^2(\Omega)}
= -2(\mbox{div}{\bf w}+ \Delta f,
(\nabla\alpha,\nabla g))_{L^2(\Omega)}\nonumber\\
- 2({\bf w}+\nabla f,((\partial_{x_1}\nabla\alpha,\nabla g),
(\partial_{x_2}\nabla\alpha,\nabla g)))_{L^2(\Omega)}
\end{eqnarray}
By (\ref{-111}) we have
$$
\mbox{div}{\bf w}+\Delta f=\frac{\mu_1}{\lambda_1+2\mu_1}\mbox{div}
\,{\bf w}-\frac{2(\nabla\mu_1,\nabla f)}{\lambda_1+2\mu_1}.$$
Therefore
\begin{eqnarray}\label{mark}
({\bf w}+\nabla f,2((\nabla\alpha,\nabla \partial_{x_1}g),(\nabla\alpha,\nabla
\partial_{x_2}g)))_{L^2(\Omega)}
= -2\left(\frac{\mu_1}{\lambda_1+2\mu_1}\mbox{div}\,{\bf w}
-\frac{2(\nabla\mu_1,\nabla f)}{\lambda_1+2\mu_1},(\nabla\alpha,\nabla g)
\right)_{L^2(\Omega)}                        \nonumber\\
- 2({\bf w}+\nabla f,((\partial_{x_1}\nabla\alpha,\nabla g),
(\partial_{x_2}\nabla\alpha,\nabla g)))_{L^2(\Omega)}
\end{eqnarray}
Using this, we obtain
\begin{eqnarray}
({\bf w}+\nabla f,{\mathcal  Z}_{\alpha}(x,D)({\bf v}+\nabla g))_{L^2(\Omega)}
                                             \nonumber\\
= \left({\bf w}+\nabla f, -\nabla\left \{\alpha \left (\frac{\lambda_2}
{\lambda_2+2\mu_2}
\mbox{div}\,{\bf v}+\frac{4(\nabla\mu_2,\nabla g)}{\lambda_2+2\mu_2}\right)
\right\}\right)_{L^2(\Omega)}
+ (({\bf w}+\nabla f),{\bf P}(x,D)({\bf v},g))_{L^2(\Omega)}         \nonumber\\- 2\left(\frac{\mu_1}{\lambda_1+2\mu_1}\mbox{div}\,{\bf w}
-\frac{2(\nabla\mu_1,\nabla f)}{\lambda_1+2\mu_1},(\nabla\alpha,\nabla g)
\right)_{L^2(\Omega)}\nonumber\\
= \left({\bf w},
-\nabla \left \{\alpha \left (\frac{\lambda_2}{\lambda_2+2\mu_2}
\mbox{div}\,{\bf v}+\frac{4(\nabla\mu_2,\nabla g)}{\lambda_2+2\mu_2}\right )
\right\}\right)_{L^2(\Omega)}+(({\bf w}+\nabla f),{\bf P}(x,D)({\bf v},g))
_{L^2(\Omega)}                                             \nonumber\\
- \left(\frac{\lambda_1+\mu_1}{\lambda_1+2\mu_1}\mbox{div}\, {\bf w}
-\frac{2(\nabla\mu_1,\nabla f)}{\lambda_1+2\mu_1},
\alpha \left \{\frac{\lambda_2}{\lambda_2+2\mu_2}
\mbox{div}\,{\bf v}+\frac{4(\nabla\mu_2,\nabla g)}{\lambda_2+2\mu_2}\right\}
\right)_{L^2(\Omega)}                   \nonumber\\
-2\left(\frac{\mu_1}{\lambda_1+2\mu_1}\mbox{div}\,{\bf w}
-\frac{2(\nabla\mu_1,\nabla f)}{\lambda_1+2\mu_1},(\nabla\alpha,\nabla g)
\right)_{L^2(\Omega)}\nonumber\\
= \left({\bf w},
-\nabla \left \{\alpha \left (\frac{\lambda_2}{\lambda_2+2\mu_2}
\mbox{div}\,{\bf v}+\frac{4(\nabla\mu_2,\nabla g)}{\lambda_2+2\mu_2}\right )
\right\}\right)_{L^2(\Omega)}
+ (({\bf w}+\nabla f),{\bf P}(x,D)({\bf v},g))_{L^2(\Omega)}
                                                   \nonumber\\
+ \left({\bf w},\nabla\left\{
\alpha\frac{\lambda_1+\mu_1}{\lambda_1+2\mu_1} \left \{\frac{\lambda_2}
{\lambda_2+2\mu_2} \mbox{div}\,{\bf v}+\frac{4(\nabla\mu_2,\nabla g)}
{\lambda_2+2\mu_2}\right\}\right\}\right)_{L^2(\Omega)}          \nonumber\\
- \left(f,
\mbox{div}\,\left\{(\nabla\mu_1)\frac{2\alpha}{\lambda_1+2\mu_1}
\left \{\frac{\lambda_2}{\lambda_2+2\mu_2} \mbox{div}\,{\bf v}+\frac{4(\nabla\mu_2,\nabla g)}{\lambda_2+2\mu_2}\right\}\right \}\right)_{L^2(\Omega)}
                                   \nonumber\\
+ \left({\bf w},\nabla \left( \frac{\mu_1}{\lambda_1+2\mu_1}
(\nabla\alpha,\nabla g) \right)\right)_{L^2(\Omega)}
- \left(f,\mbox{div}
\left( \frac{4\nabla\mu_1}{\lambda_1+2\mu_1}(\nabla\alpha,\nabla g) \right)
\right)_{L^2(\Omega)}.
\end{eqnarray}
At the final equality, we used (\ref{neapol1}) and
\begin{eqnarray*}
&&({\bf w}+\nabla f, {\bf P}(x,D)({\bf v},g)) \\
&&= ({\bf w}, {\bf P}(x,D)({\bf v},g))
- (f, \dd {\bf P}(x,D)({\bf v},g)).
\end{eqnarray*}
From (\ref{neapol1}) and (\ref{neapol2}), we obtain (\ref{pop}).
The proof of the proposition is complete.
$\blacksquare$

{\bf Proof of Theorem 1.2.}
Suppose that for positive smooth functions $\mu_j,\lambda_j$,
the Dirichlet-to-Neumann maps (\ref{z1XV}) are the same.
In \cite{ANS} it is proved
\begin{equation}\label{boleslavX}
\frac{\partial^\ell \mu_1}{\partial \nu^\ell}-\frac{\partial^\ell \mu_2}
{\partial \nu^\ell}=\frac{\partial^\ell \lambda_1}{\partial \nu^\ell}
-\frac{\partial^\ell \lambda_2}{\partial \nu^\ell}=0\quad\mbox{on}
\quad\partial\Omega, \quad\forall \ell\in \{0,\dots, 10\}.
\end{equation}
Since the domain $\Omega$ is assumed to be bounded, there
exists a ball $B(0,r)$ such that $\overline \Omega\subset B(0,r).$
Thanks to (\ref{boleslavX}) we can extend the coefficients $\mu_j,\lambda_j$
into $B(0,r)\setminus\Omega$ such that
$$
\mu_j,\lambda_j\in C^{10}(B(0,r)), \quad \mu_1-\mu_2=\lambda_1-\lambda_2
=0\quad \mbox{in} \quad B(0,r)\setminus\Omega
$$
and the functions $\mu_j,\lambda_j$ are constant in some neighborhood of
$S(0,r).$   Therefore we can assume that $\Omega=B(0,r)$, (\ref{boleslavX})
holds true and
\begin{equation}\label{oleg}
\mu_1(x)-\mu_2(x)=\lambda_1(x)-\lambda_2(x)=const \quad
\forall x\,\,\mbox{from some neighborhood of }\,\,\partial\Omega.
\end{equation}

Let $\mbox{\bf u}_1$ be the
complex geometric optics solution for the operator $\mathcal
L_{\mu_1,\lambda_1}(x,D)$
given by
\begin{equation}\label{!vo1}
\mbox{\bf u}_1=\left(\begin{matrix}w_1\\w_2\end{matrix}\right )-\nabla f,
\quad U=\left (\begin{matrix} w_1\\w_2\\f\end{matrix}\right),\quad
U=e^{\tau\Phi}( U_0-U_1)+\sum_{j=2}^\infty(-1)^jU_j e^{\tau\overline \Phi},
\end{equation}where
the function $U$ is constructed in  Proposition \ref{mursilka0}.
Then there exists a function
$\mbox{\bf u}_2$ such that
$$
\mathcal L_{\mu_2,\lambda_2}(x,D)\mbox{\bf u}_2=0\quad\mbox{in}
\,\,\Omega, \quad (\mbox{\bf u}_1-\mbox{\bf
u}_2)\vert_{\partial\Omega}=(\frac{\partial \mbox{\bf u}_1}
{\partial\nu}-\frac{\partial \mbox{\bf
u}_2}{\partial\nu})\vert_{\partial\Omega}=0.
$$

We set $ \mbox{\bf u}=\mbox{\bf u}_1-\mbox{\bf u}_2$, $\mu=\mu_2-\mu_1$ and
$\lambda=\lambda_2-\lambda_1$.
Then
\begin{equation}\label{BobA}
\mathcal L_{\mu_2,\lambda_2}(x,D)\mbox{\bf u}=\mathcal L_{\mu,\lambda}(x,D)
\mbox{\bf u}_1\quad\mbox{in}\,\,\Omega, \,\, \mbox{\bf
u}\vert_{\partial\Omega}=\frac{\partial {\bf u}}{\partial\nu}\vert
_{\partial\Omega}=0.
\end{equation}
Let $\widetilde {\mbox{\bf v}}$ be the complex geometric optics solution to
the system
\begin{equation}
\mathcal L_{\mu_2,\lambda_2}(x,D)\widetilde{\mbox{\bf v}}=0\quad\mbox{in}
\,\,\Omega
\end{equation}
which is given by formula
$$
\widetilde {\mbox{\bf v}}=\left (\begin{matrix}  v_1\\
v_2\end{matrix}\right)-\nabla g,
 \quad V=\left (\begin{matrix} v_1\\v_2\\g\end{matrix}\right),\quad
V=e^{-\tau\overline \Phi}( V_0-V_1)
+ \sum_{j=2}^\infty(-1)^jV_j e^{-\tau \Phi},
$$ 
where the function $V$ is constructed in  Proposition \ref{mursilka0}.
Taking the scalar product in $L^2(\Omega)$ of equation (\ref{BobA})
and the function  $\widetilde {\mbox{\bf v}}$ we obtain:
\begin{eqnarray}\label{magnitA}
(\mathcal L_{\mu_2,\lambda_2}(x,D)\mbox{\bf u}, \widetilde {\mbox{\bf v}})
_{L^2(\Omega)}=\int_{\partial\Omega}((\lambda_2+\mu_2)\mbox{div}\,{\bf u}
(\nu,\widetilde{\mbox{\bf v}})+2\sum_{i,j=1}^2
\mu_2\nu_i\widetilde v_i\epsilon_{ij}(\mbox{\bf u}))d\sigma\\
-\int_{\partial\Omega}2\sum_{i,j=1}^2\mu_2\nu_i u_i\epsilon_{ij}
(\widetilde{\mbox{\bf v}})d\sigma
-\int_{\partial\Omega} (\lambda_2+\mu_2)(\nu,{\bf u})\mbox{div}\,
\widetilde {\mbox{\bf v}}d\sigma+\int_\Omega(\mbox{\bf u}_1,
\mathcal L_{\mu_2,\lambda_2}(x,D)\widetilde{\mbox{\bf v}})dx\nonumber\\
=(\mathcal L_{\mu,\lambda}(x,D)\mbox{\bf u}_1,\widetilde{\bf v})_{L^2(\Omega)}
= (\mbox{\bf u}_1,\mathcal L_{\mu,\lambda}(x,D)\widetilde{\bf v})_{L^2(\Omega)}
=0.\nonumber
\end{eqnarray}

Then thanks to  (\ref{oleg}), we can apply Proposition \ref{ignalina} to
transform the left-hand side of the above equality as
$$
(({\bf w},f),{\bf H}(x,\partial_z,\partial_{\overline z}) ({\bf v},g))
_{L^2(\Omega)}= o\left(\frac 1\tau\right)\quad \mbox{as}
\,\,\tau\rightarrow +\infty.
$$
We recall that the operator ${\bf H}(x,\partial_z,\partial_{\overline z})$
is given by (\ref{neapol}).
Since the coefficients $\mu_j,\lambda_j$ are constants near $\partial\Omega$,
the operator ${\bf H}$ is compactly supported in $\Omega.$
By Proposition \ref{gavnuk}, the equality (\ref{kaput3}) holds true.

The equality (\ref{kaput3}) can be written in the form
\begin{equation}\label{pobedaA}
\frak H_{V_0,U_0}(x)+\frak Q_{V_0,U_0,q_1,q_2}(x)=0\quad\mbox{in}\,\,\Omega.
\end{equation}Computing the coefficients of the operator ${\bf H}(x,\partial_z,\partial_{\overline z})$ we have
\begin{equation}\label{!!PP}
\mathcal C_2=\left (\begin{matrix}\frac {\lambda\mu_1\mu_2}
{(\lambda_1+2\mu_1) (\lambda_2+2\mu_2)} &
-i \frac {\lambda\mu_1\mu_2}{(\lambda_1+2\mu_1) (\lambda_2+2\mu_2)} &
-4\frac{\lambda\mu_1\partial_z\mu_2}{(\lambda_1+2\mu_1)(\lambda_2+2\mu_2)}\\
-i \frac {\lambda\mu_1\mu_2}{(\lambda_1+2\mu_1) (\lambda_2+2\mu_2)} &
\frac {\lambda\mu_1\mu_2}{(\lambda_2+2\mu_2) (\lambda_1+2\mu_1)}&
4i\frac{\lambda\mu_1\partial_{z}\mu_2 }{(\lambda_1+2\mu_1)(\lambda_2+2\mu_2)}\\
-4\frac{\partial_z\mu_1\lambda\mu_2}{(\lambda_1+2\mu_1) (\lambda_2+2\mu_2)}&
4i\frac{\partial_{z}\mu_1\lambda\mu_2}{(\lambda_1+2\mu_1) (\lambda_2+2\mu_2)}&
8\frac{\lambda\partial_{ z}\mu_1\partial_{ z}\mu_2 }
{(\lambda_1+2\mu_1)(\lambda_2+2\mu_2)}-8\partial^2_{zz}\mu\end{matrix}\right)
+l.o.t.
\end{equation}
and
\begin{equation}\label{!PP}
{\bf B}_2=\left (\begin{matrix}\partial_{x_1}\left(\frac {\lambda\mu_1\mu_2}
{(\lambda_1+2\mu_1) (\lambda_2+2\mu_2)}\right)&-i\partial_{x_1}
\left(\frac {\lambda\mu_1\mu_2}{(\lambda_1+2\mu_1) (\lambda_2+2\mu_2)}\right)&
-4\partial_{x_1}\left(\frac {\lambda\mu_1\partial_z\mu_2}{(\lambda_1+2\mu_1)
(\lambda_2+2\mu_2)}\right)\\
\partial_{x_2}\left(\frac {\lambda\mu_1\mu_2}{(\lambda_1+2\mu_1)
(\lambda_2+2\mu_2)}\right)&-i\partial_{x_2}\left(\frac {\lambda\mu_1\mu_2}
{(\lambda_1+2\mu_1) (\lambda_2+2\mu_2)}\right)&-4\partial_{x_2}
\left(\frac {\lambda\mu_1\partial_z\mu_2}{(\lambda_1+2\mu_1)
(\lambda_2+2\mu_2)}\right)\\
-2\left(\frac {(\nabla \lambda,\nabla\mu_1)\mu_2}{(\lambda_1+2\mu_1)
(\lambda_2+2\mu_2)}\right)&2i\left(\frac {(\nabla \lambda,\nabla\mu_1)\mu_2}
{(\lambda_1+2\mu_1) (\lambda_2+2\mu_2)}\right)&8\left(\frac
 {(\nabla \lambda,\nabla\mu_1)\partial_z\mu_2}{(\lambda_1+2\mu_1)
(\lambda_2+2\mu_2)}\right )-16\partial^3_{zz\bar z}\mu\end{matrix}\right)+l.o.t.
\end{equation}
Here by $l.o.t.$ in (\ref{!!PP}) we mean a smooth matrix which is independent
of $\lambda$ and possibly dependent on derivatives of $\mu$ up to the
order two, and in (\ref{!PP}) we mean a smooth matrix which is possibly
dependent of $\lambda$ and possibly dependent on derivatives of $\mu$ up
to the order three.
From (\ref{!!PP}) and (\ref{!PP}), noting
$(\nabla\lambda,\nabla\mu_1)
= 2(\partial_z\lambda)(\partial_{\overline z}\mu_1)
+ 2(\partial_{\overline z}\lambda)(\partial_z\mu_1)$, we have
\begin{equation}\label{PP}
\partial^2_{\overline z \overline{z}}{\mathcal C}_2
- \partial_{\overline z}{\bf B}_2
= \left(\begin{matrix}\left(\frac {-(\partial^2_{\bar zz}\lambda)\mu_1\mu_2}
{(\lambda_1+2\mu_1) (\lambda_2+2\mu_2)}\right)&
i\left(\frac {(\partial^2_{\bar zz}\lambda)\mu_1\mu_2}{(\lambda_1+2\mu_1)
(\lambda_2+2\mu_2)}\right)&
4\left(\frac {(\partial^2_{z\bar z}\lambda)\mu_1\partial_z\mu_2}
{(\lambda_1+2\mu_1) (\lambda_2+2\mu_2)}\right)\\
i\left(\frac {-(\partial^2_{\bar zz}\lambda)\mu_1\mu_2}
{(\lambda_1+2\mu_1) (\lambda_2+2\mu_2)}\right)&
\left(\frac {(\partial^2_{\bar zz}\lambda)\mu_1\mu_2}{(\lambda_1+2\mu_1)
(\lambda_2+2\mu_2)}\right)&
4i\left(\frac {(\partial^2_{z\bar z}\lambda)\mu_1\partial_z\mu_2}
{(\lambda_1+2\mu_1) (\lambda_2+2\mu_2)}\right)\\
4\left(\frac {(\partial^2_{z\bar z} \lambda)(\partial_{\bar z}\mu_1)\mu_2}
{(\lambda_1+2\mu_1) (\lambda_2+2\mu_2)}\right)&
-4i\left(\frac {(\partial^2_{z\bar z} \lambda)(\partial_{\bar z}\mu_1)\mu_2}
{(\lambda_1+2\mu_1) (\lambda_2+2\mu_2)}\right)&
-4\left(\frac {(\partial^2_{z\bar z} \lambda)
(\partial_{\bar z}\mu_1)\partial_z\mu_2}{(\lambda_1+2\mu_1) (\lambda_2+2\mu_2)}\right )
+ 8\partial^4_{zz\bar{z}\bar{z}}\mu\end{matrix}\right)+l.o.t.
\end{equation}
Using (\ref{PP}) we write (\ref{pobedaA}) as
\begin{eqnarray}
 -\left\{\Delta\lambda\left (\frac{\mu_1(U_{0,1}+iU_{0,2})}
{\lambda_1+2\mu_1}
- 4\frac{(\partial_{\overline z}\mu_1)U_{0,3}}{\lambda_1+2\mu_1}\right)
\left(\frac{\mu_2(V_{0,1}-iV_{0,2})}{\lambda_2+2\mu_2}
- 4\frac{(\partial_{ z}\mu_2)V_{0,3}}{\lambda_2+2\mu_2}\right)
\right\}                            \nonumber\\
+ 8\Delta^2\mu V_{0,3}U_{0,3}  +P_3(x,D)\mu+P_1(x,D)\lambda
+ \frak Q_{V_0,U_0}(x)=0,
\end{eqnarray}
where $P_1,P_2$ are differential operators with smooth coefficients depending
on $U_0,V_0$ of orders three and one respectively.
For fixed functions $U,V$, the function $\frak H_{V,U}$ can be considered
as a fourth-order operator applied to the function $\mu$ and a 
second-order operator applied to the function $\lambda$ respectively.
By Proposition \ref{nikita} for each point $\widetilde x$ we can choose
functions $U_{0,\widetilde x}=(U_{0,1,\widetilde x}, U_{0,2,\widetilde x},
U_{0,3,\widetilde x})$ and $V_{0,\widetilde x}=(V_{0,1,\widetilde x},
V_{0,2,\widetilde x},V_{0,3,\widetilde x})$ satisfying (\ref{-5})
and (\ref{-55}) respectively and
$V_{0,3,\widetilde x}(\widetilde x)=U_{0,3,\widetilde x}(\widetilde x)=1.$
Applying again Proposition \ref{nikita} for each point $\widetilde x$
we can choose the function $\widetilde U_{0,\widetilde x}
=(\widetilde U_{0,1,\widetilde x}, \widetilde U_{0,2,\widetilde x},
\widetilde U_{0,3,\widetilde x})$ and $\widetilde V_{0,\widetilde x}
=(\widetilde V_{0,1,\widetilde x},\widetilde  V_{0,2,\widetilde x},
\widetilde V_{0,3,\widetilde x})$ satisfying (\ref{-5}) and (\ref{-55})
respectively and
$\widetilde V_{0,1,\widetilde x}(\widetilde x)=\widetilde U_{0,1,\widetilde x}
(\widetilde x)=1$ and $\widetilde V_{0,j,\widetilde x}(\widetilde x)
=\widetilde U_{0,j,\widetilde x}(\widetilde x)=0$ for $j=1,2.$

Therefore for any $\epsilon>0$ and any $\widetilde x$ there exists
$\delta (\widetilde x)>0$ such that
\begin{eqnarray}\label{pobeda1A}
V_{0,3,\widetilde x}(x)U_{0,3,\widetilde x}(x)\ge \frac 12
\quad \forall x\in B(\widetilde x,\delta(\widetilde x)),\\
\vert \widetilde V_{0,1,\widetilde x}( x)-1\vert
+ \vert\widetilde U_{0,1,\widetilde x}(x)-1\vert
+ \sum_{j=1}^2\vert \widetilde V_{0,1,\widetilde x}( x)\vert
+\vert\widetilde U_{0,1,\widetilde x}(x)\vert
\le \epsilon \quad \forall x\in B(\widetilde x,\delta(\widetilde x)).\nonumber
\end{eqnarray}

Then there exist partial differential operators $\widetilde P_3(x,D)$
and $\widetilde P_1(x,D)$ of orders three and one respectively such that
\begin{equation}\label{ix}
{\bf A}(x)\left (\begin{matrix} \Delta \lambda\\
\Delta ^2\mu\end{matrix}\right )+\widetilde P_3(x,D) \mu+\widetilde P_1(x,D)
\lambda=\left (
\begin{matrix}
\frak Q_{\widetilde V_{0,\widetilde x},\widetilde U_{0,\widetilde x},q_1,q_2}\\
\frak Q_{V_{0,\widetilde x},U_{0,\widetilde x},q_1,q_2}
\end{matrix}\right ) \quad \forall x\in B(\widetilde x,\delta(\widetilde x)).
\end{equation}

Taking the parameter $\epsilon$ sufficiently small we obtain that
$\det {\bf A}(x) \ne 0$ on $B(\widetilde x,\delta(\widetilde x)).$
Since $\overline \Omega$ is covered by $\cup_{\widetilde x\in
\overline\Omega}B(\widetilde x,\delta(\widetilde x))$ from such a covering,
one can take a finite subcovering
$B(\widetilde x_j,\delta(\widetilde x_j))$, $j\in \{1,\dots, J\}.$

From (\ref{ix}) there exist functions $c_{\beta,j}, b_{\beta,j},
p_{k,j,\ell}\in L^\infty(\Omega)$ such that
\begin{equation}\label{lazy5A}
\Delta^2 \mu + \sum_{\vert \beta\vert\le 3}c_{\beta,1}\partial^\beta_x \mu
+ \sum_{\vert \beta\vert\le 1}b_{\beta,1}\partial^\beta_x \lambda
+ \sum_{k=1}^J( p_{k,1,1}
\frak Q_{V_{0,\widetilde x_k},U_{0,\widetilde x_k},q_1,q_2}
+ p_{k,1,2} \frak Q_{\widetilde V_{0,\widetilde x_k},
\widetilde U_{0,\widetilde x_k},\widetilde q_1, \widetilde q_2})=0\quad \mbox{in}\,\,\Omega ,
\end{equation}
\begin{equation}\label{lazy55A}
\Delta \lambda+\sum_{\vert \beta\vert\le 3}
c_{\beta,2}\partial^\beta_x\mu
+ \sum_{\vert \beta\vert\le 1} b_{\beta,2}\partial^\beta_x \lambda
+ \sum_{k=1}^J( p_{k,2,1} \frak Q_{V_{0,\widetilde x_k},U_{0,\widetilde x_k},q_1,q_2}
+ p_{k,2,2} \frak Q
_{\widetilde V_{0,\widetilde x_k},\widetilde U_{0, \widetilde x_k}, \widetilde q_1, \widetilde q_2})=0\quad
\mbox{in}\,\,\Omega .
\end{equation}
From (\ref{kaput3}), (\ref{4.177}), (\ref{4.178}), (\ref{!zanoza1!}), (\ref{!zanoza2!}) and Propositions \ref{inga} and \ref{zanoza!},
there exist constants $C $ and $s_0$ such that
\begin{eqnarray}\label{lazy2A}
\Vert s\phi_s^\frac 12\frak Q_{\widetilde V_{0,\widetilde x_k},\widetilde U_{0,\widetilde x_k},\widetilde q_1,\widetilde q_2}
e^{s\phi_s}\Vert_{L^2(\Omega)}+\Vert s\phi_s^\frac 12\frak Q_{ V_{0,\widetilde x_k}, U_{0,\widetilde x_k},q_1,q_2}
e^{s\phi_s}\Vert_{L^2(\Omega)}\nonumber\\
\le C\left(\left\Vert \sum_{\vert\beta\vert\le 3}\partial^\beta_x \mu
e^{s\phi_s}\right\Vert_{L^2(\Omega)}
+ \left\Vert\sum_{\vert\beta\vert\le 1}\partial^\beta_x \lambda e^{s\phi_s}
\right\Vert_{L^2(\Omega)}\right)\quad\forall s\ge s_0.
\end{eqnarray}
Applying the Carleman estimate (\ref{lazy1}) to the system (\ref{lazy5A}) and
(\ref{lazy55A}) and using (\ref{lazy2A}) and (\ref{boleslavX}),
we obtain $(\mu,\lambda)\equiv 0.$
Thus the proof of the theorem is complete. $\blacksquare$


\begin{thebibliography}{99} %

\bibitem{ANS} M. Akamatsu, G. Nakamura, and S. Steinberg,
\textit{Identification of the Lam\'e coefficients from boundary
observations}, Inverse Problems, {\bf 7} (1991), 335-354.

\bibitem{BH} N. Bleistein and R. A. Handelsman,
{\it Asymptotic Expansions of Integrals,} Dover Publications, New
York, 1986.

\bibitem{C}  A. P.\ Calder\'on, \textit{On an inverse boundary value
problem,} in \emph{Seminar on Numerical Analysis and its
Applications to Continuum Physics}, 65--73, Soc. Brasil. Mat., R\'io
de Janeiro, 1980.

\bibitem{Es0} G. Eskin, \textit{Global uniqueness in the inverse
scattering problem for the Schr\"odinger operator with external
Yang-Mills potentials}, Comm. Math. Phys., {\bf 222} (2001), 503-531.

\bibitem{Es1} G.\ Eskin and J.\ Ralston,
\textit{On the inverse boundary value problem for linear isotropic
elasticity,} Inverse Problems, {\bf 18} (2002), 907--921.

\bibitem{Es2} G.\ Eskin and J.\ Ralston,
\textit{On the inverse boundary value problem for linear isotropic
elasticity and  Cauchy-Riemann systems}, in ''Inverse Problem and
Spectral Theory", Contempt. Math. {\bf 348}, Amer. Math. Soc.
Providence, RI, 2004, pp. 53--69.

\bibitem{FI} A. Fursikov and O. Imanuvilov,
{\it Local exact controllability of the two-dimensional
Navier-Stokes equations,} Mat. Sb., {\bf 187}, (1996), 102-138.

\bibitem{Heck}  H. Heck, J.N. Wang, and X. Li, \textit{Identification
of viscosity in an incompressible fluid}, Indiana Univ. Math. J.,
{\bf 56} (2007), 2489-2510.

\bibitem{Her} L. H\"ormander,  \textit{The Analysis of Linear
Partial Differential Operators IV}, Springer-Verlag, Berlin, 1985.

\bibitem{IK} M. Ikehata, \textit{Inversion formulas for the linearized
problem for an inverse boundary value problem in
elastic prospection}, SIAM J. Appl. Math., {\bf 50} (1990), 1635-1644.

\bibitem{IUY} O. Imanuvilov, G. Uhlmann and M. Yamamoto,
\textit{The Calder\'on problem with partial data in two dimensions},
J. Amer. Math. Soc., {\bf 23} (2010), 655-691.

\bibitem{IUYlame} O. Imanuvilov, G. Uhlmann and M. Yamamoto,
\textit{On uniqueness of Lam\'e coefficients from partial Cauchy
data in three dimensions}, Inverse Problems, {\bf 28} (2012), 125002.

\bibitem{IYlame} O. Imanuvilov and M. Yamamoto, \textit{
On reconstruction of Lam\'e coefficients from partial Cauchy data},
J. Inverse ll-Posed Problems, {\bf 19} (2011), 881-891.

\bibitem{IY4} O. Imanuvilov and M. Yamamoto,
\textit{Inverse boundary value problem for linear Schr\"odinger
equation in two dimensions,} arXiv:1208.3775

\bibitem{IY6} O. Imanuvilov and M. Yamamoto, {\it Inverse problem by
Cauchy data on an arbitrary sub-boundary for systems of elliptic
equations,} Inverse Problems, {\bf 28} (2012),  095015.

\bibitem{IY5} O. Imanuvilov and M. Yamamoto,
\textit{Uniqueness for inverse boundary value problems by Dirichlet-to
-Neumann map on subboundaries}, to appear in
Milano J. Mathematics (2013).

\bibitem{KA} L. Kantorovich and G. Akilov, {\it Functional Analysis,} Second Edition, Pergamant Press, Oxford-Eluford, N.Y. 1982.

\bibitem{Li} X. Li and J.-N. Wang, \textit{Determination of viscosity
in the stationary Navier-Stokes equations}, J. Differential Equations,
{\bf 242} (2007), 24-39.

\bibitem{N} A.~Nachman, \textit{Global uniqueness
for a two-dimensional inverse boundary value problem}, Ann. of
Math., \textbf{143} (1996), 71--96.

\bibitem{NU1} G. Nakamura and G. Uhlmann, \textit{Identification of
Lam\'e parameters by boundary measurements}, American
Journal of Mathematics, {\bf 115} (1993), 1161-1187.

\bibitem{NU2} G. Nakamura and G. Uhlmann, \textit{
Global uniqueness for an inverse boundary value problem arising
in elasticity}, Invent. Math., {\bf 118} (1994), 457-474, ERRATUM,
Invent. Math., {\bf 152} (2003), 205-207.

\bibitem{NU3} G. Nakamura and G. Uhlmann, \textit{Inverse boundary problems
at the boundary for an elastic system}, SIAM J. Math. Anal., {\bf 26} (1995),
263-279.

\bibitem{SU} J.~Sylvester and G.~Uhlmann,
\textit{A global uniqueness theorem for an inverse boundary value
problem}, Ann. of Math., \textbf{125}  (1987), 153--169.

\bibitem{U} G. Uhlmann, {\it Electrical impedance tomography
and Calder\"on's problem},
Inverse Problems, {\bf 25} (2009) 123011 (39pp).

\bibitem{VE}  I.\ Vekua, \textit{Generalized Analytic Functions},
Pergamon Press, Oxford, 1962.

\end{thebibliography}
\end{document}